%% file: induction-on-descent.tex
\documentclass[11pt]{article}

\usepackage[text={390pt, 592pt}, a4paper]{geometry}
\usepackage{indentfirst}
\usepackage{hyperref}
\usepackage{breakurl}
\usepackage{amssymb, amsmath}
\usepackage{array}
\usepackage{graphicx}
\usepackage[labelsep=none]{caption}

\allowdisplaybreaks

\hypersetup{pdfstartview=FitH, pdfpagemode=UseNone}

\newcounter{lemma}
\newenvironment{lemma}{\refstepcounter{lemma} \emph{Lemma \thelemma.}}{}
\newenvironment{lemma*}[1]{\emph{Lemma #1.}}{}
\newcounter{proposition}
\newenvironment{proposition}{\refstepcounter{proposition} \emph{Proposition \theproposition.}}{}
\newenvironment{proposition*}[1]{\emph{Proposition #1.}}{}
\newcounter{theorem}
\newenvironment{theorem}{\refstepcounter{theorem} \emph{Theorem \thetheorem.}}{}
\newenvironment{theorem*}[1]{\emph{Theorem #1.}}{}
\newenvironment{proof}{\emph{Proof.}}{$\square$}
\newenvironment{proof*}[1]{\emph{#1.}}{$\square$}

\newcommand{\mmod}{\mathbin{\mathrm{mmod}}}

\DeclareMathOperator{\ecf}{ecf}

\DeclareMathOperator{\BBox}{Box}
\DeclareMathOperator{\Ball}{Ball}

\DeclareMathOperator{\Min}{Min}

\DeclareMathOperator{\Comp}{Comp}
\newcommand{\CompX}{\Comp_X}
\newcommand{\CompY}{\Comp_Y}

\newcommand{\Lat}{\mathrm{Lat}}
\newcommand{\Dia}{\mathrm{Dia}}
\newcommand{\Mar}{\mathrm{Mar}}

\newcommand{\finv}{f^\downarrow}

\newcommand{\ellinv}{\ell^\downarrow}

\begin{document}

\input{induction-on-descent-01-intro}

\input{induction-on-descent-02-prelim}

\input{induction-on-descent-03-mono}

\input{induction-on-descent-04-proj}

\input{induction-on-descent-05-weave}

\input{induction-on-descent-06-ang}

\input{induction-on-descent-07-desc}

\input{induction-on-descent-08-lift}

\input{induction-on-descent-09-low}

\input{induction-on-descent-10-line}

\input{induction-on-descent-11-conn}

\input{induction-on-descent-12-rigid}

\input{induction-on-descent-13-journey}

\input{induction-on-descent-14-further}

\appendix

\input{induction-on-descent-15-biconn}

\input{induction-on-descent-16-weave-ii}

\input{induction-on-descent-17-ang-ii}

\input{induction-on-descent-18-refs}
\end{document}

%% file: induction-on-descent-01-intro.tex
\title{\textbf{Induction on Descent in Leaper Graphs}}
\author{Nikolai Beluhov}
\date{}

\maketitle

\begin{center} \parbox{0.75\textwidth}{\setlength{\parindent}{1.2em} \footnotesize \emph{Abstract}. We construct an infinite ternary tree $\mathfrak{L}$ whose root is the knight and whose vertices are all skew free leapers. We define the \emph{descent} of a skew free leaper to be its ``address'' within $\mathfrak{L}$. We introduce three transformations which relate the leaper graphs of a skew free leaper to the leaper graphs of its three children in $\mathfrak{L}$. By starting with the knight and then applying these transformations so as to advance throughout $\mathfrak{L}$, we can establish theorems about all skew free leapers. We call this proof technique \emph{induction on descent} and with its help we resolve a number of questions about leaper graphs.} \end{center}

\section{Introduction} \label{intro}

\subsection{Boards and Leapers} \label{intro:bl}

A \emph{$(p, q)$-leaper} $L$ is a fairy chess piece which generalises the knight. Two cells $(x', y')$ and $(x'', y'')$ are adjacent with respect to $L$ when $\{|x' - x''|, |y' - y''|\} = \{p, q\}$. Suppose, for concreteness, that $p \le q$. Then the knight corresponds to $p = 1$ and $q = 2$. Other leapers with traditional names include the wazir ($p = 0$ and $q = 1$), fers ($p = q = 1$), camel ($p = 1$ and $q = 3$), giraffe ($p = 1$ and $q = 4$), and zebra ($p = 2$ and $q = 3$).

We define a \emph{leaper graph} of $L$ to be a graph all of whose vertices are cells and all of whose edges join pairs of cells adjacent with respect to $L$.

Let $A$ be a board of size $m \times n$. (Thus with $m$ rows and $n$ columns.) We define the \emph{complete leaper graph} of $L$ on $A$ to be the leaper graph of $L$ on $A$ whose vertices are all cells of $A$ and whose edges are all edges of $L$ on $A$.

Note that most literature on the subject defines a ``leaper graph'' to mean a ``complete leaper graph on a rectangular board'' and then considers leaper graphs of this kind only.

A leaper is \emph{skew} when $p$ and $q$ are distinct positive integers. \cite{J}

A leaper is \emph{free} when its complete leaper graph on the infinite board $\mathbb{Z} \times \mathbb{Z}$ is connected. Or, equivalently, when its complete leaper graph is connected on every board $A$ such that both of $m$ and $n$ are sufficiently large. \cite{J} It is straightforward to see that $L$ is free if and only if $p + q$ is odd and $p$ and $q$ are relatively prime.

Most of the problems we consider only make sense for skew free leapers. (The one exception is the topic of Section \ref{rigid}. There, we briefly sketch how our analysis generalises to skew but non-free leapers.) Thus suppose that $L$ is both skew and free.

For convenience, from this point on we fix the meanings of the symbols $p$, $q$, $L$, $m$, $n$, and $A$. Throughout this work, we use $p$ and $q$ to denote two positive integers with $p < q$, $p + q$ odd, and $\gcd(p, q) = 1$; we use $L$ to denote the skew free $(p, q)$-leaper; we use $m$ and $n$ to denote two positive integers; and we use $A$ to denote the board of size $m \times n$.

We go on to introduce three different topics in the study of leaper graphs. They will appear to be unrelated at first. Eventually, however, we will discover some surprising connections between them.

\subsection{The Second Leaper Theorem} \label{intro:sl}

Let $A_\square$ be the square board of side $p + q$ and let $G_\square$ be the complete leaper graph of $L$ on $A_\square$. Then each vertex of $G_\square$ is of degree either zero or two. Consequently, $G_\square$ consists of a number of isolated vertices and a number of pairwise disjoint cycles.

The main result of the author's earlier work \cite{Be} is as follows.

\medskip

\begin{theorem*}{\textbf{SL}} (Second leaper.) Let $L$ be a skew free $(p, q)$-leaper, let $C$ be a cycle of $L$ on the square board of side $p + q$, and let $S$ be the vertex set of $C$. Then there exists a second leaper $M$, distinct from $L$, such that $S$ admits a Hamiltonian cycle $D$ of $M$. \end{theorem*}

\medskip

We call a cycle of $L$ as in the theorem a \emph{clover}.

\begin{figure}[p] \centering \begin{tabular}{>{\centering}m{0.45\textwidth} >{\centering}m{0.45\textwidth}} \includegraphics{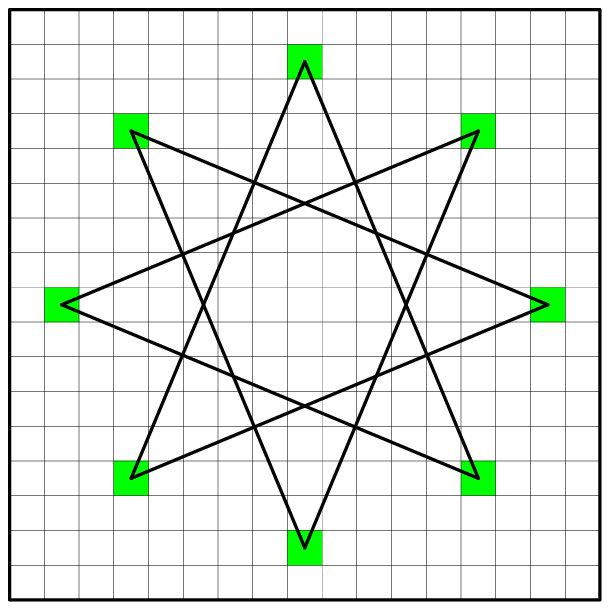} & \includegraphics{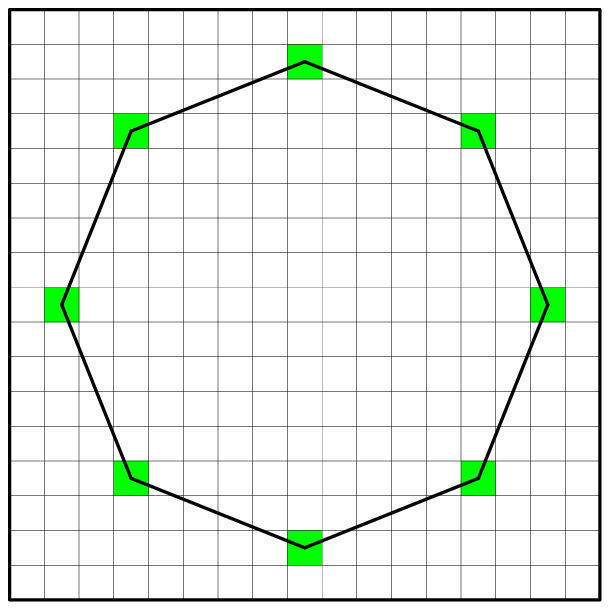} \tabularnewline \smallbreak \tabularnewline \includegraphics{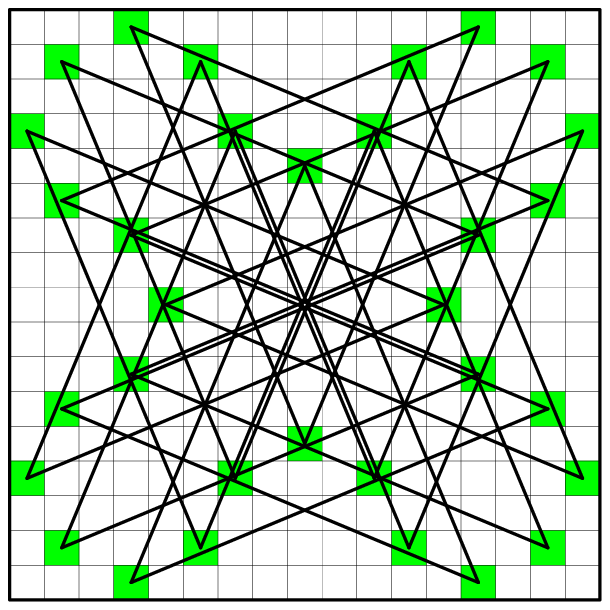} & \includegraphics{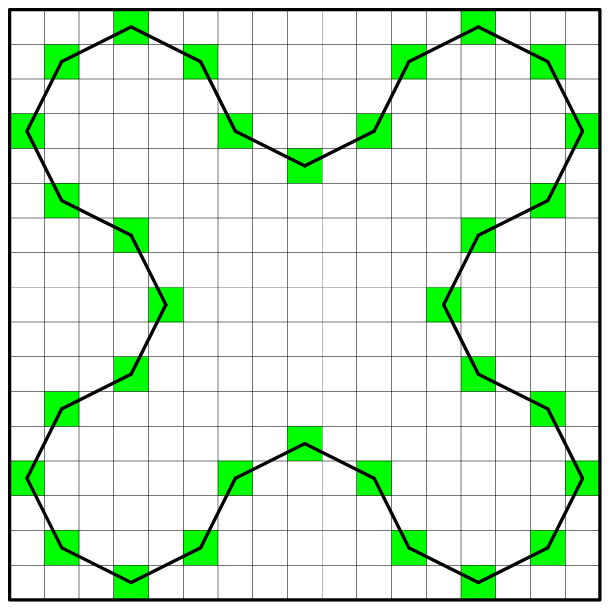} \tabularnewline \smallbreak \tabularnewline \includegraphics{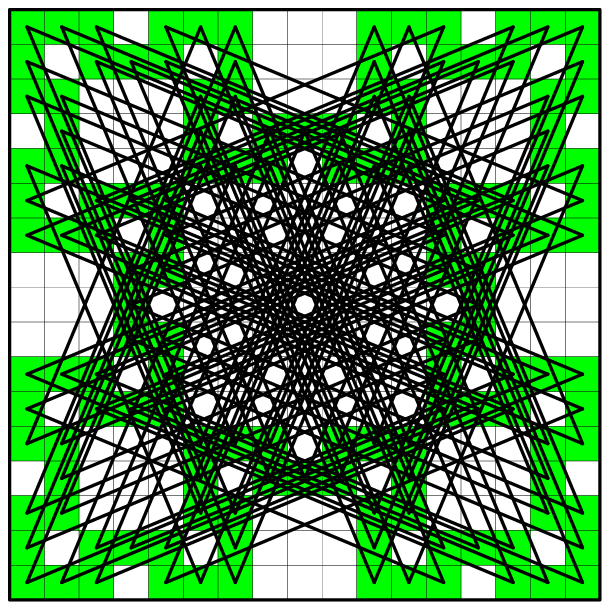} & \includegraphics{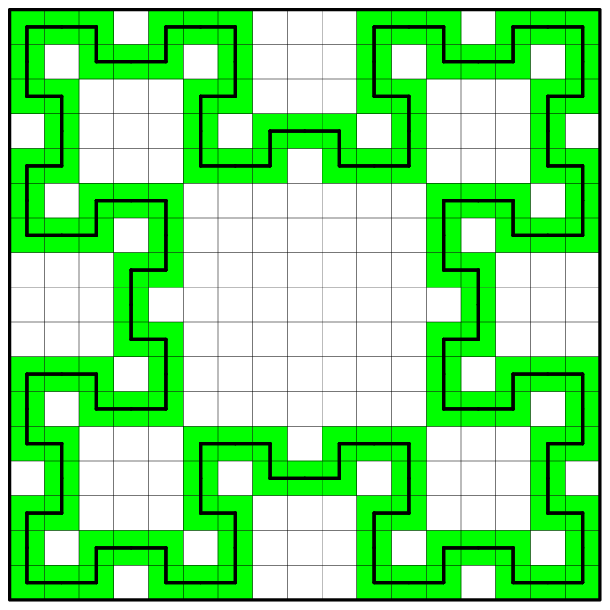} \end{tabular} \caption{} \label{sl} \end{figure}

For example, Figure \ref{sl} shows three clovers $C$ of the $(5, 12)$-leaper on the board of size $17 \times 17$ together with the corresponding cycles $D$ of the wazir, knight, and $(2, 5)$-leaper.

The proof of the second leaper theorem in \cite{Be} reveals some deeper structure, as follows.

The rational number $q/p$ admits a unique expansion of the form \[c_\kappa \pm \cfrac{1}{c_{\kappa - 1} \pm \cfrac{1}{\ddots \displaystyle \genfrac{}{}{0pt}{}{{}}{c_1}}}\] such that all of the $c_i$ are even positive integers. We call this the \emph{even continued fraction} of $L$ and we denote it by $\ecf(L)$. We also call $\kappa$ the \emph{depth} of $L$. For convenience, we abbreviate the expression $c_\kappa + \varepsilon_{\kappa - 1}/(c_{\kappa - 1} + \varepsilon_{\kappa - 2}/(\ldots c_1 \ldots))$ as $[c_\kappa, \varepsilon_{\kappa - 1}, c_{\kappa - 1}, \varepsilon_{\kappa - 2}, \ldots, c_1]$, where $\varepsilon_i \in \{-1, 1\}$ for all $i$ with $1 \le i < \kappa$.

We define the $i$-th \emph{tail} of $L$ to be the skew free $(p_i, q_i)$-leaper $L_i$ determined by $\ecf(L_i) = [c_i, \varepsilon_{i - 1}, c_{i - 1}, \varepsilon_{i - 2}, \ldots, c_1]$, for all $i$ with $1 \le i \le \kappa$. We also set $p_0 = 0$ and $q_0 = 1$ and we define the zeroth tail $L_0$ of $L$ to be the wazir.

Then we can group the clovers of $G_\square$ into $\kappa$ classes so that: (a) Two clovers are translation copies of one another if and only if they are in the same class; (b) Class $i$ contains a total of $(q_i - p_i)^2$ clovers; and (c) The vertex set $S$ of each clover $C$ in class $i$ admits also a Hamiltonian cycle $D$ of $L_i$, for all $i$ with $0 \le i < \kappa$.

\subsection{Directional Rigidity} \label{intro:rigid}

A \emph{move} of $L$ is an oriented edge of $L$. We define the \emph{direction} of a move $a \to b$ of $L$ as follows: When $b = a + (q, p)$, we say that the move points east-northeast; when $b = a + (p, q)$, north-northeast; \ldots; and, when $b = a + (q, -p)$, east-southeast.

We define the directions of the moves of a skew but not necessarily free leaper similarly.

Let $G'$ and $G''$ be leaper graphs of the skew leapers $L'$ and $L''$, respectively. Suppose that there exists a direction-preserving isomorphism $\varphi$ between $G'$ and $G''$, so that for each edge $ab$ of $G'$ the move $a \to b$ of $L'$ in $G'$ and its corresponding move $\varphi(a) \to \varphi(b)$ of $L''$ in $G''$ point in the same direction. Then we say that $G'$ and $G''$ are \emph{direction-isomorphic}.

\begin{figure}[ht!] \centering \begin{tabular}{>{\centering}m{0.4\textwidth} >{\centering}m{0.4\textwidth}} \includegraphics{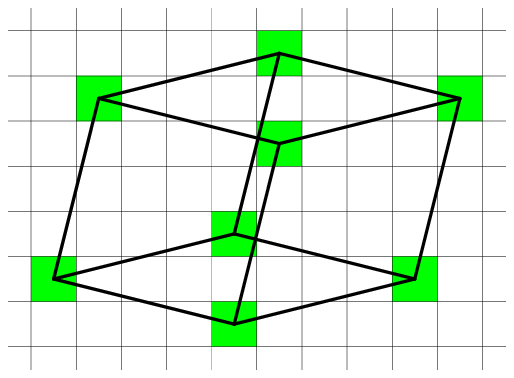} & \includegraphics{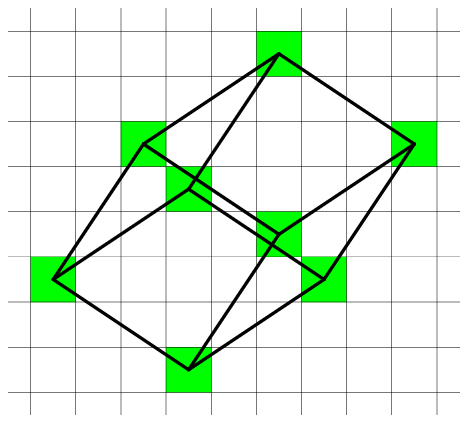} \end{tabular} \caption{} \label{cube} \end{figure}

For example, the two leaper graphs of the giraffe and the zebra in Figure \ref{cube} are direction-isomorphic.

We say that two skew leapers are \emph{proportional} when their parameters are in the same ratio. Clearly, if $L'$ and $L''$ are proportional, then every finite leaper graph of $L'$ is direction-isomorphic to some leaper graph of $L''$, and vice versa.

We say that a finite leaper graph $G$ of $L$ is \emph{directionally flexible} when it is direction-isomorphic to some leaper graph of a skew leaper which is not proportional to $L$. Otherwise, when $G$ is only direction-isomorphic to leaper graphs of skew leapers proportional to $L$, we say that it is \emph{directionally rigid}.

Essentially the same system of concepts is developed in \cite{Be} from the point of view of the equivalence classes that leaper graphs form with respect to direction-preserving isomorphisms.

One of the most natural questions raised by the notion of directional flexibility and rigidity is this: What are the boards $A$ such that the complete leaper graph of $L$ on $A$ is directionally rigid?

Let $A'$ and $A''$ be two boards of sizes $m' \times n'$ and $m'' \times n''$, respectively. We say that $A'$ \emph{fits} inside of $A''$, and we write $A' \sqsubseteq A''$, when either $m' \le m''$ and $n' \le n''$ or, conversely, $n' \le m''$ and $m' \le n''$.

Observe that if $A' \sqsubseteq A''$ and the complete leaper graph of $L$ on $A'$ is directionally rigid, then the complete leaper graph of $L$ on $A''$ is directionally rigid as well. Thus there exists a finite set of boards $\mathcal{A}$ such that the complete leaper graph of $L$ on $A$ is directionally rigid if and only if an element of $\mathcal{A}$ fits inside of $A$. An explicit set $\mathcal{A}$ of this kind would give us the characterisation we seek.

One of our main results is as follows.

\medskip

\begin{theorem} \label{thm:rigid} Let $L$ be a skew free $(p, q)$-leaper of depth $\kappa$ with tails the $(p_i, q_i)$-leapers $L_i$ and let $A^\mathcal{R}_i$ be the board of size $(p + q - p_i + 1) \times (p + q + q_i - [q_i \bmod 2])$ for all $i$ with $0 \le i \le \kappa$. Then the complete leaper graph of $L$ on $A$ is directionally rigid if and only if $A^\mathcal{R}_i \sqsubseteq A$ for some $i$. \end{theorem}

\medskip

Thus the answer to our question hinges on the even continued fraction of $L$. We give a proof of Theorem \ref{thm:rigid} in Section \ref{rigid}.

\subsection{Wazir Journeys} \label{intro:journey}

A \emph{wazir journey} of a leaper $M$ is a walk of $M$ whose endpoints are adjacent with respect to the wazir. (Or, equivalently, adjacent by side.) This is also known as a \emph{$(0, 1)$-journey} of $M$. \cite{J} Observe that $M$ is free if and only if there exists a wazir journey of $M$.

\begin{figure}[ht!] \centering \includegraphics{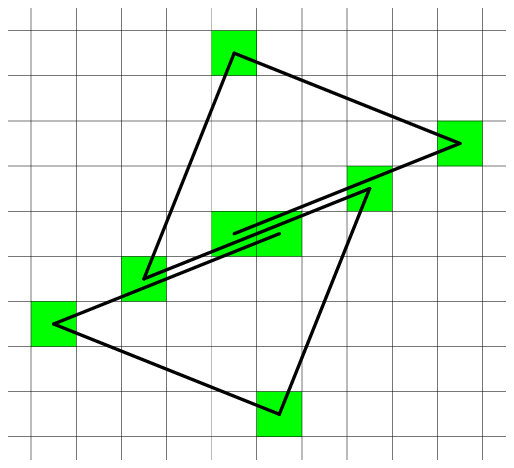} \caption{} \label{wj} \end{figure}

For example, Figure \ref{wj} shows one wazir journey of the $(2, 5)$-leaper.

Walks of this kind were introduced by Alasdair Houston, who in the 1970s posed to the members of the Fairy Chess Correspondence Circle the problem of finding the smallest number of moves in a wazir journey of a given skew free leaper. Houston solved the special case when $q = p + 1$ and George Jelliss solved the special cases when $p = 1$ and $p = 2$. \cite{J} The general form of the problem was settled by Krit Boonsiriseth in \cite{Bo}.

A different natural question suggested by the notion of a wazir journey is this: What are the boards $A$ such that there exists a wazir journey of $L$ on $A$?

As before, if $A' \sqsubseteq A''$ and there exists a wazir journey of $L$ on $A'$, then there exists a wazir journey of $L$ on $A''$ as well. Thus once again we should expect a characterisation in terms of some concrete finite set of boards an element of which must fit inside of $A$.

One of our main results is as follows.

\medskip

\begin{theorem} \label{thm:journey} Let $L$ be a skew free $(p, q)$-leaper with even continued fraction $[c_\kappa, \varepsilon_{\kappa - 1}, c_{\kappa - 1}, \varepsilon_{\kappa - 2}, \ldots, c_1]$ and tails the $(p_i, q_i)$-leapers $L_i$. (a) Let $A^\mathcal{W}_0$ be the board of size $(p + q) \times (p + q)$. (b) When $\kappa \ge 2$, let $A^\mathcal{W}_1$ be the board of size $(p + q - c_1 + 1) \times (p + q + q_2 - 1)$. (c) When $\kappa \ge 3$ and $\varepsilon_1 = 1$, let $A^\mathcal{W}_2$ be the board of size $(p + q - c_1) \times (p + q + q_3 - c_1 - 1)$. Then there exists a wazir journey of $L$ on $A$ if and only if $A^\mathcal{W}_i \sqsubseteq A$ for some $i$. \end{theorem}

\medskip

Just as with Theorem \ref{thm:rigid}, the answer to our question involves the even continued fraction of $L$. We give a proof of Theorem \ref{thm:journey} in Section \ref{journey}.

\subsection{The Structure of This Work} \label{intro:struct}

Why do the even continued fractions of skew free leapers play a key role in all three of these settings?

This work is one attempt at an answer. Our goal in it will be to develop a general toolkit which allows us to approach the second leaper theorem as well as both Theorems \ref{thm:rigid} and \ref{thm:journey} in a unified manner.

The rest of the text is structured as follows.

Section \ref{prelim} covers some basic definitions and notations.

Section \ref{mono} generalises some observations we made in the introduction and surveys a number of additional topics to which we apply our methods later on.

Sections \ref{proj}--\ref{ang} introduce some crucial concepts.

Over the course of Sections \ref{desc}--\ref{line}, we develop our general toolkit.

In Section \ref{desc}, we construct an infinite ternary tree $\mathfrak{L}$ whose root is the knight and whose vertices are all skew free leapers. We also define the \emph{descent} of $L$ to be its ``address'' within $\mathfrak{L}$. The descent of $L$ and its even continued fraction encode essentially the same information about $L$, though in slightly different ways.

Then, in Section \ref{lift}, we introduce three transformations which relate the leaper graphs of $L$ to the leaper graphs of its three children in $\mathfrak{L}$. Conversely, in Section \ref{low} we introduce three complementary transformations which relate the leaper graphs of the children of $L$ in $\mathfrak{L}$ back to the leaper graphs of $L$ itself. Section \ref{line} studies systems of leaper graphs interrelated in these ways.

By starting with the knight and then applying the transformations of Sections \ref{lift} and \ref{low} so as to advance throughout $\mathfrak{L}$, we can establish theorems about all skew free leapers. We call this proof technique \emph{induction on descent}.

Sections \ref{conn}--\ref{journey} contain applications. In Section \ref{conn}, we return to the topics surveyed in Section \ref{mono} in order to discuss the connectedness of leaper graphs. Then in Section \ref{rigid} we prove Theorem \ref{thm:rigid} and in Section \ref{journey} we prove Theorem \ref{thm:journey}.

Finally, Section \ref{further} briefly touches upon some open problems where induction on descent would likely be helpful.

%% file: induction-on-descent-02-prelim.tex
\section{Preliminaries} \label{prelim}

We formalise a cell as an ordered pair of integers and a board as the Cartesian product of two integer intervals, so that a board is a set of cells.

Throughout this work, let $A = I_X \times I_Y$ with $I_X = [x_\text{Min}; x_\text{Max}]$ and $I_Y = [y_\text{Min}; y_\text{Max}]$. Thus $m = |I_Y|$ and $n = |I_X|$.

Our formalisation of boards is such that there exist infinitely many boards of each given size. On the other hand, in the introduction we often spoke of ``the'' board of some size, in settings where it did not matter which specific board of that size it is. Let us formalise that phrasing, too.

We say that $A$ is the \emph{standard} board of size $m \times n$ when $x_\text{Min} + x_\text{Max} \in \{0, 1\}$ and $y_\text{Min} + y_\text{Max} \in \{0, 1\}$. Thus all standard boards with even height and even width share the same center of symmetry, and similarly for the other three combinations of height parity and width parity. Whenever we refer to ``the'' board of some size, as if there existed only one such board, we mean the standard board of that size.

We define the \emph{parity} of a cell $(x, y)$ to be the parity of the sum $x + y$.

We say that two boards of sizes $m \times n$ and $n \times m$ are \emph{transposes} of one another.

Given a nonempty set of cells $S$, let $\min_X(S)$ be the least $x$-coordinate of a cell of $S$, and define $\min_Y(S)$, $\max_X(S)$, and $\max_Y(S)$ similarly. Then the \emph{bounding box} of $S$ is the board $\BBox(S) = [\min_X(S); \max_X(S)] \times [\min_Y(S); \max_Y(S)]$. Clearly, the bounding box of $S$ is the $\subseteq$-smallest board which includes $S$.

Given a cell $c = (x_\circ, y_\circ)$, we define the \emph{ball} of center $c$ and radius $r$, denoted $\Ball(c, r)$, to be the set of all cells $(x, y)$ with $|x_\circ - x| \le r$ and $|y_\circ - y| \le r$. Thus $\Ball(c, r)$ is a square board of side $2r + 1$ centered at $c$.

Let $s$ and $t$ be positive integers. We write $s \mmod t$ for the distance from $s$ to the nearest multiple of $t$. Or, equivalently, $s \mmod t$ denotes the smallest nonnegative integer $r$ such that $s \equiv \pm r \pmod t$.

Define $p_0$, $q_0$, $p_1$, $q_1$, \ldots, $p_\kappa$, $q_\kappa$ as in the introduction. Then $q \mmod 2p = p_{\kappa - 1}$ for all skew free leapers $L$. Therefore, the tails of $L$ are determined inductively by $p_{i - 1} = q_i \mmod 2p_i$ and $q_{i - 1} = p_i$ for all $i$ with $1 \le i \le \kappa$.

Consider a graph $G$. For the edge of $G$ joining vertices $u$ and $v$, we write either $uv$ or $u$---$v$, whichever one reads more clearly in the setting at hand. We measure the length of a path or walk in $G$ in terms of the number of edges it traverses. We define a connected component of $G$ to be a \emph{singleton} when it consists of a single isolated vertex.

Let $\nu = u_1u_2 \ldots u_k$ be a walk in $G$ and let $v$ be a vertex of $G$ adjacent to $u_i$ for some $i$. Then the insertion of a \emph{detour} at $u_i$ via $v$ transforms $\nu$ into the walk $u_1u_2 \ldots u_{i - 1}u_ivu_iu_{i + 1}u_{i + 2} \ldots u_k$ in $G$.

Let $S$ be a subset of the vertex set of $G$. Then $S$ \emph{induces} the subgraph of $G$ with vertex set $S$ whose edges are all edges of $G$ between the elements of $S$.

We define a leaper graph of $L$ to be \emph{induced} when its edge set consists of all edges of $L$ between the cells of its vertex set.

We define a \emph{connected component} of $L$ on $A$ to be a connected component in the complete leaper graph of $L$ on $A$.

Let $e = (x', y')$---$(x'', y'')$ be an edge of $L$. We define the \emph{incline} of $e$ as follows: We say that $e$ is \emph{slight} when $|x' - x''| = q$ and $|y' - y''| = p$, and we say that $e$ is \emph{steep} otherwise, when $|x' - x''| = p$ and $|y' - y''| = q$.

Consider two edges $ab'$ and $ab''$ of $L$. Then angle $b'ab''$ is \emph{lateral} when $b'$ and $b''$ are distinct but in the same row or column, and \emph{diagonal} when $b'$ and $b''$ are distinct but on the same diagonal. Thus the seven types of angles that two edges of $L$ incident with the same cell can form are zero (when $b' = b''$), laterally acute, diagonally acute, right, laterally obtuse, diagonally obtuse, and straight (when $a$ is the midpoint of segment $b'b''$). \cite{J}

\begin{figure}[ht!] \centering \begin{tabular}{>{\centering}m{0.108\textwidth} >{\centering}m{0.108\textwidth} >{\centering}m{0.135\textwidth} >{\centering}m{0.162\textwidth} >{\centering}m{0.135\textwidth} >{\centering}m{0.162\textwidth}} \includegraphics{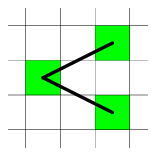} & \includegraphics{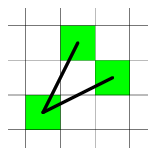} & \includegraphics{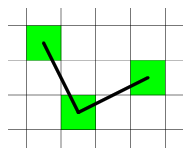} & \includegraphics{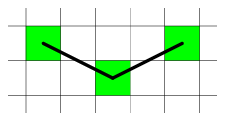} & \includegraphics{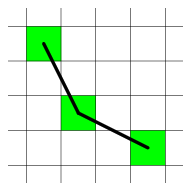} & \includegraphics{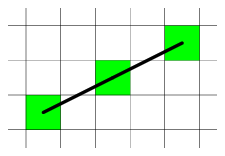} \end{tabular} \caption{} \label{angles} \end{figure}

For example, Figure \ref{angles} shows the six nonzero types of angles for the knight.

We assign a $2 \times 2$ matrix to each secondary intercardinal direction, as follows: \begin{align*} \mathcal{M}_\text{ENE} &= \left(\begin{smallmatrix} 0 & 1\\ 1 & 0 \end{smallmatrix}\right), & \mathcal{M}_\text{NNE} &= \left(\begin{smallmatrix} 1 & 0\\ 0 & 1 \end{smallmatrix}\right), & \mathcal{M}_\text{NNW} &= \left(\begin{smallmatrix} -1 & 0\\ 0 & 1 \end{smallmatrix}\right), & \mathcal{M}_\text{WNW} &= \left(\begin{smallmatrix} 0 & -1\\ 1 & 0 \end{smallmatrix}\right),\\ \mathcal{M}_\text{WSW} &= \left(\begin{smallmatrix} 0 & -1\\ -1 & 0 \end{smallmatrix}\right), & \mathcal{M}_\text{SSW} &= \left(\begin{smallmatrix} -1 & 0\\ 0 & -1 \end{smallmatrix}\right), & \mathcal{M}_\text{SSE} &= \left(\begin{smallmatrix} 1 & 0\\ 0 & -1 \end{smallmatrix}\right), & \mathcal{M}_\text{ESE} &= \left(\begin{smallmatrix} 0 & 1\\ -1 & 0 \end{smallmatrix}\right). \end{align*}

Consider a move $a' \to a''$ of $L$ pointing in direction $\delta$ with $a' = (x', y')$ and $a'' = (x'', y'')$. Then $\left(\begin{smallmatrix} x'' - x'\\ y'' - y' \end{smallmatrix}\right) = \mathcal{M}_\delta\left(\begin{smallmatrix} p\\ q \end{smallmatrix}\right)$, and we define $\Phi(a' \to a'')$ to be $\mathcal{M}_\delta$.

Consider also a walk $\alpha = a_1a_2 \ldots a_k$ of $L$ with $a_i = (x_i, y_i)$ for all $i$. Then we define $\Phi(\alpha)$ to be $\sum_{1 \le i < k} \Phi(a_i \to a_{i + 1})$, so that $\left(\begin{smallmatrix} x_k - x_1\\ y_k - y_1 \end{smallmatrix}\right) = \Phi(\alpha)\left(\begin{smallmatrix} p\\ q \end{smallmatrix}\right)$.

Let $\beta$ be a closed walk of $L$. We say that $\beta$ is \emph{balanced} when $\Phi(\beta) = \mathbf{0}$. Otherwise, we say that $\beta$ is \emph{unbalanced}.

\medskip

\begin{lemma} \label{ru} A finite leaper graph of $L$ is directionally rigid if and only if it contains an unbalanced closed walk. \end{lemma} 

\medskip

Observe that, in a finite leaper graph of $L$, all closed walks are balanced if and only if all cycles are balanced.

\medskip

\begin{proof} Let $G'$ be a leaper graph of the skew $(p', q')$-leaper $L'$ and let $G''$ be a leaper graph of the skew $(p'', q'')$-leaper $L''$ so that $L'$ and $L''$ are not proportional but $G'$ and $G''$ are direction-isomorphic.

Consider two corresponding closed walks $\alpha'$ in $G'$ and $\alpha''$ in $G''$. Since ${p' : q'} \neq {p'' : q''}$, $\Phi(\alpha') = \Phi(\alpha'')$, and $\Phi(\alpha')\left(\begin{smallmatrix} p'\\ q' \end{smallmatrix}\right) = \Phi(\alpha'')\left(\begin{smallmatrix} p''\\ q'' \end{smallmatrix}\right) = \mathbf{0}$, it follows that $\Phi(\alpha') = \Phi(\alpha'') = \mathbf{0}$.

Conversely, let $G$ be a finite leaper graph of $L$ such that every closed walk in $G$ is balanced. We claim that $G$ is then direction-isomorphic to leaper graphs of all but finitely many skew free leapers.

Clearly, it suffices to consider the case when $G$ is connected. Fix one cell $a$ of $G$. For each cell $b$ of $G$, let $\alpha_b$ be some walk from $a$ to $b$ in $G$ and define $\zeta(b) = \Phi(\alpha_b)$. Since all closed walks in $G$ are balanced, $\zeta(b)$ depends only on $b$ and not on our choice of a walk $\alpha_b$.

Observe that the $2 \times 2$ matrices $\zeta(b)$ are pairwise distinct when $b$ ranges over all cells of $G$. Therefore, for all but finitely many skew free $(r, s)$-leapers $M$ with $r < s$, the column-vectors $\zeta(b)\left(\begin{smallmatrix} r\\ s \end{smallmatrix}\right)$ will be pairwise distinct as well.

Consider one such skew free leaper $M$. We define $\varphi(b)$ to be the cell $(x, y)$ determined by $\left(\begin{smallmatrix} x\\ y \end{smallmatrix}\right) = \zeta(b)\left(\begin{smallmatrix} r\\ s \end{smallmatrix}\right)$. Then $\varphi$ is a direction-preserving isomorphism between $G$ and a leaper graph of $M$. \end{proof}

\medskip

The proof of Lemma \ref{ru} shows that if a finite leaper graph of $L$ is directionally flexible, then it is direction-isomorphic to leaper graphs of all but finitely many skew leapers with relatively prime parameters.

Clearly, the definitions of this section and the statement of Lemma \ref{ru} apply more generally to arbitrary skew but not necessarily free leapers.

We say that a leaper is \emph{half-free} when its complete leaper graph on the infinite board $\mathbb{Z} \times \mathbb{Z}$ consists of exactly two connected components. It is straightforward to see that this occurs if and only if the leaper's parameters are odd and relatively prime.

We also say that a leaper is \emph{relatively prime} when its parameters are relatively prime. Thus a leaper is relatively prime if and only if it is either free or half-free.

Let $M$ be an $(r, s)$-leaper with $d = \gcd(r, s)$. (Note that we consider only leapers with at least one nonzero parameter.) Then the $(r/d, s/d)$-leaper $M^\star$ is the unique relatively prime leaper proportional to $M$. Observe that every connected leaper graph of $M$ is a scaled copy, with a scaling factor of $d$, of some leaper graph of $M^\star$.

It is straightforward to see that the criteria which determine whether $M$ is free and whether it is half-free generalise as follows: The complete leaper graph of $M$ on the infinite board $\mathbb{Z} \times \mathbb{Z}$ consists of $d^2$ connected components when $(r + s)/d$ is odd and of $2d^2$ connected components otherwise, when $(r + s)/d$ is even.

%% file: induction-on-descent-03-mono.tex
\section{Monotonicity} \label{mono}

\subsection{Properties and Bases} \label{mono:pb}

Let $\mathcal{P}_L$ be a property involving $L$ which a board might or might not possess, and which depends only on the board's size. We write $\mathcal{P}_L(m, n)$ for ``$\mathcal{P}_L$ holds on the board of size $m \times n$'' and $\mathcal{P}_L(A)$ for ``$\mathcal{P}_L$ holds on $A$''.

Let $A'$ and $A''$ be two boards of sizes $m' \times n'$ and $m'' \times n''$, respectively. We say that $A'$ is \emph{smaller} than or \emph{congruent} to $A''$, and we write $A' \le A''$, when $m' \le m''$ and $n' \le n''$.

We say that $\mathcal{P}_L$ is \emph{monotone} when $\mathcal{P}_L(A')$ implies $\mathcal{P}_L(A'')$ for all boards $A'$ and $A''$ with $A' \le A''$.

Let $\mathcal{Q}_L$ be a monotone property. We say that $A$ is \emph{minimal} with respect to $\mathcal{Q}_L$ when $\mathcal{Q}_L$ holds on $A$ but not on any boards strictly smaller than $A$. We write $\Min(\mathcal{Q}_L)$ for the set of all standard boards minimal with respect to $\mathcal{Q}_L$, and we call this set the \emph{basis} of $\mathcal{Q}_L$.

Thus for all boards $A$ we get that $\mathcal{Q}_L(A)$ if and only if there exists a board $A^\star$ in $\Min(\mathcal{Q}_L)$ with $A^\star \le A$. Or, in other words, the basis of $\mathcal{Q}_L$ yields a complete description of all boards on which $\mathcal{Q}_L$ holds.

Observe that the basis of a monotone property is always finite. Furthermore, we can denote the boards of $\Min(\mathcal{Q}_L)$ by $A_1$, $A_2$, \ldots, $A_k$ so that, with $m_i \times n_i$ being the size of $A_i$ for all $i$, the heights and widths of these boards satisfy $m_1 < m_2 < \cdots < m_k$ and $n_1 > n_2 > \cdots > n_k$.

We say that $\mathcal{P}_L$ is \emph{symmetric} when $\mathcal{P}_L(m, n)$ if and only if $\mathcal{P}_L(n, m)$.

Suppose, from this point on, that $\mathcal{Q}_L$ is symmetric as well as monotone. Then $A_i$ and $A_{k + 1 - i}$ form a transpose pair for all $i$. In this setting, we define a \emph{reduced basis} of $\mathcal{Q}_L$ to be a subset of $\Min(\mathcal{Q}_L)$ which contains exactly one board out of each such pair. Thus every reduced basis of $\mathcal{Q}_L$ contains a total of $\lceil k/2 \rceil$ boards.

As before, a reduced basis $\mathcal{A}$ of $\mathcal{Q}_L$ yields a complete characterisation of all boards on which $\mathcal{Q}_L$ holds: For all boards $A$, we get that $\mathcal{Q}_L(A)$ if and only if there exists a board $A^\star$ in $\mathcal{A}$ with $A^\star \sqsubseteq A$.

Thus, in particular, we can also understand the motivating questions behind Theorems \ref{thm:rigid} and \ref{thm:journey} as ``What are the smallest boards inside of which we can fit a directionally rigid leaper graph of $L$?'' and ``What are the smallest boards inside of which we can fit a wazir journey of $L$?''. (Here, the formal meaning of ``smallest'' is ``$\sqsubseteq$-minimal''.)

Let us write $\mathcal{W}_L(A)$ for ``there exists a wazir journey of $L$ on $A$''.

Clearly, $\mathcal{W}_L$ is symmetric, and we already observed in the introduction that it is monotone as well. Theorem \ref{thm:journey} gives us a reduced basis of $\mathcal{W}_L$. From it, we see also that for all skew free leapers $L$ a reduced basis of $\mathcal{W}_L$ contains at most three boards.

Let us write $\mathcal{R}_L(A)$ for ``the complete leaper graph of $L$ on $A$ is directionally rigid''.

Clearly, $\mathcal{R}_L$ is symmetric, and we already observed in the introduction that it is monotone as well. Theorem \ref{thm:rigid} almost gives us a reduced basis of $\mathcal{R}_L$. To obtain one, from the set $\{A^\mathcal{R}_0, A^\mathcal{R}_1, \ldots, A^\mathcal{R}_\kappa\}$ we must remove all boards $A^\mathcal{R}_i$ such that $A^\mathcal{R}_j \sqsubseteq A^\mathcal{R}_i$ for some $j$ with $j \neq i$. That occurs if and only if either (a) $i = 1$ and $j = 0$; or (b) $i$ is odd, $j = i + 1$, $\varepsilon_1 = \varepsilon_2 = \cdots = \varepsilon_i = -1$, and $c_1 = c_2 = \cdots = c_j = 2$.

When $L$ does not satisfy the conditions of exception (b) with $i \ge 3$ and $j \ge 4$, the number of boards in a reduced basis of $\mathcal{R}_L$ is exactly equal to the depth of $L$. Thus, in particular, the size of a reduced basis of $\mathcal{R}_L$ can become arbitrarily large as $L$ ranges over all skew free leapers.

The rest of this section surveys some monotone properties that we study more carefully later on. Other monotone properties which might be amenable to our methods are described in Section \ref{further}. Lastly, in Appendix \ref{biconn} we consider one more monotone property which is only tangentially related to the main thread of our discussion.

\subsection{Isolated Vertices} \label{mono:i}

We write $\mathcal{I}_L(A)$ for ``the complete leaper graph of $L$ on $A$ does not contain isolated vertices''. This property is easily seen to be symmetric and monotone.

The next characterisation was given by Jelliss in \cite{J}.

\medskip

\begin{theorem} \label{thm:i} Let $L$ be a skew free $(p, q)$-leaper and let $A^\mathcal{I}$ be the board of size $2p \times 2q$. Then every cell of $A$ is incident with an edge of $L$ on $A$ if and only if $A^\mathcal{I} \sqsubseteq A$. \end{theorem} 

\medskip

Or, equivalently, $\{A^\mathcal{I}\}$ is a reduced basis of $\mathcal{I}_L$. The proof is straightforward.

\subsection{Connectedness} \label{mono:conn}

We write $\mathcal{C}_L(A)$ for ``the complete leaper graph of $L$ on $A$ is nontrivially connected''. (Here, ``nontrivially'' means that the graph contains at least two vertices. This rules out the board of size $1 \times 1$.)

Clearly, $\mathcal{C}_L$ is symmetric. Let us verify that it is monotone as well.

\medskip

\begin{proof*}{Proof of monotonicity for $\mathcal{C}_L$} Suppose that $\mathcal{C}_L(m, n)$. It suffices to show that then $\mathcal{C}_L(m, n + 1)$ as well.

Let $A^\star$ be a board of size $m \times (n + 1)$ and let $A'$ and $A''$ be its left and right subboards of size $m \times n$. Then the complete leaper graph of $L$ on each one of $A'$ and $A''$ is connected since $\mathcal{C}_L(m, n)$. Furthermore, $A'$ and $A''$ have at least one cell in common because $\mathcal{C}_L(m, n)$ implies $n \ge 2$. Therefore, the complete leaper graph of $L$ on $A^\star$ is connected as well. \end{proof*}

\medskip

The following characterisation was established independently by Frank Rhodes and Stephen Wilson in \cite{RW} as well as by Donald Knuth in \cite{K}.

\medskip

\begin{theorem} \label{thm:conn} Let $L$ be a skew free $(p, q)$-leaper with $p < q$ and let $A^\mathcal{C}$ be the board of size $(p + q) \times 2q$. Then the complete leaper graph of $L$ on $A$ is nontrivially connected if and only if $A^\mathcal{C} \sqsubseteq A$. \end{theorem} 

\medskip

Or, equivalently, $\{A^\mathcal{C}\}$ is a reduced basis of $\mathcal{C}_L$. We give a new proof of Theorem \ref{thm:conn} in Section \ref{conn}.

\subsection{Edge Connectedness} \label{mono:econn}

We write $\mathcal{E}_L(A)$ for ``the edges of $L$ on $A$ form a nonempty connected leaper graph''. Or, equivalently, ``all non-isolated vertices in the complete leaper graph of $L$ on $A$ belong to the same connected component of $L$ on $A$''.

Clearly, $\mathcal{E}_L$ is symmetric. Let us verify that it is monotone as well.

\medskip

\begin{proof*}{Proof of monotonicity for $\mathcal{E}_L$} Suppose that $\mathcal{E}_L(m, n)$. It suffices to show that then $\mathcal{E}_L(m, n + 1)$ as well.

Observe that $\mathcal{E}_L$ cannot hold on boards where all edges of $L$ are of the same incline. Indeed, let $e'$ be, for concreteness, a slight edge of $L$ and let $e''$ be the other slight edge of $L$ with the same bounding box as $e'$. Then no path of slight edges of $L$ can connect $e'$ and $e''$ even on the infinite board $\mathbb{Z} \times \mathbb{Z}$. Thus $\mathcal{E}_L(m, n)$ implies $m \ge q + 1$ and $n \ge q + 1$.

Let $A^\star$ be a board of size $m \times (n + 1)$ and let $A'$ and $A''$ be its left and right subboards of size $m \times n$. Since $m \ge q + 1$ and $n \ge q + 1$, we get that: (a) There is at least one edge of $L$ on $A' \cap A''$; and (b) Each edge of $L$ on $A^\star$ is contained entirely within one of $A'$ and $A''$. These observations together with $\mathcal{E}_L(A')$ and $\mathcal{E}_L(A'')$ imply that $\mathcal{E}_L(A^\star)$ as well. \end{proof*}

\medskip

We characterise the boards on which all edges of $L$ form a nonempty connected leaper graph as follows.

\medskip

\begin{theorem} \label{thm:econn} Let $L$ be a skew free $(p, q)$-leaper with $p < q$. (a) When $p \ge 2$, let $A^\mathcal{E}$ be the board of size $(p + q) \times (2p + q - [q \mmod 2p])$. (b) Otherwise, when $p = 1$, let $A^\mathcal{E}$ be the square board of side $q + 1$. Then all edges of $L$ on $A$ form a nonempty connected leaper graph if and only if $A^\mathcal{E} \sqsubseteq A$. \end{theorem} 

\medskip

Or, equivalently, $\{A^\mathcal{E}\}$ is a reduced basis of $\mathcal{E}_L$. We give a proof of Theorem \ref{thm:econn} in Section \ref{conn}.

%% file: induction-on-descent-04-proj.tex
\section{Projections} \label{proj}

Sometimes we can learn a lot about a leaper graph by focusing on one coordinate only and ignoring the other one entirely. We formalise this idea as follows.

Let $a$ and $b$ be distinct positive integers and let $I$ be a nonempty integer interval. We define the \emph{projection graph} $\Pi(a, b, I)$ to be the graph on vertex set $I$ where two vertices $u$ and $v$ are joined by an edge if and only if $|u - v| \in \{a, b\}$.

Suppose, for concreteness, that $a < b$. We call edge $uv$ \emph{short} when $|u - v| = a$ and \emph{long} otherwise, when $|u - v| = b$. Just as with leaper graphs, we call an oriented edge of $\Pi(a, b, I)$ a \emph{move}. We say that a move $u \to v$ points \emph{right} when $u < v$ and \emph{left} otherwise, when $u > v$.

Let $s = |I|$ be the size of $I$. When $s$ matters but the specific integer interval $I$ does not, we write simply $\Pi(a, b, s)$.

\begin{figure}[ht!] \centering \includegraphics{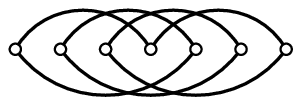} \caption{} \label{pg} \end{figure}

For example, Figure \ref{pg} shows $\Pi(3, 4, 7)$.

Let $d = \gcd(a, b)$, $a = da'$, and $b = db'$. Then each connected component of $\Pi(a, b, s)$ is a scaled copy, with a scaling factor of $d$, of either $\Pi(a', b', \lfloor s/d \rfloor)$ or $\Pi(a', b', \lceil s/d \rceil)$. Thus for most purposes it suffices to consider the case when $d = 1$ and $a$ and $b$ are relatively prime.

Throughout this work, we write $\Pi_X$ for $\Pi(p, q, I_X)$ and $\Pi_Y$ for $\Pi(p, q, I_Y)$.

Let $e = (x', y')$---$(x'', y'')$ be an edge of $L$ on $A$. Then we call edge $x'x''$ of $\Pi_X$ the \emph{$x$-projection} of $e$, and we define the \emph{$y$-projection} of $e$ similarly.

Let $\alpha = a_1a_2 \ldots a_k$ be a walk of $L$ on $A$ with $a_i = (x_i, y_i)$ for all $i$. Then we call walk $x_1x_2 \ldots x_k$ in $\Pi_X$ the \emph{$x$-projection} of $\alpha$. We define the \emph{$y$-projection} of $\alpha$ similarly. Note that in general the projections of a path of $L$ on $A$ will be walks but not necessarily paths and the projections of a cycle of $L$ on $A$ will be closed walks but not necessarily cycles.

Let $G$ be a leaper graph of $L$ on $A$. Then the \emph{$x$-projection} of $G$ is the subgraph of $\Pi_X$ whose vertex set consists of the $x$-coordinates of the vertices of $G$ and whose edge set consists of the $x$-projections of the edges of $G$. We define the \emph{$y$-projection} of $G$ similarly.

Philip Ginzboorg and Valtteri Niemi in \cite{GN} study the closely related family of graphs obtained when we orient all short edges one way and all long edges the opposite way. They also point out the connection with leaper graphs and establish a result equivalent to our Lemma \ref{pc} below.

The rest of this section collects a number of helpful lemmas about projections.

\medskip

\begin{lemma} \label{pc} Suppose that $a$ and $b$ are relatively prime. Then $\Pi(a, b, a + b)$ is a cycle. \end{lemma} 

\medskip

\begin{proof} Orient all short edges of $\Pi(a, b, a + b)$ to point to the right and all of its long edges to point to the left. Then each vertex attains unit in-degree as well as unit out-degree. Consequently, $\Pi(a, b, a + b)$ becomes the disjoint union of several oriented cycles.

Let $C$ be one of these cycles, containing $w$ short moves to the right and $z$ long moves to the left. Then $aw = bz$. Consequently, $a$ divides $z$ and $b$ divides $w$. Thus in fact $C$ contains all vertices and edges of $\Pi(a, b, a + b)$. \end{proof}

\medskip

One direct corollary of Lemma \ref{pc} is as follows.

\medskip

\begin{lemma} \label{pa} Suppose that $a$ and $b$ are relatively prime and let $C$ be a connected component of $\Pi(a, b, s)$. Then $C$ is acyclic if and only if $s \le a + b - 1$. Furthermore, $C$ is a path when $s \le a + b - 1$ and $C$ coincides with $\Pi(a, b, s)$ when $s \ge a + b$. \end{lemma} 

\medskip

Consider a walk $\nu = u_1u_2 \ldots u_k$ in $\Pi(a, b, I)$. We define the \emph{signature} of $\nu$ to be the word $\sigma = \sigma_1\sigma_2 \ldots \sigma_{k - 1}$ over the two-letter alphabet $\{\mathtt{s}, \mathtt{l}\}$ where, for all $i$, $\sigma_i = \mathtt{s}$ when edge $u_iu_{i + 1}$ of $\nu$ is short and $\sigma_i = \mathtt{l}$ otherwise, when edge $u_iu_{i + 1}$ of $\nu$ is long.

We say that two signatures $\sigma' = \sigma'_1\sigma'_2 \ldots \sigma'_{k'}$ and $\sigma'' = \sigma''_1\sigma''_2 \ldots \sigma''_{k''}$ are \emph{complements} when $k' = k''$ and $\sigma'_i \neq \sigma''_i$ for all $i$.

\medskip

\begin{lemma} \label{ps} Let $\chi$ be a walk in $\Pi_X$ and let $\upsilon$ be a walk in $\Pi_Y$. Then a walk of $L$ on $A$ whose $x$-projection is $\chi$ and whose $y$-projection is $\upsilon$ exists if and only if the signatures of $\chi$ and $\upsilon$ are complements. \end{lemma} 

\medskip

The proof is straightforward.

\medskip

\begin{lemma} \label{pf} Let $\chi = x_1x_2 \ldots x_k$ be a walk in $\Pi_X$ and let $C$ be a cycle in $\Pi_Y$ of the form $u$---$(u + p)$---$(u + p + q)$---$(u + q)$---$u$. Let $y_1$ be a vertex of $C$ and let $c = (x_1, y_1)$. Then there exists a unique walk of $L$ on $A$ which starts from $c$, whose $x$-projection is $\chi$, and whose $y$-projection is contained within $C$. \end{lemma} 

\medskip

\begin{proof} By Lemma \ref{ps}, it suffices to check that for each vertex $y_1$ of $C$ and every signature $\sigma$ there exists a unique walk in $C$ which starts from $y_1$ and whose signature is $\sigma$. This is clear since each vertex of $C$ is incident with exactly one short edge and exactly one long edge of $C$. \end{proof}

\medskip

Let $\sigma$ and $\sigma^\star$ be two signatures. We write $\sigma \overset{\iota}{\rightsquigarrow} \sigma^\star$ when we can convert $\sigma$ into $\sigma^\star$ by means of a series of transformations of the following forms: (a) Insertion of the subword $\mathtt{ss}$ at any position, possibly at the beginning or end; (b) When $\iota < 0$, substitution of any occurrence of the subword $\mathtt{l}$ with the subword $\mathtt{lll}$; and (c) Otherwise, when $\iota \ge 0$, insertion of the subword $\mathtt{ll}$ at any position, possibly at the beginning or end.

\medskip

\begin{lemma} \label{psc} Suppose that there exist a walk in $\Pi_X$ from $x'$ to $x''$ with nonempty signature $\sigma$ and a walk in $\Pi_Y$ from $y'$ to $y''$ with nonempty signature $\tau$. Suppose also that $\sigma \overset{n - 2q}{\rightsquigarrow} \sigma^\star$ and $\tau \overset{m - 2q}{\rightsquigarrow} \tau^\star$ so that signatures $\sigma^\star$ and $\tau^\star$ are complements. Then there exists a walk of $L$ on $A$ from cell $(x', y')$ to cell $(x'', y'')$. \end{lemma} 

\medskip

\begin{proof} Consider the original walk in $\Pi_X$ from $x'$ to $x''$. By inserting some detours into it, we can transform it into a walk in $\Pi_X$ from $x'$ to $x''$ with signature $\sigma^\star$. Similarly, we can also obtain a walk in $\Pi_Y$ from $y'$ to $y''$ with signature $\tau^\star$. The conclusion then follows by Lemma \ref{ps}. \end{proof}

\medskip

\begin{lemma} \label{pp} Let $G$ be a connected leaper graph of $L$ on $A$ such that the $y$-projection of $G$ is bipartite. Then, for each row of $A$, all cells of $G$ in that row are of the same parity. \end{lemma} 

\medskip

Observe that, in particular, Lemma \ref{pp} applies whenever the $y$-projection of $G$ is acyclic.

\medskip

\begin{proof} Let $a'$ and $a''$ be two cells of $G$ in the same row of $A$ and let $\alpha$ be a path in $G$ from $a'$ to $a''$. Since the $y$-coordinates of $a'$ and $a''$ coincide, the $y$-projection of $\alpha$ is a closed walk $\upsilon$ in the $y$-projection of $G$. Furthermore, since the $y$-projection of $G$ is bipartite, the length of $\upsilon$ is even. Thus $\alpha$ is of even length as well, and so $a'$ and $a''$ are of the same parity. \end{proof}

\medskip

Let $\nu$ be a closed walk in $\Pi(a, b, I)$. Then $\nu$ contains a different number of short moves to the left and right if and only if it contains a different number of long moves to the left and right. When both of these conditions are satisfied, we say that $\nu$ is \emph{unbalanced}.

Clearly, a closed walk of $L$ is unbalanced if and only if at least one of its projections is unbalanced.

\medskip

\begin{lemma} \label{pu} Let $d = \gcd(a, b)$. Then $\Pi(a, b, s)$ contains an unbalanced closed walk if and only if $s \ge a + b - d + 1$. \end{lemma} 

\medskip

\begin{proof} Let $a = da'$ and $b = db'$. Observe that $\Pi(a, b, s)$ contains an unbalanced closed walk if and only if $\Pi(a', b', \lceil s/d \rceil)$ does. When $s \le a + b - d$, by Lemma \ref{pa} we get that $\Pi(a', b', \lceil s/d \rceil)$ is acyclic, and so all closed walks in it are balanced. On the other hand, when $s = a + b - d + 1$, we get that $\lceil s/d \rceil = a' + b'$ and then $\Pi(a', b', a' + b')$ contains an unbalanced closed walk by the proof of Lemma \ref{pc}. \end{proof}

%% file: induction-on-descent-05-weave.tex
\section{Weaves} \label{weave}

We go on to define one important class of connected components in projection graphs. We consider the cases $p = 1$ and $p \ge 2$ separately.

Let $C$ be a connected component in $\Pi(p, q, s)$.

When $p = 1$, we define $C$ to be a \emph{weave} if and only if it is acyclic. Thus in this case $C$ is a weave if and only if $s \le q$ and $C$ coincides with $\Pi(1, q, s)$.

Otherwise, when $p \ge 2$, we define $C$ to be a \emph{weave} if and only if every path in $C$ traverses an even number of short edges between each pair of long ones. Or, equivalently, if and only if for every path in $C$ with signature of the form $\mathtt{l}\mathtt{s}^t\mathtt{l}$ the number of short edges $t$ in that path is even.

\begin{figure}[ht!] \centering \includegraphics{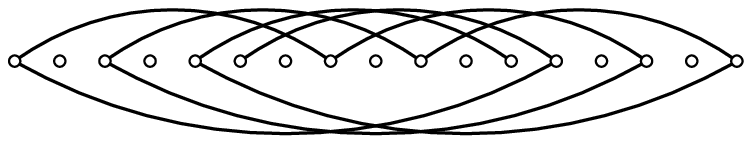} \caption{} \label{pgw} \end{figure}

For example, Figure \ref{pgw} shows one weave in $\Pi(7, 12, 17)$.

We say that a weave is \emph{simple} when all of its edges are short and \emph{compound} otherwise, when it contains both short and long edges. Thus in the case of $p = 1$ all weaves are simple.

We define also a corresponding class of connected components of leapers on rectangular boards, as follows.

Let $G$ be a connected leaper graph of $L$ on $A$ and let $G_X$ be the $x$-projection of $G$. Then $G_X$ is a connected subgraph of $\Pi_X$. We denote the connected component of $\Pi_X$ which contains $G_X$ by $\CompX(A, G)$ and we call it the \emph{$x$-completion} of $G$ on $A$. We define the \emph{$y$-completion} of $G$ on $A$ similarly, and we denote it by $\CompY(A, G)$.

When the board is clear from context, sometimes we omit it and we write simply $\CompX(G)$ or $\CompY(G)$.

Let $D$ be a connected component of $L$ on $A$. We say that $D$ is a \emph{vertical weave} of $L$ on $A$ when $\CompX(A, D)$ is a weave in $\Pi_X$. We define a \emph{horizontal weave} of $L$ on $A$ similarly. When $D$ is either a horizontal weave or a vertical weave of $L$ on $A$, we say simply that $D$ is a \emph{weave} of $L$ on $A$.

The rest of this section collects a number of helpful lemmas about weaves in projection graphs. We discuss the weaves of leapers on rectangular boards in Section \ref{ang}.

\medskip

\begin{lemma} \label{wb} Suppose that $\Pi(p, q, s)$ contains a weave. Then $s \le p + q - 1$. Thus, in particular, every weave in a projection graph is a path. \end{lemma} 

\medskip

\begin{proof} We already know this in the case of $p = 1$.

Suppose that $p \ge 2$.

When $s \ge p + q$, by Lemma \ref{pa} we get that $\Pi(p, q, s)$ is connected. Then it cannot be a weave because it contains a path with the signature $\mathtt{l}\mathtt{s}^t\mathtt{l}$ for both $t = \lfloor q/p \rfloor$ and $t = \lceil q/p \rceil$.

Otherwise, when $s \le p + q - 1$, every connected component of $\Pi(p, q, s)$ is a path by Lemma \ref{pa}. \end{proof}

\medskip

We proceed to examine the structure of compound weaves. Thus suppose, throughout the rest of this section, that $p \ge 2$.

Let $r = q \mmod 2p$ and $q = 2cp + \varepsilon r$ so that $\varepsilon \in \{-1, 1\}$ and $c$ is a positive integer. Since $p \ge 2$, with $\ecf(L) = [c_\kappa, \varepsilon_{\kappa - 1}, c_{\kappa - 1}, \varepsilon_{\kappa - 2}, \ldots, c_1]$ we get that $r = p_{\kappa - 1} = q_{\kappa - 2}$, $\varepsilon = \varepsilon_{\kappa - 1}$, and $c = c_\kappa/2$.

Consider a compound weave $W$ in $\Pi(p, q, s)$. Then $W$ is a path by Lemma \ref{wb}.

Observe that each endpoint of a long edge in a projection graph is also incident with at least one short edge. Therefore, we can group the short edges of $W$ into subpaths $W_0$, $W_1$, \ldots, $W_d$ so that a long edge $e_i$ of $W$ joins the rightmost vertex of $W_i$ and the leftmost vertex of $W_{i + 1}$ for all $i$ with $0 \le i < d$ and so that this accounts for all long edges of $W$. Let $s_i$ be the number of short edges in $W_i$, for all $i$.

\medskip

\begin{lemma} \label{cws} Suppose that $p \ge 2$. When $\varepsilon = -1$, the lengths $s_i$ in the above setting are given by $s_0 = s_d = 2c - 1$ and $s_1 = s_2 = \cdots = s_{d - 1} = 2c$. Otherwise, when $\varepsilon = 1$, they are given by $s_0 = s_1 = \cdots = s_d = 2c$. \end{lemma} 

\medskip

\begin{proof} Let $u$ be any vertex incident with a long edge in $\Pi(p, q, s)$. Then the longest path of short edges in $\Pi(p, q, s)$ which contains $u$ has length at least $\lfloor q/p \rfloor$. On the other hand, every path of short edges in $\Pi(p, q, s)$ has length at most $\lfloor (s - 1)/p \rfloor$. By Lemma \ref{wb}, the latter does not exceed $\lceil q/p \rceil$.

Thus for all $i$ we get that $s_i \in \{\lfloor q/p \rfloor, \lceil q/p \rceil\}$. However, in fact $\{\lfloor q/p \rfloor, \lceil q/p \rceil\} = \{2c, 2c + \varepsilon\}$. Since $W$ is a weave, all of $s_1$, $s_2$, \ldots, $s_{d - 1}$ must be even. Therefore, $s_1 = s_2 = \cdots = s_{d - 1} = 2c$.

We are left to determine $s_0$ and $s_d$. By the definition of $W_0$, its leftmost vertex is not incident with a long edge of $\Pi(p, q, s)$. Consequently, the difference between the two outermost vertices of $W_0$ must be less than $q$. This rules out the possibility of $s_0 = \lceil q/p \rceil$. Similarly, necessarily $s_d = \lfloor q/p \rfloor$ as well. \end{proof}

\medskip

In light of Lemma \ref{cws}, we denote the vertices of $W$ as follows: When $\varepsilon = 1$, we let $W_i = w_{i, 0}w_{i, 1} \ldots w_{i, 2c}$ for all $i$. Otherwise, when $\varepsilon = -1$, we let $W_i = w_{i, 0}w_{i, 1} \ldots w_{i, 2c}$ for all $i$ with $0 < i < d$ as well as $W_0 = w_{0, 1}w_{0, 2} \ldots w_{0, 2c}$ and $W_d = w_{d, 0}w_{d, 1} \ldots w_{d, 2c - 1}$. We enumerate the vertices of each subpath from left to right, so that the long edges of $W$ become $e_i = w_{i, 2c}w_{i + 1, 0}$ for all $i$ with $0 \le i < d$.

Our notation ensures that, for all $i$ and $j$, the subpath of $W$ which leads from $w_{i, j}$ to $w_{i + 1, j}$ contains a total of $2c$ short moves to the right and one long move to the left. This fact has some useful implications, as follows.

Consider two vertices $w_{i', j'}$ and $w_{i'', j''}$ of $W$. Then \[w_{i', j'} - w_{i'', j''} = -\varepsilon(i' - i'')r + (j' - j'')p,\] since $w_{i + 1, j} - w_{i, j} = 2cp - q = -\varepsilon r$ and $w_{i, j + 1} - w_{i, j} = p$ for all $i$ and $j$ such that the left-hand side is well-defined.

Furthermore, let $\nu$ be a walk connecting $w_{i', j'}$ and $w_{i'', j''}$ within $W$. Then the number of short moves in $\nu$ is of the same parity as $j' + j''$; the number of long moves in $\nu$ is of the same parity as $i' + i''$; and the total length of $\nu$ is of the same parity as $i' + i'' + j' + j''$.

\medskip

\begin{lemma} \label{cwl} Suppose that $p \ge 2$. Then every compound weave in a projection graph of the form $\Pi(p, q, s)$ contains fewer than $p/r$ long edges. \end{lemma} 

\medskip

\begin{proof} Consider a compound weave $W$ in $\Pi(p, q, I)$.

Suppose first that $\varepsilon = -1$. Then $w_{d, 2c - 1} - q < \min I \le w_{1, 0}$. Consequently, $q > w_{d, 2c - 1} - w_{1, 0} = (d - 1)r + (2c - 1)p$. This works out to $dr < p$.

Suppose, otherwise, that $\varepsilon = 1$. Then $w_{0, 0} - p < \min I \le w_{d, 0}$. Consequently, $p > w_{0, 0} - w_{d, 0} = dr$ once again. \end{proof}

\medskip

It is straightforward to verify that, when $p \ge 2$, for every positive integer smaller than $p/r$ there does exist a compound weave in a projection graph of the form $\Pi(p, q, s)$ with that many long edges.

We collect some more results regarding weaves in Appendix \ref{weave-ii}.

%% file: induction-on-descent-06-ang.tex
\section{Angularity} \label{ang}

We define the \emph{angularity graph} $\mathcal{N}$ of $L$ on $A$ as follows: The vertices of $\mathcal{N}$ are all edges of $L$ on $A$, and two edges of $L$ on $A$ are joined by an edge in $\mathcal{N}$ if and only if they form an acute angle.

We say that a set of edges of $L$ on $A$ is \emph{angular} when it is the vertex set of a connected subgraph of $\mathcal{N}$.

Let $G$ be a leaper graph of $L$ on $A$ without isolated vertices. We say that $G$ is \emph{angular} when its edge set is angular. We also say that $G$ is an \emph{angular component} of $L$ on $A$ when its edges are the vertices of a connected component of $\mathcal{N}$. Thus the angular components of $L$ on $A$ form a partitioning of the edges of $L$ on $A$.

\begin{figure}[ht!] \centering \includegraphics{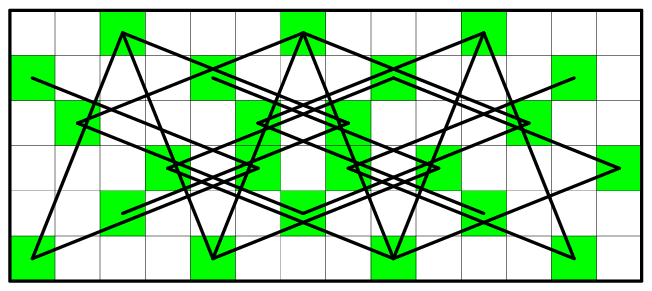} \caption{} \label{ac} \end{figure}

For example, Figure \ref{ac} shows one of the four angular components of the $(2, 5)$-leaper on the board of size $6 \times 14$.

We say that a walk $\alpha = a_1a_2 \ldots a_k$ of $L$ on $A$ is \emph{angular} when each angle of the form $a_{i - 1}a_ia_{i + 1}$ in it is either acute or zero. Thus two edges of $L$ on $A$ are in the same angular component if and only if some angular walk contains both of them. We also say that $\alpha$ is \emph{angular-closed} when: (a) It is closed; (b) It is angular; and (c) Its first edge and its last edge form either an acute or a zero angle as well.

When angle $b'ab''$ of $L$ on $A$ is either acute or zero, edges $ab'$ and $ab''$ are in the same angular component by definition. We go on to analyse the other types of angles as well.

\medskip

\begin{lemma} \label{rdo} Let $b'ab''$ be either a right or a diagonally obtuse angle of $L$ on $A$. Then edges $ab'$ and $ab''$ are in the same angular component. \end{lemma} 

\medskip

Consequently, two edges of $L$ on $A$ incident with the same cell and of different inclines are always in the same angular component.

\medskip

\begin{proof} When angle $b'ab''$ is right, up to symmetry $b' = a + (q, p)$ and $b'' = a + (-p, q)$. Then cell $c = a + (p, q)$ is in $A$ as well, and both angles $b'ac$ and $cab''$ are acute.

When angle $b'ab''$ is diagonally obtuse, up to symmetry $b' = a + (q, -p)$ and $b'' = a + (-p, q)$. Then both cells $c' = a + (q, p)$ and $c'' = a + (p, q)$ are in $A$ as well, and all three angles $b'ac'$, $c'ac''$, and $c''ab''$ are acute. \end{proof}

\medskip

\begin{lemma} \label{los} Let $b'ab''$ be either a laterally obtuse or a straight angle of $L$ on $A$. Suppose, for concreteness, that both edges $ab'$ and $ab''$ are slight. Suppose, furthermore, that $\CompY(b'ab'')$ is a non-weave. Then edges $ab'$ and $ab''$ are in the same angular component. \end{lemma} 

\medskip

\begin{proof} First we consider the case when $p = 1$.

Since $\CompY(b'ab'')$ is a non-weave, $m \ge q + 1$. Furthermore, since both edges $ab'$ and $ab''$ are slight, $n \ge 2q + 1$. We claim that then all edges of $L$ on $A$ are in the same angular component.

Indeed, consider any subboard $A^\star$ of $A$ of size $(q + 1) \times (q + 1)$. Since the edges of $L$ on $A^\star$ form a cycle all of whose angles are acute, they are in the same angular component. Furthermore, whenever $A'$ and $A''$ are two subboards of $A$ of size $(q + 1) \times (q + 1)$ whose central cells are adjacent by side, at least one edge of $L$ is contained within both of them.

Suppose, throughout the rest of the proof, that $p \ge 2$.

Let $a$, $b'$, and $b''$ be, up to symmetry, in columns $z$, $z - q$, and $z + q$ of $A$, respectively.

Since $p \ge 2$ and $\CompY(b'ab'')$ is a non-weave, there exists a walk $\upsilon = y_1y_2 \ldots y_k$ in $\Pi_Y$ such that the signature of $\upsilon$ begins $\mathtt{l}\mathtt{s}^t\mathtt{l}\ldots$ for some odd positive integer $t$ and $a$ is in row $y_k$ of $A$. By inserting some short-edge detours into it as needed, we can also ensure without loss of generality that $\upsilon$ does not contain two consecutive long edges.

\begin{figure}[ht!] \centering \includegraphics{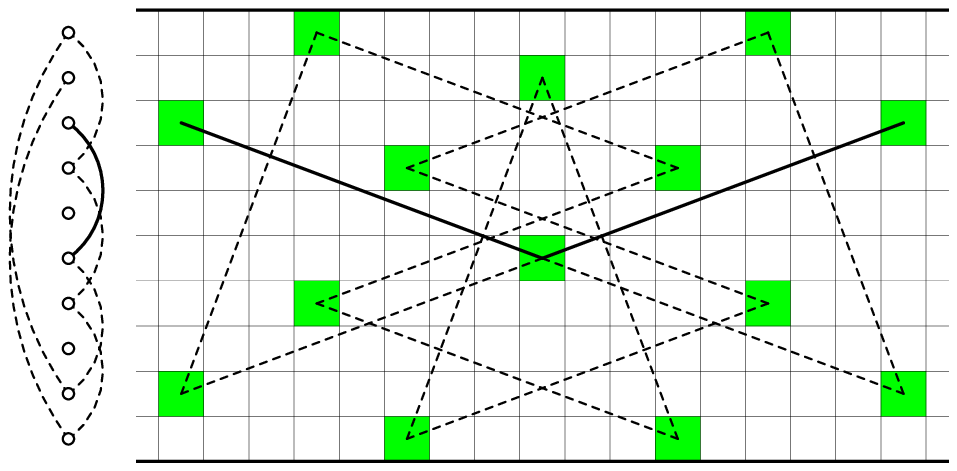} \caption{} \label{lo-nw} \end{figure} 

Consider the cycles $C' = z$---$(z + p)$---$(z + p - q)$---$(z - q)$---$z$ and $C'' = z$---$(z - p)$---$(z - p + q)$---$(z + q)$---$z$ in $\Pi_X$. By Lemma \ref{pf}, there exists a unique walk $\alpha' = a'_1a'_2 \ldots a'_k$ of $L$ on $A$ such that the $x$-projection of $\alpha'$ is contained within $C'$, the $y$-projection of $\alpha'$ coincides with $\upsilon$, and $a'_k = a$. Define $\alpha''$ similarly. Thus $\alpha'$ and $\alpha''$ are reflections of one another with respect to the vertical line through $a$.

Since $\upsilon$ does not contain two consecutive long edges, no angle of $\alpha'$ and $\alpha''$ is either laterally obtuse or straight. By Lemma \ref{rdo}, it follows that all edges of $\alpha'$ are in the same angular component which also contains edge $ab'$. Similarly, all edges of $\alpha''$ are in the same angular component which also contains edge $ab''$. We are left to show that these two angular components coincide.

Let $a'_i = (x'_i, y_i)$ and $a''_i = (x''_i, y_i)$ for all $i$. Since $t$ is odd, either $x'_2 \in \{z, z + p\}$ or $x'_{t + 2} \in \{z, z + p\}$. We consider the former case for concreteness, and the latter one is analogous.

Observe that $x'_2 = z$ implies $a'_2 = a''_2$ and $x'_2 = z + p$ implies $a'_1 = a''_1$. Either way, edges $a'_1a'_2$ and $a''_1a''_2$ form a laterally acute angle whose angular component contains both of $\alpha'$ and $\alpha''$. \end{proof}

\medskip

For example, Figure \ref{lo-nw} shows the construction in the proof carried out with a laterally obtuse angle of the $(3, 8)$-leaper.

Clearly, every angular leaper graph of $L$ on $A$ is connected. Thus, in particular, every angular component of $L$ on $A$ is connected as well. However, it is false in general that every non-singleton connected component of $L$ on $A$ is angular. (Though if a non-singleton connected component of $L$ on $A$ is indeed angular, then of course it must also be an angular component of $L$ on $A$.)

We proceed to examine the structure of the non-angular connected components of $L$ on $A$. Note that, since our definition of angularity applies only to leaper graphs without isolated vertices, a singleton connected component of $L$ on $A$ is neither angular nor non-angular.

The main result of this section is as follows.

\medskip

\begin{lemma} \label{naw} Every non-angular connected component of a skew free leaper on a rectangular board is a weave. \end{lemma} 

\medskip

\begin{proof} Let $C$ be a non-angular connected component of $L$ on $A$. Then $C$ contains two edges $ab'$ and $ab''$ in different angular components. By Lemma \ref{rdo}, angle $b'ab''$ is either laterally obtuse or straight. Thus up to symmetry both edges $ab'$ and $ab''$ are slight. Consequently, by Lemma \ref{los} we obtain that $\CompY(C)$ is a weave. \end{proof}

\medskip

Consider a non-angular connected component $C$ of $L$ on $A$. As in the proof of Lemma \ref{naw}, there exist two edges $ab'$ and $ab''$ of $C$ in different angular components such that either both of them are slight or both of them are steep. In the former case, it follows that $C$ is a horizontal weave, $m \le p + q - 1$ by Lemma \ref{wb}, and $n \ge 2q + 1$ because angle $b'ab''$ is contained within $A$. Conversely, in the latter case it follows that $C$ is a vertical weave, $m \ge 2q + 1$, and $n \le p + q - 1$.

Consequently, either $m < n$ and $C$ is horizontal weave or $m > n$ and $C$ is a vertical weave. Thus, in particular, $C$ cannot be simultaneously a horizontal weave and a vertical weave.

We collect some more results regarding angularity in Appendix \ref{ang-ii}.

%% file: induction-on-descent-07-desc.tex
\section{Descents} \label{desc}

We define three \emph{lifting transformations} which map skew free leapers onto skew free leapers, as follows. \cite{Be}

(a) $f(L)$ is the $(p, 2p + q)$-leaper.

(b) $g(L)$ is the $(q, 2q - p)$-leaper.

(c) $h(L)$ is the $(q, p + 2q)$-leaper.

Let $\ell$ be one of the three lifting transformations. Then we call $\ell(L)$ the \emph{$\ell$-child} of $L$.

When $L$ is distinct from the knight, there exist a unique lifting transformation $\ell$ and a unique skew free leaper $K$ such that $\ell(K) = L$. Then we say that $L$ is of \emph{type} $\ell$ and we call $K$ the \emph{parent} of $L$. Explicitly:

(a) $L$ is of type $f$ if and only if $3p < q$. Then $K$ is the $(p, q - 2p)$-leaper.

(b) $L$ is of type $g$ if and only if $2p > q$. Then $K$ is the $(2p - q, p)$-leaper.

(c) $L$ is of type $h$ if and only if $2p < q < 3p$. Then $K$ is the $(q - 2p, p)$-leaper.

Thus in all cases $K$ is the leaper with parameters $p$ and $|2p - q|$.

Observe that the sum of the parameters of the parent is always smaller than the sum of the parameters of the child. Therefore, starting from an arbitrary skew free leaper and moving from child to parent repeatedly, eventually we must reach the knight. This allows us to place all skew free leapers in an infinite ternary tree $\mathfrak{L}$ with the knight at its root.

\begin{figure}[t!] \centering \includegraphics{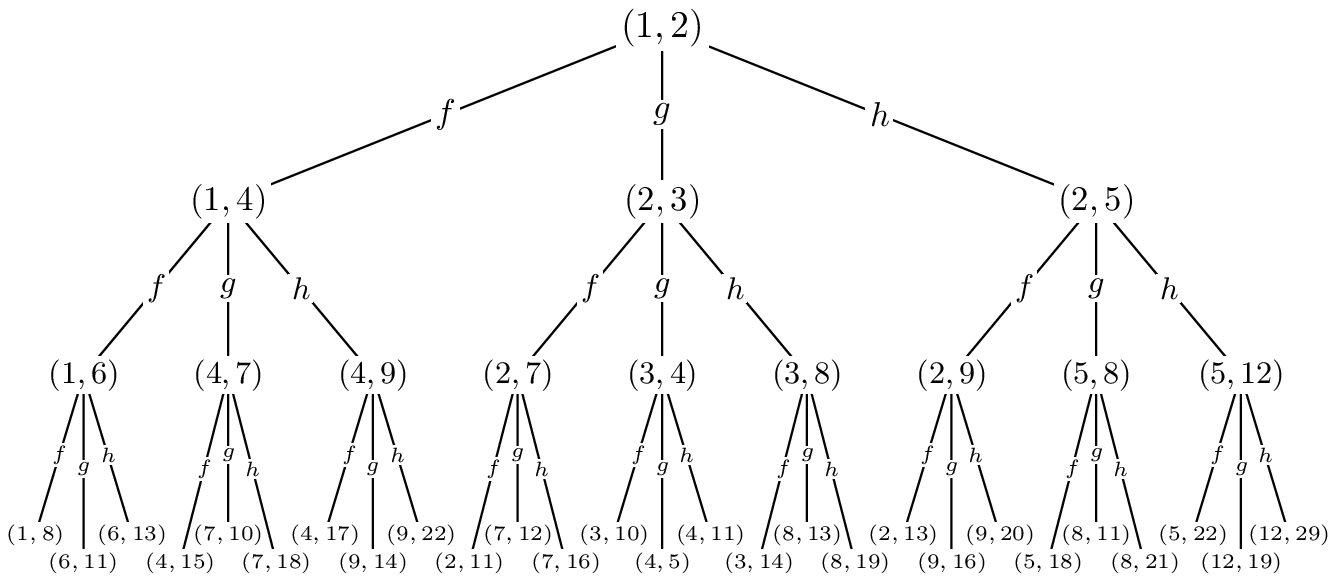} \caption{} \label{tree} \end{figure}

For example, Figure \ref{tree} shows the first four levels of $\mathfrak{L}$.

The \emph{descent} of a skew free leaper is a finite word over the alphabet $\{\mathtt{f}, \mathtt{g}, \mathtt{h}\}$, defined inductively as follows. The descent of the knight is the empty word, and if $w$ is the descent of $L$, then the descents of its three children are $\mathtt{f}w$, $\mathtt{g}w$, and $\mathtt{h}w$, respectively.

The descent of a skew free leaper is its ``address'' within $\mathfrak{L}$. It encodes the unique sequence of lifting transformations which leads from the knight to $L$. Thus we obtain a bijection between all skew free leapers and all finite words over the alphabet $\{\mathtt{f}, \mathtt{g}, \mathtt{h}\}$.

Let $\ecf(L) = [c_\kappa, \varepsilon_{\kappa - 1}, c_{\kappa - 1}, \varepsilon_{\kappa - 2}, \ldots, c_1]$. Then the even continued fractions of the children of $L$ are as follows.

(a) $\ecf(f(L)) = [2 + c_\kappa, \varepsilon_{\kappa - 1}, c_{\kappa - 1}, \varepsilon_{\kappa - 2}, \ldots, c_1]$.

(b) $\ecf(g(L)) = [2, -1, c_\kappa, \varepsilon_{\kappa - 1}, c_{\kappa - 1}, \varepsilon_{\kappa - 2}, \ldots, c_1]$.

(c) $\ecf(h(L)) = [2, 1, c_\kappa, \varepsilon_{\kappa - 1}, c_{\kappa - 1}, \varepsilon_{\kappa - 2}, \ldots, c_1]$.

Consequently, from the even continued fraction of $L$ we can obtain its descent by replacing each coefficient $c_i$ with the subword $\mathtt{f}^{c_i/2 - 1}$, each minus sign with the letter $\mathtt{g}$, and each plus sign with the letter $\mathtt{h}$.

Let $L'$ and $L''$ be two skew free leapers. We say that $L'$ is an \emph{ancestor} of $L''$ when there exists a sequence of lifting transformations $\ell_1$, $\ell_2$, \ldots, $\ell_k$ with $\ell_k(\ell_{k - 1}(\ldots \ell_1(L') \ldots)) = L''$. Then we also say that $L''$ is a \emph{descendant} of $L'$. Equivalently, $L'$ is an ancestor of $L''$ and $L''$ is a descendant of $L'$ if and only if the descent of $L'$ is a suffix of the descent of $L''$.

A tail is a special case of an ancestor. In terms of descents, $L'$ is a tail of $L''$ if and only if the descent of $L'$ is a suffix of the descent of $L''$ which is not preceded by the letter $\mathtt{f}$.

Before we continue, let us look into how some of the concepts in the introduction and in this section extend to skew half-free leapers.

Let $N$ be a skew half-free $(r, s)$-leaper with $r < s$. Then the rational number $s/r$ admits a unique expansion of the form $d_\lambda \pm 1/(d_{\lambda - 1} \pm 1/(\ldots d_1 \ldots))$ such that $d_1$ is an odd positive integer with $d_1 \ge 3$ and $d_i$ is an even positive integer for all $i$ with $2 \le i \le \lambda$. We call this the \emph{even continued fraction} of $N$ and we call $\lambda$ its \emph{depth}. We define the \emph{tails} of $N$ in the same way as with skew free leapers, except that we do not add a zeroth tail.

The three lifting transformations map skew half-free leapers onto skew half-free leapers. Furthermore, when $N$ is distinct from the camel, there exist a unique lifting transformation $\ell$ and a unique skew half-free leaper $K$ with $\ell(K) = N$.

Thus the notions of parent, child, and descent all carry over directly. Consequently, all skew half-free leapers form an infinite ternary tree $\mathfrak{L}^\text{Half}$ with the camel at its root in the same way as all skew free leapers form $\mathfrak{L}$. This allows us to use induction on descent with skew half-free leapers in the same manner as with skew free leapers.

The material in Sections \ref{lift}--\ref{line} also carries over to skew half-free leapers without any substantial changes. We will return to these points in Sections \ref{rigid} and \ref{further}.

%% file: induction-on-descent-08-lift.tex
\section{Lifting} \label{lift}

We go on to extend the definitions of the three lifting transformations so that they can act on a few other kinds of objects beyond leapers.

Let $\ell$ be a lifting transformation. Given $\ell$ and $L$, we define two positive integers, the \emph{lateral parameter} $\Lat$ and the \emph{diagonal parameter} $\Dia$ of $\ell$ and $L$, as in Table~\ref{ld}.

\begin{table}[!ht]
	\centering
	\begin{tabular}{c | c | c}
		$\ell$ & $\Lat$ & $\Dia$\\
		\hline
		$f$ & $p + q$ & $p$\\
		$g$ & $q - p$ & $q$\\
		$h$ & $p + q$ & $q$
	\end{tabular}
	\caption{}
	\label{ld}
\end{table} 

We define one more positive integer, the \emph{margin parameter} $\Mar$ of $\ell$ and $L$, by $\Mar = \min\{\Lat, \Dia\}$.

Thus for all three lifting transformations the $\ell$-child of $L$ is the $(\Dia, \Lat + \Dia)$-leaper. For convenience, let $\widehat{p} = \Dia$, $\widehat{q} = \Lat + \Dia$, and $M = \ell(L)$, so that $M$ is the skew free $(\widehat{p}, \widehat{q})$-leaper.

Given $\ell$ and $L$, let $J_X = [x_\text{Min} - \Mar; x_\text{Max} + \Mar]$ and $J_Y = [y_\text{Min} - \Mar; y_\text{Max} + \Mar]$. Then we define $\ell(L, A)$ to be the board $B = J_X \times J_Y$.

Let $\widehat{m}$ and $\widehat{n}$ be the height and width of $B$, respectively. Then the sides of $A$ and $B$ are related by $\widehat{m} = m + 2\Mar$ and $\widehat{n} = n + 2\Mar$.

When the leaper is clear from context, sometimes we omit it and we write simply $\ell(A) = B$.

Observe that $A' \le A''$ if and only if $\ell(A') \le \ell(A'')$ and $A' \sqsubseteq A''$ if and only if $\ell(A') \sqsubseteq \ell(A'')$. Furthermore, $A$ is standard if and only if $\ell(A)$ is standard.

Let $a$ be a cell and let $b'b''$ be an edge of $M$. We say that $a$ and $b'b''$ are \emph{matched} when one of $ab'$ and $ab''$ is an edge of the $(0, \Lat)$-leaper and the other one is an edge of the $(\Dia, \Dia)$-leaper. Thus on the infinite board $\mathbb{Z} \times \mathbb{Z}$ each cell is matched with eight edges of $M$ and each edge of $M$ is matched with two cells.

\begin{figure}[ht!] \centering \begin{tabular}{>{\centering}m{0.3\textwidth} >{\centering}m{0.225\textwidth} >{\centering}m{0.3\textwidth}} \includegraphics{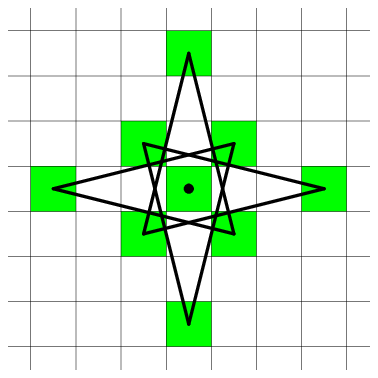} & \includegraphics{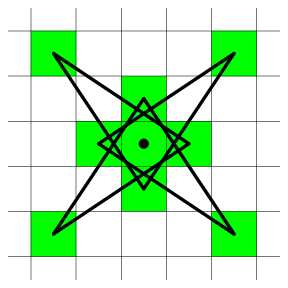} & \includegraphics{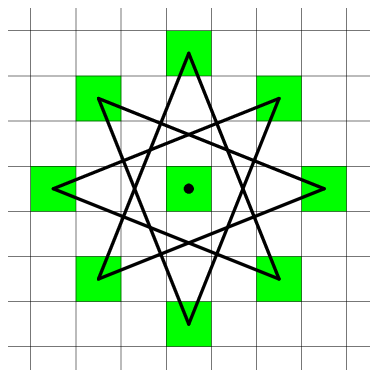} \end{tabular} \caption{} \label{match-c} \end{figure}

For example, each part of Figure \ref{match-c} shows one cell and the eight edges of $M$ matched with it. In all three parts, $L$ is the knight. On the left, $\ell = f$ and $M$ is the giraffe; in the middle, $\ell = g$ and $M$ is the zebra; and, on the right, $\ell = h$ and $M$ is the $(2, 5)$-leaper.

\begin{figure}[ht!] \centering \begin{tabular}{>{\centering}m{0.3\textwidth} >{\centering}m{0.225\textwidth} >{\centering}m{0.3\textwidth}} \includegraphics{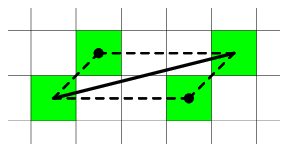} & \includegraphics{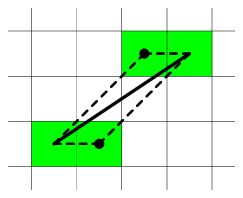} & \includegraphics{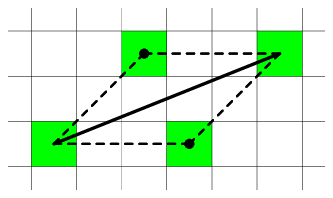} \end{tabular} \caption{} \label{match-e} \end{figure}

Conversely, Figure \ref{match-e} shows one edge of $M$ and the two cells matched with it in the same three settings.

Given a cell $a$, we define cells $a^\text{E}$, $a^\text{NE}$, \ldots, $a^\text{SE}$ by \begin{align*} a^\text{E} &= a + (\Lat, 0), & a^\text{NE} &= a + (\Dia, \Dia),\\ a^\text{N} &= a + (0, \Lat), & a^\text{NW} &= a + (-\Dia, \Dia),\\ a^\text{W} &= a + (-\Lat, 0), & a^\text{SW} &= a + (-\Dia, -\Dia),\\ a^\text{S} &= a + (0, -\Lat), & a^\text{SE} &= a + (\Dia, -\Dia). \end{align*}

(The superscripts correspond to the primary and secondary cardinal directions.)

Then the edges of $M$ matched with $a$ are $a^\text{E}a^\text{NW}$, $a^\text{NW}a^\text{S}$, \ldots, $a^\text{SW}a^\text{E}$. Observe that they form a cycle of length eight all of whose angles are acute.

Conversely, let $b'b''$ be an edge of $M$. Up to symmetry, $b'' = b' + (\widehat{q}, \widehat{p})$. Let $a' = b' + (\Lat, 0) = b'' + (-\Dia, -\Dia)$ and $a'' = b' + (\Dia, \Dia) = b'' + (-\Lat, 0)$. Then $a'$ and $a''$ are the two cells matched with edge $b'b''$. Observe that $a'a''$ is an edge of $L$.

Let $e'$ be an edge of $L$ and let $e''$ be an edge of $M$. When both endpoints of $e'$ are matched with $e''$, we say that $e'$ and $e''$ are \emph{cross-edges}. When $\ell$ and $L$ are fixed, every edge of $L$ has a unique cross-edge of $M$ and, conversely, every edge of $M$ has a unique cross-edge of $L$. Thus we obtain a bijection between the edges of $L$ and the edges of $M$ on the infinite board $\mathbb{Z} \times \mathbb{Z}$.

Let $S$ be a set of cells on $A$. Then we define $\ell(L, A, S)$ to be the set $T$ of all edges of $M$ on $B$ matched with a cell in $S$.

Finally, let $G$ be a leaper graph of $L$ on $A$. Then we define $\ell(L, A, G)$ to be the leaper graph $H$ of $M$ on $B$ formed by all edges of $M$ on $B$ matched with a vertex of $G$. Thus also lifting the vertex set of $G$ yields the edge set of $H$.

When the leaper and board are clear from context, sometimes we omit them and we write simply $\ell(S) = T$ and $\ell(G) = H$.

When $A$ is the square board of side $p + q$ and $G$ is a clover of $L$ on $A$, the definitions of $f(G)$, $g(G)$, and $h(G)$ can be greatly simplified. This special case is considered in \cite{Be} in the context of the second leaper theorem.

For convenience, throughout the rest of this section the symbols $\ell$, $L$, $M$, $A$, $B$, $S$, $T$, $G$, and $H$ will retain their current meanings.

\medskip

\begin{lemma} \label{cross} Let $e'$ be an edge of $L$ on $A$. Then the cross-edge $e''$ of $e'$ is an edge of $M$ on $B$. \end{lemma}

\medskip

\begin{proof} Let $e' = ab$. By the definition of a cross-edge, one endpoint of $e''$ is in $\Ball(a, \Lat) \cap \Ball(b, \Dia)$ and the other one is in $\Ball(a, \Dia) \cap \Ball(b, \Lat)$. Thus one endpoint of $e''$ is in $\Ball(a, \Mar)$ and the other one is in $\Ball(b, \Mar)$. Consequently, both endpoints of $e''$ are in $B$. \end{proof}

\medskip

The cross-edge correspondence reveals a close connection between the angular components of $L$ on $A$ and the angular components of $M$ on $B$.

\medskip

\begin{lemma} \label{ace} Suppose that $G$ is an angular leaper graph of $L$ on $A$. Then the cross-edges of the edges of $G$ are all in the same angular subgraph of $H$. \end{lemma} 

\medskip

Thus, in particular, if two edges are in the same angular component of $L$ on $A$ then their cross-edges will be in the same angular component of $M$ on $B$.

\medskip

\begin{proof} It suffices to consider the case when the two edges of $L$ on $A$ form an acute angle. Denote them by $ab'$ and $ab''$. We claim that in fact there exists an angular path of $M$ on $B$ which consists entirely of edges matched with $a$ and which contains the cross-edges of $ab'$ and $ab''$. We consider two cases for angle $b'ab''$ multiplied by three cases for the lifting transformation $\ell$.

\smallskip

\emph{Case 1}. Angle $b'ab''$ is laterally acute. Then up to symmetry $b' = a + (q, -p)$ and $b'' = a + (q, p)$.

\smallskip

\emph{Case 1.1}. $\ell = f$. The path is $a^\text{SW}a^\text{E}a^\text{NW}$.

\smallskip

\emph{Case 1.2}. $\ell = g$. The path is $a^\text{N}a^\text{SE}a^\text{W}a^\text{NE}a^\text{S}$.

\smallskip

\emph{Case 1.3}. $\ell = h$. The path is $a^\text{NE}a^\text{S}a^\text{NW}a^\text{E}a^\text{SW}a^\text{N}a^\text{SE}$.

\smallskip

\emph{Case 2}. Angle $b'ab''$ is diagonally acute. Then up to symmetry $b' = a + (q, p)$ and $b'' = a + (p, q)$.

\smallskip

\emph{Case 2.1}. $\ell = f$. The path is $a^\text{NW}a^\text{E}a^\text{SW}a^\text{N}a^\text{SE}$.

\smallskip

\emph{Case 2.2}. $\ell = g$. The path is $a^\text{S}a^\text{NE}a^\text{W}$.

\smallskip

\emph{Case 2.3}. $\ell = h$. The path is $a^\text{SE}a^\text{N}a^\text{SW}a^\text{E}a^\text{NW}$. \end{proof}

\medskip

We proceed to study the effect of the three lifting transformations on some fundamental properties of leaper graphs. The case of $f$ is more complicated, and so we begin with $g$ and $h$.

\medskip

\begin{lemma} \label{lg} Suppose that $\ell = g$ and $G$ is a connected leaper graph of $L$ on $A$ which contains at least two cells. Then $H$ is an angular leaper graph of $M$ on $B$. Thus, in particular, $H$ is connected as well. \end{lemma} 

\medskip

The condition that $G$ contains at least two cells can be relaxed somewhat. With minor modifications, the proof continues to go through when $G$ consists of a single cell but both of $\CompX(A, G)$ and $\CompY(A, G)$ are non-singletons.

\medskip

\begin{proof} It suffices to check that: (a) For every cell $a$ of $G$, the edges of $g(\{a\})$ form an angular leaper graph of $M$ on $B$; and (b) For every edge $a'a''$ of $G$, the sets of edges $g(\{a'\})$ and $g(\{a''\})$ have at least one element in common.

Claim (b) is clear: The cross-edge of $a'a''$ lies within $B$ by Lemma \ref{cross}, and it belongs to both sets $g(\{a'\})$ and $g(\{a''\})$.

Claim (a) is not so straightforward.

When $\ell = g$, we get that $\Mar = \Lat$. Thus all four of $a^\text{E}$, $a^\text{N}$, $a^\text{W}$, and $a^\text{S}$ are in $B$. Up to symmetry, there are six cases to consider as to which ones out of $a^\text{NE}$, $a^\text{NW}$, $a^\text{SW}$, and $a^\text{SE}$ are in $B$ as well.

\smallskip

\emph{Case 1}. None of $a^\text{NE}$, $a^\text{NW}$, $a^\text{SW}$, and $a^\text{SE}$ are in $B$. But then $a$ cannot be incident with any edges of $L$ on $A$, and so this case cannot occur.

\smallskip

\emph{Case 2}. Only $a^\text{NE}$ is in $B$. Then $g(\{a\})$ consists of two edges of $M$ which form the acute angle $a^\text{S}a^\text{NE}a^\text{W}$.

\smallskip

\emph{Case 3}. Cells $a^\text{NE}$ and $a^\text{NW}$ are in $B$ whereas cells $a^\text{SW}$ and $a^\text{SE}$ are outside of $B$. Then $g(\{a\})$ consists of four edges of $M$ which form the angular path $a^\text{E}a^\text{NW}a^\text{S}a^\text{NE}a^\text{W}$.

\smallskip

\emph{Case 4}. Cells $a^\text{NE}$ and $a^\text{SW}$ are in $B$ whereas cells $a^\text{NW}$ and $a^\text{SE}$ are outside of $B$. This case cannot occur because cells $a^\text{NW}$ and $a^\text{SE}$ are in the bounding box of cells $a^\text{NE}$ and $a^\text{SW}$.

\smallskip

\emph{Case 5}. Only $a^\text{NE}$ is outside of $B$. This case cannot occur, either, for the same reason as Case 4.

\smallskip

\emph{Case 6}. All of $a^\text{NE}$, $a^\text{NW}$, $a^\text{SW}$, and $a^\text{SE}$ are in $B$. Then $g(\{a\})$ consists of eight edges of $M$ which form a cycle all of whose angles are acute. \end{proof}

\medskip

\begin{lemma} \label{lh} Suppose that $\ell = h$ and $G$ is a connected leaper graph of $L$ on $A$ which contains at least two cells. Then $H$ is an angular leaper graph of $M$ on $B$. Thus, in particular, $H$ is connected as well. \end{lemma} 

\medskip

Just as with Lemma \ref{lg}, the proof continues to go through with minor modifications when $G$ consists of a single cell but both of $\CompX(A, G)$ and $\CompY(A, G)$ are non-singletons.

\medskip

\begin{proof} It suffices to establish the analogues of claims (a) and (b) from the proof of Lemma \ref{lg}. The verification of claim (b) remains unchanged, and so we turn to claim (a).

When $\ell = h$, we get that $\Mar = \Dia$. Thus all four of $a^\text{NE}$, $a^\text{NW}$, $a^\text{SW}$, and $a^\text{SE}$ are in $B$. Up to symmetry, there are six cases to consider as to which ones out of $a^\text{E}$, $a^\text{N}$, $a^\text{W}$, and $a^\text{S}$ are in $B$ as well.

\smallskip

\emph{Case 1}. None of $a^\text{E}$, $a^\text{N}$, $a^\text{W}$, and $a^\text{S}$ are in $B$. But then $a$ cannot be incident with any edges of $L$ on $A$, and so this case cannot occur.

\smallskip

\emph{Case 2}. Only $a^\text{E}$ is in $B$. Then $h(\{a\})$ consists of two edges of $M$ which form the acute angle $a^\text{SW}a^\text{E}a^\text{NW}$.

\smallskip

\emph{Case 3}. Cells $a^\text{E}$ and $a^\text{N}$ are in $B$ whereas cells $a^\text{W}$ and $a^\text{S}$ are outside of $B$. Then $h(\{a\})$ consists of four edges of $M$ which form the angular path $a^\text{SE}a^\text{N}a^\text{SW}a^\text{E}a^\text{NW}$.

\smallskip

\emph{Case 4}. Cells $a^\text{E}$ and $a^\text{W}$ are in $B$ whereas cells $a^\text{N}$ and $a^\text{S}$ are outside of $B$. This case cannot occur, either, for the same reason as Case 1.

\smallskip

\emph{Case 5}. Only $a^\text{E}$ is outside of $B$. Then $h(\{a\})$ consists of six edges of $M$ which form the angular path $a^\text{NW}a^\text{S}a^\text{NE}a^\text{W}a^\text{SE}a^\text{N}a^\text{SW}$.

\smallskip

\emph{Case 6}. All of $a^\text{E}$, $a^\text{N}$, $a^\text{W}$, and $a^\text{S}$ are in $B$. Then $h(\{a\})$ consists of eight edges of $M$ which form a cycle all of whose angles are acute. \end{proof}

\medskip

When we attempt to put together a similar lemma involving $f$, we run into a problem: The analogue of claim (a) from the proofs of Lemmas \ref{lg} and \ref{lh} does not in general continue to hold. We get around this obstacle by strengthening our premise.

\medskip

\begin{lemma} \label{lf} Suppose that $\ell = f$ and $G$ is an angular connected component of $L$ on $A$. Then $H$ is an angular leaper graph of $M$ on $B$. Thus, in particular, $H$ is connected as well. \end{lemma} 

\medskip

\begin{proof} We follow the same approach as with Lemmas \ref{lg} and \ref{lh}, except that we modify the analogue of claim (a) as follows: (a$^\star$) For every cell $a$ of $G$, the edges of $f(\{a\})$ are all in the same angular subgraph of $H$.

Once again, the verification of claim (b) remains unchanged, and so we turn to claim (a$^\star$).

When $\ell = f$, we get that $\Mar = \Dia$. Thus all four of $a^\text{NE}$, $a^\text{NW}$, $a^\text{SW}$, and $a^\text{SE}$ are in $B$. Just as in the proof of Lemma \ref{lh}, up to symmetry there are six cases to consider as to which ones out of $a^\text{E}$, $a^\text{N}$, $a^\text{W}$, and $a^\text{S}$ are in $B$ as well.

We handle Cases 1--3 and 5--6 in the exact same way as before.

With this, we are left to take care of Case 4, when cells $a^\text{E}$ and $a^\text{W}$ are in $B$ whereas cells $a^\text{N}$ and $a^\text{S}$ are outside of $B$. Then the edges of $M$ in $f(\{a\})$ form the two disjoint laterally acute angles $a^\text{SW}a^\text{E}a^\text{NW}$ and $a^\text{NE}a^\text{W}a^\text{SE}$.

Let $a = (x, y)$. Since $G$ is angular, $a$ is incident with an edge of $L$ on $A$, and so up to symmetry $y + p \in I_Y$. Furthermore, $a^\text{E} \in B$ implies $x + q \in I_X$ and $a^\text{W} \in B$ similarly implies $x - q \in I_X$.

Therefore, both cells $b' = a + (q, p)$ and $b'' = a + (-q, p)$ are in $A$.

Since $G$ is a connected component of $L$ on $A$, it contains both edges $ab'$ and $ab''$. Observe, however, that their cross-edges are $a^\text{E}a^\text{NW}$ and $a^\text{NE}a^\text{W}$, respectively. Since $G$ is also an angular leaper graph of $L$ on $A$, by Lemma \ref{ace} it follows that $a^\text{E}a^\text{NW}$ and $a^\text{NE}a^\text{W}$ are in the same angular subgraph of $H$. \end{proof}

%% file: induction-on-descent-09-low.tex
\section{Lowering} \label{low}

For each lifting transformation $\ell$, we define also a complementary \emph{lowering transformation} $\ellinv$, as follows.

For skew free leapers $L$ and $M$, we let $\ellinv(M) = L$ if and only if $\ell(L) = M$. Thus $\ellinv(M)$ is well-defined if and only if $M$ is of type $\ell$.

Fix two skew free leapers $L$ and $M$ with $\ell(L) = M$ and $\ellinv(M) = L$. For boards $A$ and $B$, we let $\ellinv(M, B) = A$ if and only if $\ell(L, A) = B$.

When the leaper is clear from context, sometimes we omit it and we write simply $\ellinv(B) = A$.

Define $\Lat$, $\Dia$, $\Mar$, $\widehat{p}$, and $\widehat{q}$ as in Section \ref{lift}. Then $\ellinv(M, B)$ is well-defined if and only if both sides of $B$ exceed $2\Mar$. Since $2\Mar = 2\min\{\Lat, \Dia\} < \Lat + \Dia = \widehat{q}$, if $\ellinv(M, B)$ is not well-defined then all edges of $M$ on $B$ are of the same incline. Conversely, if there exist edges of $M$ on $B$ of both inclines then $\ellinv(M, B)$ is well-defined.

Fix two boards $A$ and $B$ with $\ell(L, A) = B$ and $\ellinv(M, B) = A$. Let also $T$ be a set of edges of $M$ on $B$. Then we define $\ellinv(M, B, T)$ to be the set $S$ of all cells of $A$ matched with an edge in $T$.

Finally, let $H$ be a leaper graph of $M$ on $B$. Then we define $\ellinv(M, B, H)$ to be the induced leaper graph of $L$ on $A$ whose vertices are all cells of $A$ matched with an edge of $H$. Thus also lowering the edge set of $H$ yields the vertex set of $G$.

Once again, when the leaper and board are clear from context, sometimes we omit them and we write simply $\ellinv(T) = S$ and $\ellinv(H) = G$.

Observe that lifting and lowering are exact inverses for leapers and boards but not in general for cell sets, edge sets, and leaper graphs. Thus $\ell(S) = T$ and $\ellinv(T) = S$ do not necessarily imply one another. Similarly, neither do $\ell(G) = H$ and $\ellinv(H) = G$.

Let $b'ab''$ be an acute angle of $M$. We define the \emph{core} of $b'ab''$ to be the unique cell matched with both edges $ab'$ and $ab''$.

\medskip

\begin{lemma} \label{core} Let $b'ab''$ be an acute angle of $M$ on $B$. Then its core is a cell of~$A$. \end{lemma}

\medskip

\begin{proof} We examine two cases for the type of angle $b'ab''$ multiplied by three cases for the lifting transformation $\ell$ in order to confirm that its core $c$ always satisfies $\Ball(c, \Mar) \subseteq \BBox(b'ab'')$. Since also $\BBox(b'ab'') \subseteq B$, it follows that $c$ is indeed a cell of $A$. \end{proof}

\medskip

The main result of this section is as follows.

\medskip

\begin{lemma} \label{alc} Let $H$ be an angular leaper graph of $M$ on $B$. Then $\ellinv(H)$ is a connected leaper graph of $L$ on $A$. \end{lemma} 

\medskip

This is, of course, the lowering complement of Lemmas \ref{lg}--\ref{lf}. Note that, if $H$ consists of a single edge, then $\ellinv(H)$ might happen to be the null graph.

\medskip

\begin{proof} It suffices to check that: (a) For every edge $e$ of $M$ on $B$, the cells of $\ellinv(\{e\})$ induce a connected leaper graph of $L$ on $A$; and (b) For every acute angle $b'ab''$ of $M$ on $B$, the sets of cells $\ellinv(\{ab'\})$ and $\ellinv(\{ab''\})$ have at least one element in common.

For claim (a), when $\ellinv(\{e\})$ is empty or it contains a single cell, there is nothing to prove. Otherwise, when $\ellinv(\{e\})$ contains two cells, they are the endpoints of the cross-edge of $e$.

For claim (b), the core of angle $b'ab''$ is in $A$ by Lemma \ref{core}, and it also belongs to both sets $\ellinv(\{ab'\})$ and $\ellinv(\{ab''\})$. \end{proof}

\medskip

There are some natural situations where lifting and lowering are in fact exact inverses of one another for leaper graphs.

\medskip

\begin{lemma} \label{inverse} Let $C$ be a connected component of $L$ on $A$. Suppose that $\ell$ and $C$ satisfy the conditions of the appropriate lemma out of Lemmas \ref{lg}--\ref{lf}. Then $D = \ell(C)$ is an angular component of $M$ on $B$ and $\ellinv(D) = C$. \end{lemma}

\medskip

Explicitly, ``the conditions of the appropriate lemma'' are as follows: If $\ell = f$, then $C$ must be angular; and, if $\ell \in \{g, h\}$, then $C$ must contain at least two cells. Furthermore, as we remarked in Section \ref{lift}, with $\ell \in \{g, h\}$ we can relax these conditions somewhat. For $g$ and $h$, Lemma \ref{inverse} continues to hold when $C$ is a singleton but both of $\CompX(A, C)$ and $\CompY(A, C)$ are non-singletons.

\medskip

\begin{proof} By Lemmas \ref{lg}--\ref{lf}, we get that $D$ is an angular leaper graph of $M$ on $B$. Let $D^\star$ be the angular component of $M$ on $B$ which includes $D$. By Lemma \ref{alc}, we obtain that $C^\star = \ellinv(D^\star)$ is a connected leaper graph of $L$ on $A$.

Since $C$ is a non-singleton, each one of its cells is incident with an edge. By Lemma \ref{cross}, it follows that each cell of $C$ is matched with an edge of $D$. Therefore, all cells of $C$ are also in $C^\star$. Since $C$ is a connected component of $L$ on $A$ and $C^\star$ is an induced connected leaper graph of $L$ on $A$, we conclude that $C$ and $C^\star$ coincide.

We are left to show that $D$ and $D^\star$ coincide as well. This is clear when $D^\star$ consists of a single edge. Otherwise, suppose that $D^\star$ contains at least two edges. Then each edge of $D^\star$ is part of an acute angle in $D^\star$ because $D^\star$ is an angular component of $M$ on $B$. By Lemma \ref{core}, it follows that each edge of $D^\star$ is matched with a cell of $A$. Therefore, all edges of $D^\star$ are also in $\ell(C^\star) = \ell(C) = D$. \end{proof}

%% file: induction-on-descent-10-line.tex
\section{Lineages} \label{line}

Let $\mathfrak{T}$ be a subtree of $\mathfrak{L}$ with root $L$. To each skew free leaper of $\mathfrak{T}$, we assign a board and a leaper graph of that leaper on that board, as follows.

First we choose a board $A$ and a leaper graph $G$ of $L$ on $A$ and we assign them to $L$. Then we determine the rest of the assignments by induction on descent. Whenever $L' \xrightarrow{\ell} L''$ is an arrow of $\mathfrak{T}$ (so that $\ell(L') = L''$) and $L'$ is assigned board $A'$ and leaper graph $G'$, furthermore we assign to $L''$ the board $A'' = \ell(L', A')$ and the leaper graph $G'' = \ell(L', A', G')$.

We call a system of leapers, boards, and leaper graphs of this form a \emph{lineage}. We also call $L$, $A$, and $G$ the \emph{originators} of the lineage.

\medskip

\begin{lemma} \label{offset} Consider an arbitrary lineage. Then there exist two integer constants $w$ and $z$ such that, for every skew free $(r, s)$-leaper $N$ in the lineage, the height $u$ and the width $v$ of the board which the lineage assigns to $N$ satisfy $u = r + s + w$ and $v = r + s + z$. \end{lemma}

\medskip

\begin{proof} By induction on descent. \end{proof}

\medskip

Lemma \ref{offset} explains why in Theorems \ref{thm:rigid} and \ref{thm:journey} we see board sizes of the form $(p + q + w) \times (p + q + z)$ with both of $w$ and $z$ expressed in terms of certain tails of $L$.

We define the \emph{perfect} subtree of $\mathfrak{L}$ with root $L$ to be the one whose vertex set consists of all skew free leapers $N$ such that $L$ is a tail of $N$. Or, equivalently, whose vertex set consists of $L$, $g(L)$ together with all of its descendants, and $h(L)$ together with all of its descendants.

We define the lineage on subtree $\mathfrak{T}$ originated by $L$, $A$, and $G$ to be \emph{perfect} when $\mathfrak{T}$ is perfect, $A$ is standard, and $G$ is a connected component of $L$ on $A$ such that neither one of $\CompX(A, G)$ and $\CompY(A, G)$ is a simple weave.

Thus, in particular, all boards in a perfect lineage are standard.

Note that we allow $G$ to be a singleton, and we also allow each one of its completions $\CompX(A, G)$ and $\CompY(A, G)$ to be a compound weave. The reason for this choice of constraints will become clear in the proof of Theorem \ref{thm:line}.

We consider perfect lineages originated by the wazir, too, even though it is not a skew free leaper.

Observe first that, in order to specify a perfect lineage originated by the wazir, it suffices to indicate its originating board since that would determine its originating wazir graph as well.

Formally, $f$ preserves the wazir whereas both of $g$ and $h$ map it onto the knight. Furthermore, with the wazir, both of $g$ and $h$ share the same lateral, diagonal, and margin parameters, and so they act upon all boards and wazir graphs in the same way.

This line of reasoning eventually leads us to the following definition: The perfect lineage originated by the wazir with the board $A$ of size $m \times n$ assigns to all skew free leapers the same boards and leaper graphs as the non-perfect lineage on $\mathfrak{L}$ originated by the knight with the board $A^+ = g(\text{Wazir}, A) = h(\text{Wazir}, A)$ of size $(m + 2) \times (n + 2)$ and the knight graph formed by all knight edges on $A^+$.

Perfect lineages occur naturally in the three settings we considered in the introduction. For directional rigidity and wazir journeys, we discuss this in detail in Sections \ref{rigid} and \ref{journey}, respectively, in our proofs of sufficiency for Theorems \ref{thm:rigid} and \ref{thm:journey}. Here, we take a look at the second leaper theorem through the same lens.

Observe that every clover is an angular leaper graph. Lifting a clover yields another clover and lowering a clover of a skew free leaper distinct from the knight yields either another clover or a singleton. Thus we can group all clovers into perfect lineages, as follows.

We define the \emph{\textbf{SL}-lineage} of the wazir to be the perfect lineage originated by the wazir with the board of size $1 \times 1$.

We also define an \emph{\textbf{SL}-lineage} of the skew free leaper $L$ to be a perfect lineage originated by $L$ with the square board of side $p + q$ and a singleton connected component of $L$ on that board. Thus each skew free leaper $L$ originates $(q - p)^2$ distinct \textbf{SL}-lineages.

This ensures that each clover belongs to exactly one \textbf{SL}-lineage.

We can now restate the key idea of the proof of the second leaper theorem in terms of lineages: Let $N$ be a free leaper and let $C$ be a clover in the \textbf{SL}-lineage of $N$. Then there exists a Hamiltonian cycle of $N$ on the cells of $C$.

The main result of this section is as follows.

\medskip

\begin{theorem} \label{thm:line} Consider an arbitrary perfect lineage. Let $L$ be a skew free leaper in it and let $A$ and $G$ be the board and the leaper graph which the lineage assigns to $L$. Then $G$ is a connected component of $L$ on $A$. Furthermore, let $L \xrightarrow{\ell} M$ be an arrow in that perfect lineage and let $B$ and $H$ be the board and the leaper graph which the lineage assigns to $M$. Then $\ellinv(H) = G$. \end{theorem}

\medskip

Thus, in particular, within a perfect lineage lifting and lowering are exact inverses of one another.

Note that we do not rely on Theorem \ref{thm:line} in our proof of Theorems \ref{thm:rigid} and \ref{thm:journey}. However, it is clearly of interest in its own right.

Before we can approach Theorem \ref{thm:line}, we must establish some lemmas. The first one of them, Lemma \ref{wf}, plays a crucial role also in the proofs of Theorems \ref{thm:rigid} and \ref{thm:journey}. The other two, Lemmas \ref{wg} and \ref{wh}, do not appear elsewhere in this work.

Define $\ell$, $\widehat{p}$, $\widehat{q}$, $M$, $\widehat{m}$, $\widehat{n}$, and $B$ as in Sections \ref{lift} and \ref{low}.

\medskip

\begin{lemma} \label{wf} Suppose that $\ell = f$. Consider a matched pair of a cell $a$ of $A$ and an edge $e$ of $M$ on $B$. Then $\CompX(A, a)$ is a weave if and only if $\CompX(B, e)$ is a weave. \end{lemma} 

\medskip

\begin{proof} Observe that $n \ge p + q$ if and only if $\widehat{n} \ge \widehat{p} + \widehat{q}$. When both of these inequalities hold, by Lemma \ref{wb} neither $\Pi_X$ nor $\Pi(\widehat{p}, \widehat{q}, J_X)$ contains a weave. Suppose, throughout the rest of the proof, that $n < p + q$ and $\widehat{n} < \widehat{p} + \widehat{q}$.

Then by Lemma \ref{pa} both of $\CompX(A, a)$ and $\CompX(B, e)$ are paths. Let \[\CompX(A, a) = x_{0, 0}x_{0, 1} \ldots x_{0, s_0}x_{1, 0}x_{1, 1} \ldots x_{1, s_1} \ldots x_{d, 0}x_{d, 1} \ldots x_{d, s_d}\] so that each $x_{i, j} \to x_{i, j + 1}$ is a short move to the right and each $x_{i, s_i} \to x_{i + 1, 0}$ is a long move to the left.

Let $x_{i, -1} = x_{i, 0} - p$ and $x_{i, s_i + 1} = x_{i, s_i} + p$ for all $i$, and consider the path \begin{gather*} x_{0, -1}x_{0, 0}x_{0, 1} \ldots x_{0, s_0}x_{0, s_0 + 1}x_{1, -1}x_{1, 0}x_{1, 1} \ldots\\ x_{1, s_1}x_{1, s_1 + 1} \ldots x_{d, -1}x_{d, 0}x_{d, 1} \ldots x_{d, s_d}x_{d, s_d + 1} \end{gather*} in $\Pi(\widehat{p}, \widehat{q}, J_X)$.

Since both endpoints of this path are of degree one in $\Pi(\widehat{p}, \widehat{q}, J_X)$, by Lemma \ref{pa} we get that it is a connected component of $\Pi(\widehat{p}, \widehat{q}, J_X)$. Furthermore, since it contains the $x$-projection of $e$, we conclude that in fact our path coincides with $\CompX(B, e)$.

Thus $\CompX(A, a)$ is a weave if and only if $s_i$ is even for all $i$ with $0 < i < d$ and $\CompX(B, e)$ is weave if and only if $s_i + 2$ is even for all $i$ with $0 < i < d$. \end{proof}

\medskip

\begin{lemma} \label{wg} Suppose that $\ell = g$. Consider a matched pair of a cell $a$ of $A$ and an edge $e$ of $M$ on $B$. Then $\CompX(A, a)$ is a simple weave if and only if $\CompX(B, e)$ is a weave. \end{lemma} 

\medskip

\begin{proof} We begin as in the proof of Lemma \ref{wf}. However, we construct a path in $\Pi(\widehat{p}, \widehat{q}, J_X)$ out of $\CompX(A, a)$ somewhat differently.

For all $i$ and $j$, let $x^-_{i, j} = x_{i, j} - q$ and $x^+_{i, j} = x_{i, j} + q$. First we replace each vertex $x_{i, j}$ in $\CompX(A, a)$ with the fragment $x^-_{i, j}x_{i, j}x^+_{i, j}$. Then in the sequence of vertices thus obtained we furthermore delete the first one and the last one as well as all fragments of the form $x_{i, s_i}x^+_{i, s_i}x^-_{i + 1, 0}x_{i + 1, 0}$. The net result is the path \begin{gather*} x_{0, 0}x^+_{0, 0}x^-_{0, 1}x_{0, 1}x^+_{0, 1} \ldots x^-_{0, s_0 - 1}x_{0, s_0 - 1}x^+_{0, s_0 - 1}x^-_{0, s_0}x^+_{1, 0}x^-_{1, 1}x_{1, 1}x^+_{1, 1} \ldots\\ x^-_{1, s_1 - 1}x_{1, s_1 - 1}x^+_{1, s_1 - 1}x^-_{1, s_1} \ldots x^+_{d, 0}x^-_{d, 1}x_{d, 1}x^+_{d, 1} \ldots x^-_{d, s_d - 1}x_{d, s_d - 1}x^+_{d, s_d - 1}x^-_{d, s_d}x_{d, s_d} \end{gather*} in $\Pi(\widehat{p}, \widehat{q}, J_X)$.

Just as in the proof of Lemma \ref{wf}, this path coincides with $\CompX(B, e)$.

Consider an arbitrary subpath of $\CompX(B, e)$ with signature of the form $\mathtt{l}\mathtt{s}^t\mathtt{l}$. Then $t = 2$ unless that subpath is of the form $x^+_{i, s_i - 1}x^-_{i, s_i}x^+_{i + 1, 0}x^-_{i + 1, 1}$, when $t = 1$. Therefore, $\CompX(B, e)$ is a weave if and only if all edges of $\CompX(A, a)$ are short. \end{proof}

\medskip

\begin{lemma} \label{wh} Suppose that $\ell = h$. Consider a matched pair of a cell $a$ of $A$ and an edge $e$ of $M$ on $B$. Then $\CompX(A, a)$ is a simple weave if and only if $\CompX(B, e)$ is a weave. \end{lemma} 

\medskip

\begin{proof} Once again, we begin as in the proof of Lemma \ref{wf}. This time around, though, we construct a path in $\Pi(\widehat{p}, \widehat{q}, J_X)$ out of $\CompX(A, a)$ as follows.

Define $x^-_{i, j}$ and $x^+_{i, j}$ as in the proof of Lemma \ref{wg}. First we replace each vertex $x_{i, j}$ in $\CompX(A, a)$ with the fragment $x^+_{i, j}x_{i, j}x^-_{i, j}$. Then in the sequence of vertices thus obtained we furthermore delete all fragments of the form $x^-_{i, s_i}x^+_{i + 1, 0}$. The net result is the path \begin{gather*} x^+_{0, 0}x_{0, 0}x^-_{0, 0} \ldots x^+_{0, s_0 - 1}x_{0, s_0 - 1}x^-_{0, s_0 - 1}x^+_{0, s_0}x_{0, s_0}x_{1, 0}x^-_{1, 0}x^+_{1, 1}x_{1, 1}x^-_{1, 1} \ldots\\ x^+_{1, s_1 - 1}x_{1, s_1 - 1}x^-_{1, s_1 - 1}x^+_{1, s_1}x_{1, s_1} \ldots x_{d, 0}x^-_{d, 0}x^+_{d, 1}x_{d, 1}x^-_{d, 1} \ldots x^+_{d, s_d}x_{d, s_d}x^-_{d, s_d} \end{gather*} in $\Pi(\widehat{p}, \widehat{q}, J_X)$.

Just as in the proof of Lemma \ref{wf}, this path coincides with $\CompX(B, e)$.

Consider an arbitrary subpath of $\CompX(B, e)$ with signature of the form $\mathtt{l}\mathtt{s}^t\mathtt{l}$. Then $t = 2$ unless that subpath is of the form $x^-_{i, s_i - 1}x^+_{i, s_i}x_{i, s_i}x_{i + 1, 0}x^-_{i + 1, 0}x^+_{i + 1, 1}$, when $t = 3$. Therefore, $\CompX(B, e)$ is a weave if and only if all edges of $\CompX(A, a)$ are short. \end{proof}

\medskip

\begin{proof*}{Proof of Theorem \ref{thm:line}} We will show additionally that $G$ is a non-weave, except possibly when it is the originating leaper graph of our lineage.

For convenience, we assume that the lineage is originated by a skew free leaper. The necessary modifications in the case when it is originated by the wazir are few and straightforward.

We proceed by induction on descent.

In the base case, $G$ is a connected component of $L$ on $A$ by the definition of a perfect lineage.

For the induction step, suppose that $G$ is a connected component of $L$ on $A$ and also that either $G$ is a non-weave or $G$ is the lineage's originating leaper graph.

When $\ell = f$, we cannot be at the root of our lineage, and so $G$ cannot be a weave. Thus by Lemma \ref{naw} we get that $G$ is angular.

Otherwise, when $\ell \in \{g, h\}$, since neither one of the completions $\CompX(A, G)$ and $\CompY(A, G)$ is a simple weave, neither one of them can be a singleton, either.

Therefore, in all cases it follows by Lemma \ref{inverse} that $H = \ell(G)$ is an angular component of $M$ on $B$ and $G = \ellinv(H)$. (Note that we use the stronger version of Lemma \ref{inverse} with the weaker conditions on $G$, as outlined in Section \ref{low}.)

Lastly, let $H^\star$ be the connected component of $M$ on $B$ which includes $H$. By Lemmas \ref{wf}--\ref{wh}, we obtain that $H^\star$ cannot be a weave. Thus by Lemma \ref{naw} it follows that $H^\star$ is an angular connected component of $M$ on $B$. Or, in other words, $H$ and $H^\star$ coincide. \end{proof*}

%% file: induction-on-descent-11-conn.tex
\section{Connectedness} \label{conn}

To begin with, we apply our methods to Theorems \ref{thm:conn} and \ref{thm:econn}.

\medskip

\begin{proof*}{Proof of necessity for Theorem \ref{thm:conn}} Suppose that $A^\mathcal{C} \not \sqsubseteq A$, and so either $\min\{m, n\} < p + q$ or $\max\{m, n\} < 2q$. Suppose also, for concreteness, that $m \le n$.

When $m < p + q$, the complete leaper graph of $L$ on $A$ is disconnected by Lemmas \ref{pa} and \ref{pp}. On the other hand, when $n < 2q$, the complete leaper graph of $L$ on $A$ is disconnected because it contains an isolated vertex. \end{proof*}

\medskip

Our argument in the case of $m < p + q$ is not substantially different from the ones in \cite{RW} and \cite{K}. However, our phrasing in terms of projections perhaps makes the central idea somewhat clearer.

\medskip

\begin{proof*}{Proof of sufficiency for Theorem \ref{thm:conn}} We proceed by induction on descent.

Our base case is the knight, and the knight graph on the board of size $3 \times 4$ is indeed connected.

For the induction step, let $\ell$ be a lifting transformation with $\ell(L) = M$ and suppose that the sufficiency part of Theorem \ref{thm:conn} holds for $L$.

Define the board $B^\mathcal{C}$ relative to $M$ in the same way as we defined the board $A^\mathcal{C}$ relative to $L$, so that $B^\mathcal{C}$ is the board of size $(\widehat{p} + \widehat{q}) \times 2\widehat{q}$.

Let also $A^\star = \ellinv(B^\mathcal{C})$, so that $A^\star$ is the board of size $(p + q) \times 2\max\{\Lat, \Dia\}$.

Consider the complete leaper graph $G^\star$ of $L$ on $A^\star$. By the induction hypothesis and the monotonicity of $\mathcal{C}_L$, we get that $G^\star$ is connected. Furthermore, by Lemmas \ref{wb} and \ref{naw} we get that $G^\star$ is angular, too. Thus by Lemmas \ref{lg}--\ref{lf} we conclude that the leaper graph $H^\star = \ell(G^\star)$ of $M$ on $B^\mathcal{C}$ is connected as well.

On the other hand, it is straightforward to verify that each cell of $B^\mathcal{C}$ is incident with an edge of $M$ on $B^\mathcal{C}$ and each edge of $M$ on $B^\mathcal{C}$ is part of an acute angle of $M$ on $B^\mathcal{C}$. Thus by Lemma \ref{core} in fact $H^\star$ is the complete leaper graph of $M$ on $B^\mathcal{C}$. \end{proof*}

\medskip

Theorem \ref{thm:econn} requires more work.

\medskip

\begin{lemma} \label{eca} Suppose that $\mathcal{E}_L(m, n)$. Then $m \ge p + q$ and $n \ge p + q$. \end{lemma} 

\medskip

\begin{proof} Our proof of monotonicity for $\mathcal{E}_L$ shows that $\mathcal{E}_L(m, n)$ implies $m \ge q + 1$ and $n \ge q + 1$. Then each one of the four corner cells of $A$ and their neighbours by side is incident with an edge of $L$ on $A$. By Lemmas \ref{pa} and \ref{pp}, the sides of $A$ must satisfy $m \ge p + q$ and $n \ge p + q$ for all of these edges to be in the same connected component of $L$ on $A$. \end{proof}

\medskip

\begin{proof*}{Proof of Theorem \ref{thm:econn}} We proceed by induction on descent.

Our base case is the knight. Since $\mathcal{E}_\text{Knight}$ is symmetric and monotone, it suffices to check that it does not hold on boards of height at most two but it does hold on the board of size $3 \times 3$. Both of these claims are straightforward.

For the induction step, let $\ell$ be a lifting transformation with $\ell(L) = M$ and suppose that Theorem \ref{thm:econn} holds for $L$.

Define the board $B^\mathcal{E}$ relative to $M$ in the same way as we defined the board $A^\mathcal{E}$ relative to $L$.

Consider a board $B$ of size $\widehat{m} \times \widehat{n}$. By Lemma \ref{eca}, if either $\widehat{m} < \widehat{p} + \widehat{q}$ or $\widehat{n} < \widehat{p} + \widehat{q}$, then Theorem \ref{thm:econn} holds for $M$ on $B$. Thus suppose, throughout the rest of the proof, that $\widehat{m} \ge \widehat{p} + \widehat{q}$ and $\widehat{n} \ge \widehat{p} + \widehat{q}$.

Then $\ellinv(B)$ is well-defined. Let $A = \ellinv(B)$. Observe that, consequently, $m \ge p + q$ and $n \ge p + q$ as well.

From this point on, we consider the cases of $\ell = f$ and $\ell \in \{f, g\}$ separately.

\smallskip

\emph{Case 1}. $\ell = f$. Let $S$ be the set of all cells of $A$ incident with an edge of $L$ on $A$ and let $T$ be the set of all edges of $M$ on $B$. It is straightforward to verify that $m \ge p + q$ and $n \ge p + q$ imply $f(S) = T$ and $\finv(T) = S$.

Suppose first that all cells of $S$ are in the same connected component $C$ of $L$ on $A$. Then $C$ is angular by Lemmas \ref{wb} and \ref{naw}. Consequently, all edges of $T$ are in the same connected component of $M$ on $B$ by Lemma \ref{lf}.

Conversely, suppose that all edges of $T$ are in the same connected component $D$ of $M$ on $B$. Then $D$ is angular by Lemmas \ref{wb} and \ref{naw}. Consequently, all cells of $S$ are in the same connected component of $L$ on $A$ by Lemma \ref{alc}.

Or, in summary, $\mathcal{E}_L(A)$ if and only if $\mathcal{E}_M(B)$. By the induction hypothesis, it follows that $\mathcal{E}_M(B)$ if and only if $B^\mathcal{E} = f(A^\mathcal{E}) \sqsubseteq f(A) = B$.

\smallskip

\emph{Case 2}. $\ell \in \{g, h\}$. Let $S$ be the set of all cells of $A$ and let $T$ be the set of all edges of $M$ on $B$. Once again, $m \ge p + q$ and $n \ge p + q$ imply $\ell(S) = T$ and $\ellinv(T) = S$.

We establish just as in Case 1 that all cells of $S$ are in the same connected component of $L$ on $A$ if and only if all edges of $T$ are in the same connected component of $M$ on $B$. The only difference is that we must refer to Lemmas \ref{lg} and \ref{lh} in place of Lemma \ref{lf}.

Thus $\mathcal{C}_L(A)$ if and only if $\mathcal{E}_M(B)$. By Theorem \ref{thm:conn}, it follows that $\mathcal{E}_M(B)$ if and only if $B^\mathcal{E} = \ell(A^\mathcal{C}) \sqsubseteq \ell(A) = B$. \end{proof*}

\medskip

We use Theorem \ref{thm:conn} in Case 2 of the proof. But we can also make do without it, as follows.

The cells of $A$ are all in the same connected component of $L$ if and only if no cell of $A$ is an isolated vertex in the complete leaper graph of $L$ on $A$ and the edges of $L$ on $A$ form a connected leaper graph. The former occurs if and only if $A^\mathcal{I} \sqsubseteq A$ by Theorem \ref{thm:i}. The latter occurs if and only if $A^\mathcal{E} \sqsubseteq A$ by Theorem \ref{thm:econn} for $L$. (This is our induction hypothesis.) Finally, both of $A^\mathcal{I}$ and $A^\mathcal{E}$ fit inside of $A$ if and only if $A^\mathcal{C}$ fits inside of $A$.

Thus Theorem \ref{thm:econn} admits a self-contained proof which does not reference Theorem \ref{thm:conn}.

The exact same argument shows also that in fact Theorems \ref{thm:i} and \ref{thm:econn} together imply Theorem \ref{thm:conn}.

Observe that our proof of necessity for Theorem \ref{thm:conn} is direct and does not use induction on descent. The necessity part of Theorem \ref{thm:econn} admits a direct proof as well, as follows.

\medskip

\begin{proof*}{Direct proof of necessity for Theorem \ref{thm:econn}} Let $\xi$ be the width of $A^\mathcal{E}$. Since $\mathcal{E}_L$ is symmetric and monotone, and in light of Lemma \ref{eca}, it suffices to demonstrate that if $p \ge 2$ then the edges of $L$ on the standard board $A_\maltese$ of size $(\xi - 1) \times (\xi - 1)$ do not form a connected leaper graph.

Let $q = 2cp + \varepsilon r$ as in Section \ref{weave}. Then also let $S$ be the set of all cells of $A_\maltese$ of the form either $(cp, ip)$ with $c + i$ even or $(cp + \varepsilon r, ip)$ with $c + i$ odd and let $S_\maltese$ be the union of the reflections of $S$ with respect to the four axes of symmetry of $A_\maltese$.

\begin{figure}[t!] \centering \includegraphics{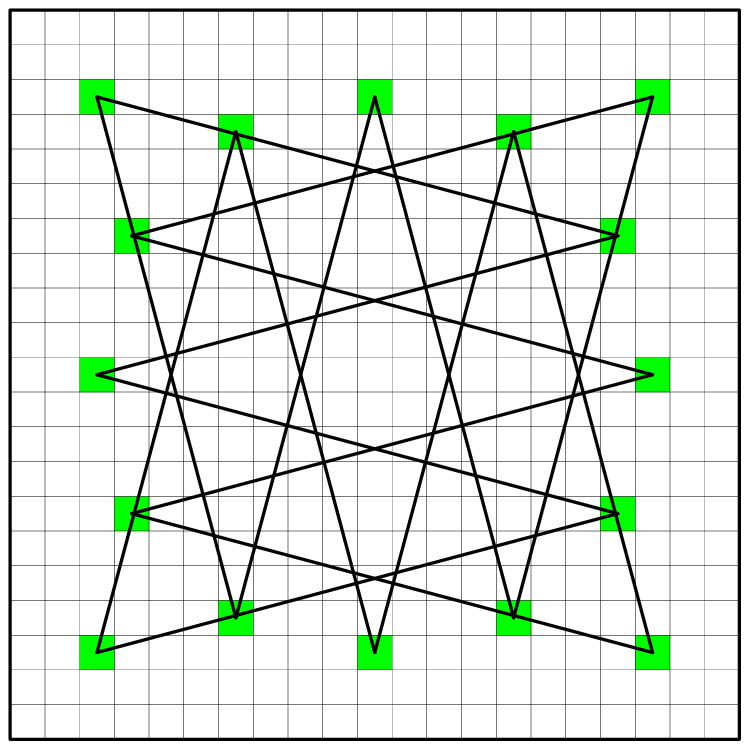} \caption{} \label{e-conn} \end{figure} 

It is straightforward to verify that the cells of $S_\maltese$ are the vertices of a connected component of $L$ on $A_\maltese$ which does not contain all edges of $L$ on $A_\maltese$. (In fact, this connected component is also a clover of $L$ on the standard board of size $(p + q) \times (p + q)$.) \end{proof*}

\medskip

For example, the complete leaper graph of the $(4, 15)$-leaper on the board of size $21 \times 21$ contains exactly two non-singleton connected components. One of them, as described in the proof, is shown in Figure \ref{e-conn}.

%% file: induction-on-descent-12-rigid.tex
\section{Directional Rigidity} \label{rigid}

We are ready to approach Theorem \ref{thm:rigid}.

Throughout this section, we will rely on the criterion of Lemma \ref{ru}: A finite leaper graph of $L$ is directionally rigid if and only if it contains an unbalanced closed walk.

Observe also that every unbalanced closed walk of $L$ must necessarily contain edges of both inclines.

First we study the directionally rigid weaves of skew free leapers.

\medskip

\begin{lemma} \label{rw} Suppose that there exists a directionally rigid horizontal weave $U$ of $L$ on $A$. Then $p \ge 2$, $m \ge q + 1$, and $n \ge p + 2q - (q \bmod 2)$. \end{lemma} 

\medskip

\begin{proof} Both of $p \ge 2$ and $m \ge q + 1$ follow because $U$ contains edges of both inclines.

Let $\alpha$ be an unbalanced closed walk in $U$. Then at least one projection of $\alpha$ is unbalanced as well. Since $\CompY(U)$ is acyclic by Lemma \ref{wb}, we conclude that it must be the $x$-projection $\chi$ of $\alpha$ which is the unbalanced one.

Since the $y$-projection of $\alpha$ is contained within the projection weave $\CompY(U)$, it traverses an even number of short edges between each pair of long ones. Conversely, $\chi$ must traverse an even number of long edges between each pair of short ones.

Let us replace each maximal subwalk of long moves in $\chi$ with a path of long moves in $\Pi(p, 2q, I_X)$. This transforms $\chi$ into an unbalanced closed walk $\chi^\star$ in $\Pi(p, 2q, I_X)$. By Lemma \ref{pu}, it follows that $n = |I_X| \ge p + 2q - \gcd(p, 2q) + 1 = p + 2q - (q \bmod 2)$. \end{proof}

\medskip

\begin{lemma} \label{ro} There exists a directionally rigid connected component of $L$ on the board of size $(q + 1) \times (p + 2q - [q \bmod 2])$. \end{lemma} 

\medskip

Observe also that the connected component of Lemma \ref{ro} is a horizontal weave if and only if $p \ge 2$.

\medskip

\begin{proof} Let $m = q + 1$ and $n = p + 2q - (q \bmod 2)$.

By Lemma \ref{pu}, there exists an unbalanced closed walk $\chi^\star$ in $\Pi(p, 2q, I_X)$. By replacing each long move in it with a path of two long moves in $\Pi_X$, we obtain an unbalanced closed walk $\chi$ in $\Pi_X$ with signature $\mathtt{l}^{d_0}\mathtt{s}\mathtt{l}^{d_1}\mathtt{s}\mathtt{l}^{d_2} \ldots \mathtt{s}\mathtt{l}^{d_k}$ such that $d_i$ is an even nonnegative integer for all $i$.

By going around $\chi$ twice if needed, we can ensure without loss of generality that $k$ is even as well.

Consider the subgraph $G = (y_\text{Min} + p)$---$y_\text{Min}$---$y_\text{Max}$---$(y_\text{Max} - p)$ of $\Pi_Y$. Since $k$ is even and also $d_i$ is even for all $i$, we can find a closed walk $\upsilon$ in $G$ with signature $\mathtt{s}^{d_0}\mathtt{l}\mathtt{s}^{d_1}\mathtt{l}\mathtt{s}^{d_2} \ldots \mathtt{l}\mathtt{s}^{d_k}$. Then by Lemma \ref{ps} there exists a closed walk $\alpha$ of $L$ on $A$ whose $x$-projection is $\chi$ and whose $y$-projection is $\upsilon$. \end{proof}

\medskip

\begin{figure}[t!] \centering \includegraphics{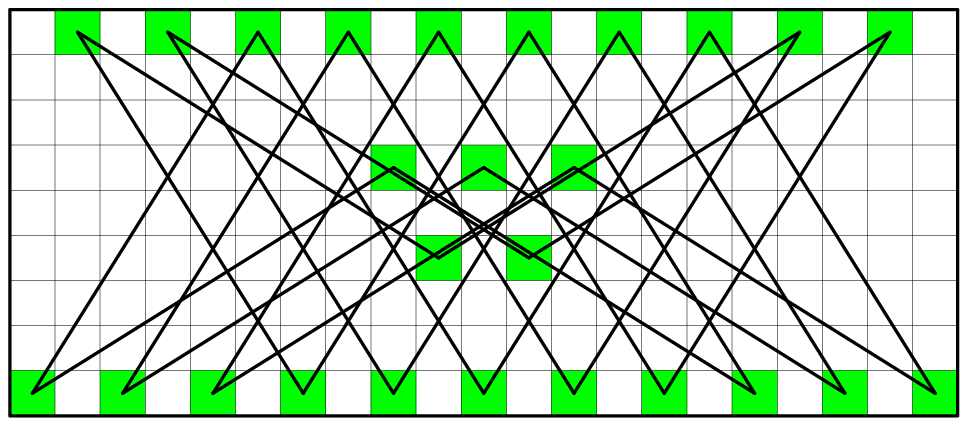} \caption{} \label{uc} \end{figure}

For example, Figure \ref{uc} shows one unbalanced closed walk of the $(5, 8)$-leaper on the board of size $9 \times 21$ constructed as in the proof.

We proceed to study the effect of lifting and lowering on closed walks.

Let $\ell$ be a lifting transformation with $\ell(L) = M$, let $\alpha$ be a closed walk of $L$, and let $\beta$ be a closed walk of $M$.

We say that $\alpha$ and $\beta$ are \emph{matched} when $\alpha = a_1a_2 \ldots a_ka_1$ and $\beta = b_{1, 1}b_{1, 2} \ldots\allowbreak b_{1, j_1}b_{2, 1}b_{2, 2} \ldots b_{2, j_2} \ldots b_{k, 1}b_{k, 2} \ldots b_{k, j_k}b_{1, 1}$ so that: (a) For all $i$ and all $j$ with $1 \le j < j_i$, cell $a_i$ and edge $b_{i, j}b_{i, j + 1}$ are matched; and (b) For all $i$, edge $a_ia_{i + 1}$ and edge $b_{i, j_i}b_{i + 1, 1}$ are cross-edges. We allow the indices to wrap around cyclically when $i = k$ in condition (b), so that $a_{k + 1}$ denotes the same cell as $a_1$ and $b_{k + 1, 1}$ denotes the same cell as $b_{1, 1}$.

Lemmas \ref{wa}--\ref{wlow} below recast some of the results of Sections \ref{ang}, \ref{lift}, and \ref{low} from the point of view of the closed walks of skew free leapers.

\medskip

\begin{lemma} \label{wa} Suppose that $p \ge 2$. Let $\alpha$ be a closed walk of $L$ on $A$ such that the connected component of $L$ on $A$ which contains $\alpha$ is a non-weave. Then there exists an angular-closed walk $\alpha^\star$ of $L$ on $A$ with $\Phi(\alpha) = \Phi(\alpha^\star)$. \end{lemma} 

\medskip

\begin{proof} Let $ab'$ and $ab''$ be two consecutive edges of $\alpha$. Note that we consider the first edge and the last edge of $\alpha$ to be consecutive as well.

When angle $b'ab''$ is either right or diagonally obtuse, we can insert some detours at $a$ as in the proof of Lemma \ref{rdo} so as to replace it with a series of angles each one of which is either acute or zero.

Otherwise, when angle $b'ab''$ is either laterally obtuse or straight, we construct a closed walk $\gamma = ac_1c_2 \ldots c_sa$ of $L$ on $A$ as in the proof of Lemma \ref{los} so that each angle formed by two consecutive edges of the walk $b'ac_1c_2 \ldots c_sab''$ is either acute, zero, right, or diagonally obtuse. Then we splice $\gamma$ into $\alpha$ at $a$ and we rework all right and diagonally obtuse angles of the longer closed walk thus obtained by inserting some additional detours into it as in the previous case.

It is straightforward to see that the construction in the proof of Lemma \ref{los} yields a balanced $\gamma$. Therefore, all of our augmentations preserve $\Phi(\alpha)$. \end{proof}

\medskip

\begin{lemma} \label{mwf} Let $\alpha$ be an angular-closed walk of $L$ on $A$. Then there exists an angular-closed walk $\beta$ of $M$ on $f(A)$ such that $\alpha$ and $\beta$ are matched. \end{lemma} 

\medskip

\begin{proof} This is a direct corollary of the proof of Lemma \ref{ace}. \end{proof}

\medskip

\begin{lemma} \label{mwgh} Let $\alpha$ be a closed walk of $L$ on $A$. Suppose that $\ell \in \{g, h\}$. Then there exists an angular-closed walk $\beta$ of $M$ on $\ell(A)$ such that $\alpha$ and $\beta$ are matched. \end{lemma} 

\medskip

\begin{proof} This is a direct corollary of the proofs of Lemmas \ref{lg} and \ref{lh}. \end{proof}

\medskip

\begin{lemma} \label{wlow} Let $\beta$ be an angular-closed walk of $M$ on $\ell(A)$ such that at least two distinct cells of $A$ are cores of acute angles in $\beta$. Then there exist a closed walk $\alpha$ of $L$ on $A$ and a cyclic shift $\beta^\star$ of $\beta$ such that $\alpha$ and $\beta^\star$ are matched. \end{lemma} 

\medskip

Observe also that, if at most one cell of $A$ is the core of an acute angle in $\beta$, then $\beta$ is balanced.

\medskip

\begin{proof} This is a direct corollary of the proof of Lemma \ref{alc}. \end{proof}

\medskip

Let $\mathcal{M}_f = \left(\begin{smallmatrix} 1 & 0\\ -2 & 1 \end{smallmatrix}\right)$, $\mathcal{M}_g = \left(\begin{smallmatrix} 2 & -1\\ 1 & 0 \end{smallmatrix}\right)$, and $\mathcal{M}_h = \left(\begin{smallmatrix} -2 & 1\\ 1 & 0 \end{smallmatrix}\right)$.

\medskip

\begin{lemma} \label{mwm} Let $\alpha$ be a closed walk of $L$ and let $\beta$ be a closed walk of $M$ so that $\alpha$ and $\beta$ are matched. Then $\Phi(\beta) = \Phi(\alpha)\mathcal{M}_\ell$. Thus, in particular, $\alpha$ is balanced if and only if $\beta$ is balanced. \end{lemma} 

\medskip

\begin{proof} Fix a reference cell $c$ and let $C$ be the induced leaper graph of $M$ on vertices $c^\text{E}$, $c^\text{NE}$, \ldots, $c^\text{SE}$. Thus $C$ is a cycle of length eight.

Let $\alpha = a_1a_2 \ldots a_ka_1$ and $\beta = b_{1, 1}b_{1, 2} \ldots b_{1, j_1}b_{2, 1}b_{2, 2} \ldots b_{2, j_2} \ldots b_{k, 1}b_{k, 2} \ldots b_{k, j_k}\allowbreak b_{1, 1}$ as in the definition of matched closed walks.

For each $i$, let $T_i$ be the translation which maps $a_i$ onto $c$ and let also $T_i$ map the subwalk $b_{i, 1}b_{i, 2} \ldots b_{i, j_i}b_{i + 1, 1}$ of $\beta$ onto the walk $\beta_i = b^\star_{i, 1}b^\star_{i, 2} \ldots b^\star_{i, j_i}b^\star_{i, j_i + 1}$ in $C$. Let also $\gamma_i$ be the unique path in $C$ from $b^\star_{i, j_i + 1}$ to $b^\star_{i + 1, 1}$ which goes counterclockwise around $c$. We allow the indices to wrap around cyclically, as in the definition of matched closed walks, when $i = k$.

Then the concatenation $\gamma$ of $\beta_1$, $\gamma_1$, $\beta_2$, $\gamma_2$, \ldots, $\beta_k$, $\gamma_k$ is a closed walk in $C$.

Since $C$ is directionally flexible, $\sum_i \Phi(\beta_i) + \sum_i \Phi(\gamma_i) = \Phi(\gamma) = \mathbf{0}$. Consequently, $\Phi(\beta) = \sum_i \Phi(\beta_i) = -\sum_i \Phi(\gamma_i)$.

On the other hand, by direct calculation we see that $\Phi(\gamma_i) = -\Phi(a_i \to a_{i + 1})\mathcal{M}_\ell$ in all cases for $\ell$, the direction of move $a_i \to a_{i + 1}$, and the direction of move $b_{i, j_i} \to b_{i + 1, 1}$. Therefore, $\Phi(\beta) = \sum_i \big[\Phi(a_i \to a_{i + 1})\mathcal{M}_\ell\big] = \big[\sum_i \Phi(a_i \to a_{i + 1})\big]\mathcal{M}_\ell = \Phi(\alpha)\mathcal{M}_\ell$. \end{proof}

\medskip

We are only left to put the pieces together.

\medskip

\begin{proof*}{Proof of sufficiency for Theorem \ref{thm:rigid}} We define the \emph{$\mathcal{R}$-lineage} of the wazir to be the perfect lineage originated by the wazir with the board of size $2 \times 1$. Observe that the knight graph in this lineage contains an unbalanced angular-closed walk.

We also define an \emph{$\mathcal{R}$-lineage} of the skew free leaper $L$ to be a perfect lineage originated by $L$ with a board and a directionally rigid connected component of $L$ on that board as in Lemma \ref{ro}.

By Lemma \ref{offset}, the $\mathcal{R}$-lineages which contain $L$ assign to it precisely the boards $A^\mathcal{R}_i$.

On the other hand, by induction on descent with the help of Lemmas \ref{mwf}, \ref{mwgh}, and \ref{mwm}, it follows that every leaper graph of a skew free leaper in an $\mathcal{R}$-lineage contains an unbalanced closed walk of that leaper and is therefore indeed directionally rigid.~\end{proof*}

\medskip

It is straightforward to see that the number of directionally rigid connected components of $L$ on the board of Lemma \ref{ro} is one when either $p = 1$ or $p$ is even and two otherwise, when $p \ge 3$ and $p$ is odd. This is also the number of $\mathcal{R}$-lineages which $L$ originates. Furthermore, in the latter case reflection with respect to the horizontal axis of symmetry of the board swaps the two connected components whereas reflection with respect to the vertical axis of symmetry of the board preserves each one of them.

\medskip

\begin{proof*}{Proof of necessity for Theorem \ref{thm:rigid}} We establish one slightly stronger result: Suppose that there exists a directionally rigid connected component of $L$ on $A$. When that connected component is a weave, $A^\mathcal{R}_\kappa \sqsubseteq A$. Otherwise, when that connected component is a non-weave, $A^\mathcal{R}_i \sqsubseteq A$ for some $i$ with $i < \kappa$.

We proceed by induction on descent.

Our base case is the knight, and with it the verification is straightforward.

For the induction step, let $\ell$ be a lifting transformation with $\ell(L) = M$ and suppose that our strengthening of the necessity part of Theorem \ref{thm:rigid} holds for $L$.

Define the boards $B^\mathcal{R}_i$ relative to $M$ in the same way as we defined the boards $A^\mathcal{R}_i$ relative to $L$. When $\ell = f$, we get that $f(A^\mathcal{R}_i) = B^\mathcal{R}_i$ if and only if $i < \kappa$. Otherwise, when $\ell \in \{g, h\}$, we get that $\ell(A^\mathcal{R}_i) = B^\mathcal{R}_i$ for all $i$ with $0 \le i \le \kappa$.

Let $D$ be a directionally rigid connected component of $M$ on a board $B$ and let $\beta$ be an unbalanced closed walk in $D$.

When $D$ is a weave, we are done immediately by Lemma \ref{rw}.

Suppose, throughout the rest of the proof, that $D$ is a non-weave.

Observe that $\ellinv(B)$ is well-defined because an unbalanced closed walk necessarily contains edges of both inclines. Let $\ellinv(B) = A$.

If $\widehat{p} = 1$, then since the clover of $M$ on the board of size $(\widehat{q} + 1) \times (\widehat{q} + 1)$ is directionally flexible it follows that $B^\mathcal{R}_0 \sqsubseteq B$.

Suppose, otherwise, that $\widehat{p} \ge 2$. By Lemma \ref{wa}, we can assume without loss of generality that $\beta$ is closed-angular. Consequently, by Lemmas \ref{wlow} and \ref{mwm} there exists an unbalanced closed walk $\alpha$ of $L$ on $A$ matched with some cyclic shift $\beta^\star$ of $\beta$. Let $C$ be the connected component of $L$ on $A$ which contains $\alpha$.

We consider the case when $\ell = f$ first. Then $C$ is a non-weave by Lemma \ref{wf}, and so $A^\mathcal{R}_i \sqsubseteq A$ for some $i$ with $i < \kappa$ by the induction hypothesis. Therefore, $B^\mathcal{R}_i = f(A^\mathcal{R}_i) \sqsubseteq f(A) = B$, as needed.

We continue with the case when $\ell \in \{g, h\}$. Then $A^\mathcal{R}_i \sqsubseteq A$ for some $i$ by the induction hypothesis. (Whether $C$ is a weave or not does not matter anymore.) Therefore, $B^\mathcal{R}_i = \ell(A^\mathcal{R}_i) \sqsubseteq \ell(A) = B$ once again. The proof is complete. \end{proof*}

\medskip

From the proof it is clear that in fact the $\mathcal{R}$-lineages we constructed account for all directionally rigid connected components of $L$ on the boards of Theorem \ref{thm:rigid}.

The question of characterising all directionally rigid complete leaper graphs on rectangular boards continues to make sense for skew but non-free leapers.

Let $N$ be a skew $(r, s)$-leaper with $r < s$. Let $d = \gcd(r, s)$ and let $N^\star$ be the unique relatively prime leaper proportional to $N$. Then the complete leaper graph of $N$ on the board of size $m \times n$ is directionally rigid if and only if the complete leaper graph of $N^\star$ on the board of size $\lceil m/d \rceil \times \lceil m/d \rceil$ is directionally rigid.

Thus it suffices to consider the special case when $d = 1$ and $N$ is either free or half-free. Since the former case is settled by Theorem \ref{thm:rigid}, suppose, from this point on, that $N$ is half-free.

\medskip

\begin{theorem*}{\ref{thm:rigid}$^\text{Half}$} Let $N$ be a skew half-free $(r, s)$-leaper of depth $\lambda$ with tails the $(r_i, s_i)$-leapers $N_i$, let $A^\mathcal{R}_0$ be the square board of side $r + s$, and let $A^\mathcal{R}_i$ be the board of size $(r + s - r_i + 1) \times (r + s + s_i)$ for all $i$ with $1 \le i \le \lambda$. Then the complete leaper graph of $N$ on $A$ is directionally rigid if and only if $A^\mathcal{R}_i \sqsubseteq A$ for some $i$. \end{theorem*}

\medskip

The proof is fully analogous to the proof of Theorem \ref{thm:rigid}. Observe also that the boards of Theorem $\ref{thm:rigid}^\text{Half}$ with $i \neq 1$ form a reduced basis of $\mathcal{R}_N$.

%% file: induction-on-descent-13-journey.tex
\section{Wazir Journeys} \label{journey}

We conclude our series of applications with a proof of Theorem \ref{thm:journey}.

We define a \emph{wazir-neighbourly} connected component of $L$ on $A$ to be one which contains two cells adjacent with respect to the wazir. (Or, equivalently, adjacent by side.) Thus there exists a wazir journey of $L$ on $A$ if and only if there exists a wazir-neighbourly connected component of $L$ on $A$.

Observe that a wazir journey of $L$ contains edges of both inclines. Furthermore, every wazir journey of $L$ is of odd length because its endpoints are of opposite parities.

First we study the wazir-neighbourly weaves of $L$.

Let $\kappa$ be the depth of $L$ and $\ecf(L) = [c_\kappa, \varepsilon_{\kappa - 1}, c_{\kappa - 1}, \varepsilon_{\kappa - 2}, \ldots, c_1]$.

\medskip

\begin{lemma} \label{ww} Suppose that there exists a wazir-neighbourly horizontal weave $U$ of $L$ on $A$. Then either (a) $\kappa = 2$; or (b) $\kappa = 3$ and $\varepsilon_1 = 1$. \end{lemma} 

\medskip

\begin{proof} By Lemmas \ref{pp} and \ref{wb}, all cells of $U$ in the same row of $A$ are of the same parity. Thus $U$ cannot contain two cells adjacent by side in the same row, and so it must contain two cells adjacent by side in the same column. Consequently, $W = \CompY(U)$ contains two vertices which differ by one.

Since a wazir journey of $L$ contains edges of both inclines, $W$ must be a compound weave. We let $q = 2cp + \varepsilon r$, we define $d$ to be the number of long edges in $W$, and we use the same notation for the vertices of $W$ as in Section \ref{weave}.

Then $w_{i', j'} - w_{i'', j''} = 1$ if and only if $\varepsilon(i' - i'')r + 1 = (j' - j'')p$. Since $|i' - i''| \le d < p/r$ by Lemma \ref{cwl}, it follows that either (a) $j' - j'' = 0$, when $|i' - i''| = 1$ and $r = 1$; or (b) $j' - j'' = 1$, when $|i' - i''| = d$ and $dr + 1 = p$.

Case (a) implies immediately that $\kappa = 2$.

In case (b), however, we must furthermore show that $d$ is even.

Since the length of a wazir journey of $L$ is always odd, the length of its $y$-projection must be odd as well. On the other hand, the length of a walk within $W$ connecting $w_{i', j'}$ and $w_{i'', j''}$ is necessarily of the same parity as $|i' - i''| + |j' - j''| = d + 1$. \end{proof}

\medskip

\begin{lemma} \label{wwai} Suppose that $\kappa = 2$ and there exists a wazir-neighbourly horizontal weave $U$ of $L$ on $A$. Let $m_\text{I} = q + 1$ and $n_\text{I} = p + 2q - 1$. Then $m \ge m_\text{I}$ and $n \ge n_\text{I}$. \end{lemma} 

\medskip

\begin{proof} That $m \ge q + 1 = m_\text{I}$ follows because every wazir journey of $L$ contains edges of both inclines.

Let $\alpha$ be a wazir journey of $L$ in $U$. By the proof of Lemma \ref{ww}, the endpoints of $\alpha$ are in the same column of $A$.

Let the $y$-projection of $\alpha$ be $\upsilon$ with signature $\mathtt{s}^{s_0}\mathtt{l}\mathtt{s}^{s_1}\mathtt{l} \ldots \mathtt{s}^{s_k}$ so that all of the $s_i$ are nonnegative integers. Since $U$ is a horizontal weave, $s_i$ is even for all $i$ with $0 < i < k$.

Observe that $s_0 + s_k$ is even as well. Indeed, we can assume without loss of generality that $\upsilon$ leads from $w_{i', j'}$ to $w_{i'', j''}$ in the setting of the proof of Lemma \ref{ww}. Thus $j' = j''$. However, $s_0$ is of the same parity as $j'$ and $s_k$ is of the same parity as $j''$.

Let $\chi$ be the $x$-projection of $\alpha$. Since the endpoints of $\alpha$ are in the same column of $A$, we get that $\chi$ is a closed walk in $\Pi_X$. Furthermore, since $\alpha$ is a wazir journey of $L$, both it and $\chi$ must be of odd length, and so $\chi$ cannot be balanced.

By our observations about the signature of $\upsilon$, there exists a cyclic shift of $\chi$ such that each maximal subwalk of long moves in it is of even length. We convert that cyclic shift of $\chi$ into an unbalanced closed walk $\chi^\star$ in $\Pi(p, 2q, I_X)$ just as in the proof of Lemma \ref{rw}. By Lemma \ref{pu}, it follows that $n = |I_X| \ge p + 2q - \gcd(p, 2q) + 1$. Lastly, $\gcd(p, 2q) = 2$ since $\kappa = 2$ implies that $p$ is even. \end{proof}

\medskip

\begin{lemma} \label{wwaii} In the setting of Lemma \ref{wwai}, there exists a wazir-neighbourly horizontal weave of $L$ on the board $A_\text{I}$ of size $m_\text{I} \times n_\text{I}$. \end{lemma} 

\medskip

For example, Figure \ref{ww-a} shows the unique wazir-neighbourly connected component of the $(2, 7)$-leaper on the board of size $8 \times 15$.

\medskip

\begin{proof} Let $A_\text{I} = E_X \times E_Y$.

\begin{figure}[t!] \centering \includegraphics{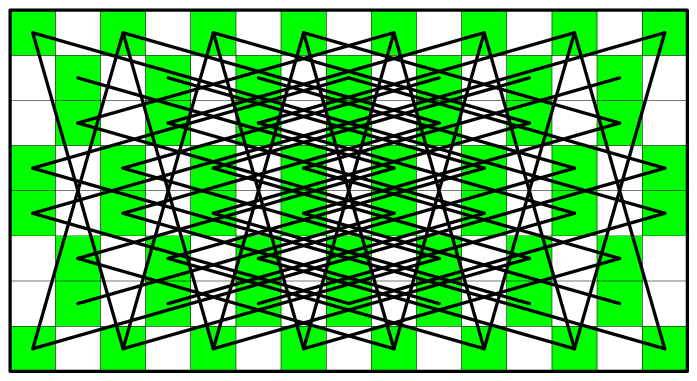} \caption{} \label{ww-a} \end{figure} 

Clearly, the connected component of $\Pi(p, q, E_Y)$ which contains its unique long edge is a compound weave $W$. Furthermore, there exists a walk $\upsilon$ in $W$ whose endpoints differ by one and whose signature is of the form $\mathtt{s}^{s_0}\mathtt{l}\mathtt{s}^{s_1}$ so that both of $s_0$ and $s_1$ are nonnegative integers and $s_0 + s_1$ is even.

Since $\kappa = 2$, we get that $p$ is even and $q$ is odd. By Lemma \ref{pc}, it follows that there exists a closed walk $\chi^\star$ in $\Pi(p, 2q, E_X)$ with $q$ short moves and $p/2$ long moves.

Just as in the proof of Lemma \ref{ro}, we convert $\chi^\star$ into an unbalanced closed walk $\chi$ in $\Pi(p, q, E_X)$ with $q$ short moves and $p$ long moves such that each maximal subwalk of long moves in it is of even length.

Thus there exists a cyclic shift of $\chi$ with signature of the form $\mathtt{l}^{t_0}\mathtt{s}\mathtt{l}^{t_1}\mathtt{s} \ldots \mathtt{l}^{t_q}$ such that all of the $t_i$ are nonnegative integers, $t_0$ is of the same parity as $s_0$, $t_q$ is of the same parity as $s_1$, and $t_i$ is even for all $i$ with $0 < i < q$.

Therefore, by Lemma \ref{psc}, there exists a wazir-neighbourly horizontal weave of $L$ on $A_\text{I}$. \end{proof}

\medskip

\begin{lemma} \label{wwbi} Suppose that $\kappa = 3$, $\varepsilon_1 = 1$, and there exists a wazir-neighbourly horizontal weave $U$ of $L$ on $A$. Let $m_\text{II} = p + q - c_1$ and $n_\text{II} = p + 2q - c_1 - 1$. Then $m = m_\text{II}$ and $n \ge n_\text{II}$. \end{lemma} 

\medskip

\begin{proof} We use the same notation as in the proof of Lemma \ref{ww}.

We begin with the height of $A$.

Suppose first that $\varepsilon_2 = -1$. Then $w_{d, 2c - 1} - q < y_\text{Min} \le w_{1, 0}$. Since $dr + 1 = p$ by the proof of Lemma \ref{ww}, it follows that $w_{d, 2c - 1} - w_{1, 0} = (d - 1)r + (2c - 1)p = q - 1$ and $y_\text{Min} = w_{1, 0}$. Similarly, because of $w_{0, 1} + q > y_\text{Max} \ge w_{d - 1, 2c}$ we conclude that $y_\text{Max} = w_{d - 1, 2c}$. Thus $m = y_\text{Max} - y_\text{Min} + 1 = (d - 2)r + 2cp + 1 = m_\text{II}$.

The case when $\varepsilon_2 = 1$ is analogous, except that we use the bounds $w_{0, 0} - p < y_\text{Min} \le w_{d, 0}$ and $w_{d, 2c} + p > y_\text{Max} \ge w_{0, 2c}$ instead.

We continue with the width of $A$.

Let $\alpha$ be a wazir journey of $L$ in $U$. By the proof of Lemma \ref{ww}, the endpoints of $\alpha$ are in the same column of $A$.

Let the $y$-projection of $\alpha$ be $\upsilon$ with signature $\mathtt{s}^{s_0}\mathtt{l}\mathtt{s}^{s_1}\mathtt{l} \ldots \mathtt{s}^{s_k}$ so that all of the $s_i$ are nonnegative integers. Since $U$ is a horizontal weave, $s_i$ is even for all $i$ with $0 < i < k$.

Observe also that $s_0 + s_k$ is odd. Indeed, we can assume without loss of generality that $\upsilon$ leads from $w_{i', j'}$ to $w_{i'', j''}$ in the setting of the proof of Lemma \ref{ww}. Thus $j' - j'' = 1$. However, $s_0$ is of the same parity as $j'$ and $s_k$ is of the same parity as $j''$.

Let $\chi$ be the $x$-projection of $\alpha$. As in the proof of Lemma \ref{wwai}, we get that $\chi$ is a closed walk in $\Pi_X$. This time around, though, one maximal subwalk of long moves in $\chi$ is odd and all of the others are even.

By the same method as in the proofs of Lemmas \ref{rw} and \ref{wwai}, we convert $\chi$ into a walk $\chi^\star$ in $\Pi(p, 2q, I_X)$ whose endpoints differ by $q$. We are left to demonstrate that the existence of such a walk implies $n = |I_X| \ge n_\text{II}$.

Since $\kappa = 3$, we get that $p$ is odd, $q$ is even, and $\gcd(p, 2q) = 1$. By Lemma \ref{pa}, it follows that $\Pi(p, 2q, I_X)$ is connected when $n \ge p + 2q - 1$. Thus let $n = p + 2q - n_\diamond$ so that, without loss of generality, $n_\diamond$ is a positive integer.

Then we can orient $\chi^\star$ so that all of its short moves point to the right and all of its long moves point to the left. With that orientation, let $\chi^\star$ lead from $u'$ to $u''$ with $|u' - u''| = q$ and let also $u'_0 \to u''_0$, $u'_1 \to u''_1$, \ldots, $u'_k \to u''_k$ be all long moves of $\chi^\star$ in the order in which they occur.

We consider the case when $\varepsilon_2 = -1$ in detail. The opposite case, when $\varepsilon_2 = 1$, is analogous.

Since $\kappa = 3$ and $\varepsilon_2 = -1$, we get that $q \equiv -c_1 \pmod p$. By translating $I_X$ as needed, we can ensure without loss of generality that $u' \equiv 0 \pmod p$ and $u'' \equiv \pm c_1 \pmod p$. Thus $u'_i \equiv 2c_1i \pmod p$ for all $i$ and $u'_k = u''_k + 2q \equiv u'' + 2q \equiv (\pm 1 - 2)c_1 \pmod p$.

Consequently, each element of the set $S = \{0, 2c_1, 4c_1, \ldots, (p - 3)c_1\}$ is congruent to some $u'_i$ modulo $p$.

On the other hand, since each $u'_i$ is the right vertex of a long edge in $\Pi(p, 2q, I_X)$, it follows that no $u'_i$ is congruent modulo $p$ to an element of $[x_\text{Max} + 1; x_\text{Max} + n_\diamond]$. Therefore, $n_\diamond$ cannot exceed the greatest size of an integer interval which avoids all elements of $S$ modulo $p$.

Since $\kappa = 3$ and $\varepsilon_1 = 1$, we get that $p = c_1c_2 + 1$. Thus modulo $p$ the set $S$ becomes $\{0, 2c_1, 4c_1, \ldots, -1, 2c_1 - 1, 4c_1 - 1, \ldots, -2, 2c_1 - 2, 4c_1 - 2, \ldots, -c_1 + 1, c_1 + 1, 3c_1 + 1, \ldots, (c_2 - 3)c_1 + 1\}$. Consequently, the $\subseteq$-maximal integer intervals which avoid all elements of $S$ modulo $p$ are $[1; c_1]$, $[2c_1 + 1; 3c_1]$, \ldots, $[(c_2 - 4)c_1 + 1; (c_2 - 3)c_1]$, $[(c_2 - 2)c_1 + 1; (c_2 - 1)c_1 + 1]$ as well as their translations by multiples of $p$.

The last one of these is of size $c_1 + 1$ and the others are all of size $c_1$. Therefore, $n_\diamond \le c_1 + 1$ and $n \ge n_\text{II}$. \end{proof}

\medskip

\begin{lemma} \label{wwbii} In the setting of Lemma \ref{wwbi}, there exists a wazir-neighbourly horizontal weave of $L$ on the board $A_\text{II}$ of size $m_\text{II} \times n_\text{II}$. \end{lemma} 

\medskip

For example, Figure \ref{ww-b} shows one of the two wazir-neighbourly connected components of the $(5, 8)$-leaper on the board of size $11 \times 18$. The second one is its mirror image.

\medskip

\begin{proof} Let $A_\text{II} = E_X \times E_Y$.

\begin{figure}[t!] \centering \includegraphics{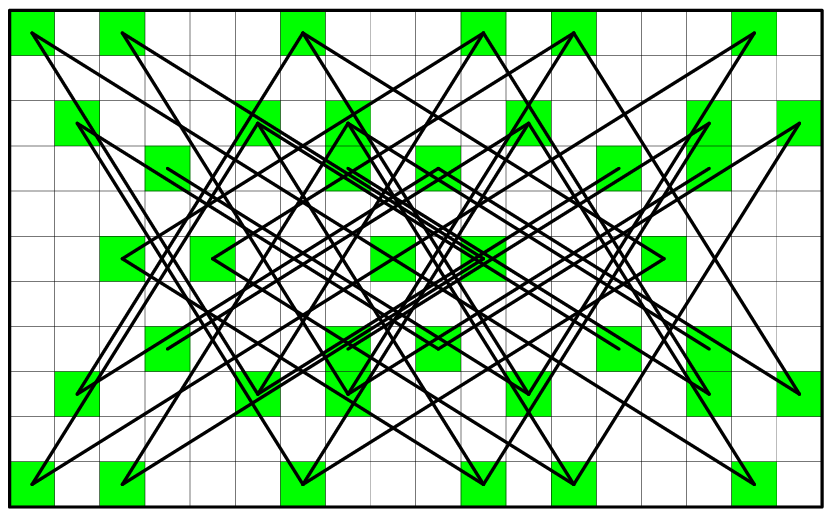} \caption{} \label{ww-b} \end{figure} 

It is straightforward to see that the connected component of $\Pi(p, q, E_Y)$ which contains $\min E_Y$ and $\max E_Y$ is a compound weave $W$ with $c_2$ long edges.

Once again, we use the same notation for the vertices of $W$ as in Section \ref{weave} and the proof of Lemma \ref{ww}. Then $W$ contains a walk $\upsilon$ whose endpoints $w_{i', j'}$ and $w_{i'', j''}$ differ by one with $|i' - i''| = c_2$ and $j' - j'' = 1$.

As in the proof of Lemma \ref{wwbi}, it follows that the signature of $\upsilon$ is of the form $\mathtt{s}^{s_0}\mathtt{l}\mathtt{s}^{s_1}\mathtt{l} \ldots \mathtt{s}^{s_k}$ so that $s_i$ is an even nonnegative integer for all $i$ with $0 < i < k$ and $s_0 + s_k$ is odd. Observe also that $k$ must be even because, being the number of long moves in $\upsilon$, it is of the same parity as $|i' - i''|$.

By inverting the reasoning in the proof of Lemma \ref{wwbi}, we see that there exists a walk $\chi^\star$ in $\Pi(p, 2q, E_X)$ whose endpoints differ by $q$. With the help of the same method as in the proofs of Lemmas \ref{ro} and \ref{wwaii}, we convert $\chi^\star$ into a closed walk $\chi$ in $\Pi(p, q, E_X)$ such that one maximal subwalk of long moves in $\chi$ is odd and all of the others are even.

Thus there exists a cyclic shift of $\chi$ with signature of the form $\mathtt{l}^{t_0}\mathtt{s}\mathtt{l}^{t_1}\mathtt{s} \ldots \mathtt{l}^{t_\tau}$ such that all of the $t_i$ are nonnegative integers, $t_0$ is of the same parity as $s_0$, $t_\tau$ is of the same parity as $s_k$, and $t_i$ is even for all $i$ with $0 < i < \tau$. Observe also that $\tau$ must be even because $p$ is odd and $q$ is even by virtue of $\kappa = 3$ whereas $\chi$ is a closed walk containing a total of $\tau$ short moves.

Therefore, by Lemma \ref{psc}, there exists a wazir-neighbourly horizontal weave of $L$ on $A_\text{II}$. \end{proof}

\medskip

We proceed to study the effect of lifting and lowering on the wazir journeys of skew free leapers.

We say that two edges $a'b'$ and $a''b''$ of $L$ form a \emph{wazir-suitable configuration} when $a' + (0, 1) = a''$, cells $b'$ and $b''$ are in the same column, and both edges $a'b'$ and $a''b''$ are slight. We also say that a walk of $L$ is \emph{wazir-suitable} when its first edge and its last edge form a wazir-suitable configuration.

Let $\ell$ be a lifting transformation with $\ell(L) = M$.

\medskip

\begin{lemma} \label{wjf} Suppose that $\ell = f$ and the leaper graph $G$ of $L$ on $A$ contains a wazir-suitable angular walk of $L$. Then $f(G)$ contains a wazir-suitable angular walk of $M$. \end{lemma} 

\medskip

\begin{proof} When $\ell = f$, if two edges of $L$ form a wazir-suitable configuration then so do their cross-edges. Since we can assume without loss of generality that the vertices and edges of $G$ coincide with the vertices and edges of that wazir-suitable angular walk of $L$, the result follows by Lemma \ref{ace}. \end{proof}

\medskip

\begin{lemma} \label{wjgh} Suppose that $\ell \in \{g, h\}$ and the leaper graph $G$ of $L$ on $A$ contains a wazir-suitable walk of $L$. Then $\ell(G)$ contains a wazir-suitable angular walk of~$M$. \end{lemma} 

\medskip

\begin{proof} Consider two edges $a'b'$ and $a''b''$ of $L$ on $A$ in a wazir-suitable configuration with $a'$, $b'$, $a''$, and $b''$ as in the definition. When $\ell \in \{g, h\}$, cells $a'$ and $a''$ are matched with two edges of $M$ on $B$ which form a wazir-suitable configuration as well. Since we can assume without loss of generality that the vertices and edges of $G$ coincide with the vertices and edges of that wazir-suitable walk of $L$, the result follows by Lemmas \ref{lg} and \ref{lh}. \end{proof}

\medskip

\begin{lemma} \label{wjlow} Suppose that the angular connected component $H$ of $M$ on $\ell(A)$ contains a wazir journey of $M$. Then $\ellinv(H)$ contains a wazir journey of $L$. \end{lemma} 

\medskip

\begin{proof} Let $B = \ell(A)$ and $G = \ellinv(H)$.

Since $G$ is a connected leaper graph of $L$ on $A$ by Lemma \ref{alc}, it suffices to show that $G$ contains two cells adjacent by side.

By transposing $A$ and $B$ if needed, we can assume without loss of generality that $H$ contains two cells $a'$ and $a''$ with $a' + (0, 1) = a''$. Let $a'b'$ and $a''b''$ be two edges of $M$ on $B$ incident with them. We distinguish two cases, as follows.

\smallskip

\emph{Case 1}. We can choose $b'$ and $b''$ so that moves $a' \to b'$ and $a'' \to b''$ of $M$ point in the same direction.

Let $c'$ and $d'$ be the two cells matched with edge $a'b'$ and let $c'' = c' + (0, 1)$ and $d'' = d' + (0, 1)$ be the two cells matched with edge $a''b''$ so that $c'$ lies below $d'$ and $c''$ lies below $d''$.

Since each edge of $H$ is part of an acute angle in $H$, by Lemma \ref{core} it follows that at least one of $c'$ and $d'$ is in $A$ and at least one of $c''$ and $d''$ is in $A$.

Suppose, for the sake of contradiction, that both of $c'$ and $d''$ are outside of $A$. Then $c''$ must be in the lowermost row of $A$ and $d'$ must be in the topmost row of $A$. Consequently, $m \le q$ and we arrive at a contradiction with $c''$ and $d'$ being in the same connected leaper graph $G$ of $L$ on $A$.

Therefore, up to symmetry, $c'$ is in $A$. Since at least one of $c''$ and $d''$ is in $A$ as well, and $c''$ is in the bounding box of $c'$ and $d''$, we conclude that in fact both cells $c'$ and $c''$ are in $A$ and also in $G$.

\smallskip

\emph{Case 2}. We cannot choose $b'$ and $b''$ as in Case 1.

Thus both cells $b' + (0, 1)$ and $b'' - (0, 1)$ are outside of $B$. Consequently, $b'$ must lie above $a'$ in the topmost row of $B$ and $b''$ must lie below $a''$ in the lowermost row of $B$.

Suppose first that both edges $a'b'$ and $a''b''$ are steep. Then $\widehat{m} = 2\widehat{q}$. On the other hand, since cells $a'$ and $a''$ are in the same connected component of $M$ on $B$, by Lemmas \ref{pa} and \ref{pp} it follows also that $\widehat{n} \ge \widehat{p} + \widehat{q}$. But then $H$ coincides with the complete leaper graph of $M$ on $B$ by Theorem \ref{thm:conn}, and the result is immediate.

Suppose, otherwise, that at least one of edges $a'b'$ and $a''b''$ is slight. Up to symmetry, $b' = a' + (\widehat{q}, \widehat{p})$. Thus cell $a'' + (\widehat{q}, -\widehat{p})$ is in the bounding box of cells $b'$ and $b''$, and so we can assume without loss of generality that in fact $b'' = a'' + (\widehat{q}, -\widehat{p})$.

It follows that $\widehat{m} = 2\widehat{p}$. Since also $\widehat{m} > 2\Mar$, we conclude that $\ell = g$. Therefore, both cells $a' + (\Lat, 0)$ and $a'' + (\Lat, 0)$ are in $G$. \end{proof}

\medskip

We are only left to put the pieces together. The overall structure of the proof will be the same as with Theorem \ref{thm:rigid}.

\medskip

\begin{proof*}{Proof of sufficiency for Theorem \ref{thm:journey}} We define the \emph{$\mathcal{W}$-lineage} of the wazir to be the perfect lineage originated by the wazir with the board of size $1 \times 1$. Observe that the knight graph in this lineage contains a wazir-suitable angular walk.

We also define a \emph{$\mathcal{W}$-lineage} of the skew free leaper $L$ as follows: (a) When $\kappa = 2$, a $\mathcal{W}$-lineage of $L$ is a perfect lineage originated by $L$ with a board and a wazir-neighbourly connected component of $L$ on that board as in Lemma \ref{wwaii}; and (b) When $\kappa = 3$ and $\varepsilon_1 = 1$, a $\mathcal{W}$-lineage of $L$ is a perfect lineage originated by $L$ with a board and a wazir-neighbourly connected component of $L$ on that board as in Lemma \ref{wwbii}.

The proofs of Lemmas \ref{wwaii} and \ref{wwbii} show that each one of these connected components contains two cells adjacent by side and in the same column. Since $n_\text{I} \ge 2q$ and $n_\text{II} \ge 2q$, it follows that each one of them contains also a wazir-suitable walk.

By Lemma \ref{offset}, the $\mathcal{W}$-lineages which contain $L$ assign to it precisely the boards $A^\mathcal{W}_i$.

On the other hand, by induction on descent with the help of Lemmas \ref{wjf} and \ref{wjgh}, it follows that every leaper graph of a skew free leaper in a $\mathcal{W}$-lineage does indeed contain a wazir journey of that leaper. \end{proof*}

\medskip

It is not too difficult to see that the number of wazir-neighbourly connected components of $L$ on the corresponding board is one in case (a) and two in case (b). This is also the number of $\mathcal{W}$-lineages that a skew free leaper $L$ satisfying the conditions of one of cases (a) and (b) originates. Furthermore, in case (b) reflection with respect to the horizontal axis of symmetry of the board preserves each one of the two connected components whereas reflection with respect to the vertical axis of symmetry of the board swaps them.

\medskip

\begin{proof*}{Proof of necessity for Theorem \ref{thm:journey}} We establish one slightly stronger result: Define $\overline{0} = 0$, $\overline{1} = 2$, and $\overline{2} = 3$. Suppose that there exists a wazir-neighbourly connected component of $L$ on $A$. When that connected component is a weave, $\kappa \in \{2, 3\}$ and $A^\mathcal{W}_i \sqsubseteq A$ for the $i$ with $\overline{i} = \kappa$. Otherwise, when that connected component is a non-weave, $A^\mathcal{W}_i \sqsubseteq A$ for some $i$ with $\overline{i} < \kappa$.

We proceed by induction on descent.

Our base case is the knight, and with it the verification is straightforward.

For the induction step, let $\ell$ be a lifting transformation with $\ell(L) = M$ and suppose that our strengthening of the necessity part of Theorem \ref{thm:journey} holds for $L$.

Define the boards $B^\mathcal{W}_i$ relative to $M$ in the same way as we defined the boards $A^\mathcal{W}_i$ relative to $L$. When $\ell = f$, we get that $f(A^\mathcal{W}_i) = B^\mathcal{W}_i$ if and only if $\overline{i} < \kappa$. Otherwise, when $\ell \in \{g, h\}$, we get that $\ell(A^\mathcal{W}_i) = B^\mathcal{W}_i$ for all $i$ such that the left-hand side is well-defined.

Let $D$ be a wazir-neighbourly connected component of $M$ on a board $B$.

When $D$ is a weave, we are done immediately by Lemmas \ref{ww}, \ref{wwai}, and \ref{wwbi}.

Suppose, throughout the rest of the proof, that $D$ is a non-weave.

Observe that $\ellinv(B)$ is well-defined because a wazir journey necessarily contains edges of both inclines. Let $\ellinv(B) = A$.

Since $D$ is a non-weave, it is an angular connected component of $M$ on $B$ by Lemma \ref{naw}. Thus by Lemma \ref{wjlow} there exists a wazir journey of $L$ in $\ellinv(D)$. Let $C$ be the connected component of $L$ on $A$ which contains that wazir journey of $L$.

We consider the case when $\ell = f$ first. Then $C$ is a non-weave by Lemma \ref{wf}, and so $A^\mathcal{W}_i \sqsubseteq A$ for some $i$ with $\overline{i} < \kappa$ by the induction hypothesis. Therefore, $B^\mathcal{W}_i = f(A^\mathcal{W}_i) \sqsubseteq f(A) = B$, as needed.

We continue with the case when $\ell \in \{g, h\}$. Then $A^\mathcal{W}_i \sqsubseteq A$ for some $i$ by the induction hypothesis. (Whether $C$ is a weave or not does not matter anymore.) Therefore, $B^\mathcal{W}_i = \ell(A^\mathcal{W}_i) \sqsubseteq \ell(A) = B$ once again. The proof is complete. \end{proof*}

\medskip

Just as in Section \ref{rigid}, from the proof it is clear that in fact the $\mathcal{W}$-lineages we constructed account for all wazir-neighbourly connected components of $L$ on the boards of Theorem \ref{thm:journey}.

%% file: induction-on-descent-14-further.tex
\section{Further Work} \label{further}

Let $M$ be an $(r, s)$-leaper and let $N$ be a $(u, v)$-leaper. An \emph{$M$-journey} of $N$ is a walk of $N$ whose endpoints are adjacent with respect to $M$. This is also known as an \emph{$(r, s)$-journey} of $N$. \cite{J}

For example, when $M$ and $N$ form a parent-and-child pair as in Section \ref{desc}, there exist both a three-move $M$-journey of $N$ and a three-move $N$-journey of $M$.

It is straightforward to see that an $M$-journey of $N$ exists if and only if $M$ is at most as free as $N$, in the sense that the number of connected components in the complete leaper graph of $M$ on the infinite board $\mathbb{Z} \times \mathbb{Z}$ is greater than or equal to the analogous number for $N$. We give expressions for these numbers in terms of $r$, $s$, $u$, and $v$ in Section \ref{prelim}.

We write $\mathcal{J}^M_N(A)$ for ``there exists an $M$-journey of $N$ on $A$''. Clearly, $\mathcal{J}^M_N$ is both symmetric and monotone for all $M$ and $N$.

For example, if $M$ is the wazir, then $\mathcal{J}^\text{Wazir}_N$ is the same property as $\mathcal{W}_N$.

Let $d = \gcd(u, v)$. For an $M$-journey of $N$ to exist, $d$ must divide $r$ and $s$ as well. When it does, let $M^\star$ be the $(r/d, s/d)$-leaper and let $N^\star$ be the $(u/d, v/d)$-leaper. Then $\mathcal{J}^M_N(m, n)$ if and only if $\mathcal{J}^{M^\star}_{N^\star}(\lceil m/d \rceil, \lceil n/d \rceil)$. Thus it suffices to consider the case when $d = 1$, $N$ is either free or half-free, and the considerations of Section \ref{desc} apply to $N$.

Given one concrete leaper $M$, what can we say about the bases of the properties $\mathcal{J}^M_N$ as $N$ varies over all relatively prime leapers?

For example, when $M$ is the wazir, by Theorem \ref{thm:journey} we get that: (a) The number of boards in such a basis is bounded from above by an absolute constant; and (b) The sizes of the boards in such a basis are of the form $(u + v + w_i) \times (u + v + z_i)$ so that both the total number of boards and all of the $w_i$ and $z_i$ depend only on a tail of $N$ whose depth is bounded from above by an absolute constant as well.

What other leapers $M$ behave in each one of these two ways? How can we find or estimate the corresponding absolute constants?

We generalise the notion of a journey as follows.

We define a \emph{pattern} $S$ to be a finite set of cells. We say that $N$ \emph{realises} $S$ on $A$ when there exists a translation copy of $S$ all of whose cells are in the same connected component of $N$ on $A$.

Let us write $\mathcal{S}^S_N(A)$ for ``$N$ realises $S$ on $A$''. Clearly, $\mathcal{S}^S_N$ is monotone for all patterns $S$ and all leapers $N$ but it is not in general symmetric.

Once again, it suffices to consider the case when $N$ is either free or half-free, by the same reasoning as above.

We can model journeys in terms of two-cell patterns. Observe that $\mathcal{J}^M_N(A)$ if and only if either $\mathcal{S}^{\{(0, 0), (r, s)\}}_N(A)$ or $\mathcal{S}^{\{(0, 0), (s, r)\}}_N(A)$. Furthermore, the basis of each one of $\mathcal{S}^{\{(0, 0), (r, s)\}}_N$ and $\mathcal{S}^{\{(0, 0), (s, r)\}}_N$ consists of the transposes of the boards in the basis of the other. Thus if we know the basis of $\mathcal{S}^{\{(0, 0), (r, s)\}}_N$ then we can also determine the bases of $\mathcal{S}^{\{(0, 0), (s, r)\}}_N$ and $\mathcal{J}^M_N$.

Given one concrete pattern $S$, what can we say about the bases of the properties $\mathcal{S}^S_N$ as $N$ varies over all relatively prime leapers?

In the special case when $S = \{(0, 0), (0, 1)\}$, these bases are as follows. Clearly, non-free leapers cannot realise this pattern. On the other hand, the proof of Theorem \ref{thm:journey} shows that $L$ does realise it on all boards of the form $A^\mathcal{W}_i$. Therefore, these boards are also the elements of the basis of $\mathcal{S}^{\{(0, 0), (0, 1)\}}_L$.

The same questions that we raised because of Theorem \ref{thm:journey} regarding journeys continue to make sense here as well: What are the patterns $S$ which exhibit either one or both of the behaviours (a) and (b) described above? Furthermore, how can we find or estimate the corresponding absolute constants?

This concludes our discussion of induction on descent in the study of leaper graphs.

%% file: induction-on-descent-15-biconn.tex
\section{Biconnectedness} \label{biconn}

A graph $G$ is \emph{biconnected} when it contains at least three vertices and, for each vertex $v$ of $G$, the subgraph $G \setminus v$ of $G$ is connected. (We obtain $G \setminus v$ from $G$ by deleting $v$ together with all edges incident with it.) A biconnected graph is necessarily connected as well.

Knuth writes in \cite{K} that it would be ``of interest [\ldots]\ to characterize the smallest boards on which [complete]\ leaper graphs are biconnected''. We proceed to resolve this question completely.

We write $\mathcal{B}_L(A)$ for ``the complete leaper graph of $L$ on $A$ is biconnected''.

Clearly, $\mathcal{B}_L$ is symmetric. Let us verify that it is monotone as well. The proof proceeds along the same lines as our proofs of monotonicity for $\mathcal{C}_L$ and $\mathcal{E}_L$.

\medskip

\begin{proof*}{Proof of monotonicity for $\mathcal{B}_L$} Suppose that $\mathcal{B}_L(m, n)$. It suffices to show that then $\mathcal{B}_L(m, n + 1)$ as well.

Let $A^\star$ be a board of size $m \times (n + 1)$ and let $A'$ and $A''$ be its left and right subboards of size $m \times n$. Then $\mathcal{B}_L(m, n)$ implies $m \ge 2$, $n \ge 2$, and $|A' \cap A''| \ge 2$.

Let $G'$, $G''$, and $G^\star$ be the complete leaper graphs of $L$ on $A'$, $A''$, and $A^\star$, respectively. Let also $a$ be an arbitrary cell of $A^\star$. Then both of $G' \setminus a$ and $G'' \setminus a$ are connected since $\mathcal{B}_L(A')$ and $\mathcal{B}_L(A'')$. (When $v$ is not a vertex of $G$, we define $G \setminus v$ to be $G$ itself.) Furthermore, $G' \setminus a$ and $G'' \setminus a$ have at least one vertex in common because $A' \cap A''$ contains at least two cells. Therefore, $G^\star \setminus a$ is connected as well. \end{proof*}

\medskip

We characterise all biconnected complete leaper graphs on rectangular boards as follows.

\medskip

\begin{theorem} \label{thm:biconn} Let $L$ be a skew free $(p, q)$-leaper with $p < q$ and let $A^\mathcal{B}$ be the board of size $(p + q) \times 2q$. Then the complete leaper graph of $L$ on $A$ is biconnected if and only if $A^\mathcal{B} \sqsubseteq A$. \end{theorem}

\medskip

Or, equivalently, $\{A^\mathcal{B}\}$ is a reduced basis of $\mathcal{B}_L$.

Observe that $A^\mathcal{B}$ coincides with the board $A^\mathcal{C}$ of Theorem \ref{thm:conn}. Thus one more way to state Theorem \ref{thm:biconn} is this: The complete leaper graph of $L$ on $A$ is biconnected if and only if it is nontrivially connected.

\medskip

\begin{proof*}{Proof of Theorem \ref{thm:biconn}} Let $G^\mathcal{B}$ be the complete leaper graph of $L$ on $A^\mathcal{B}$.

The necessity part of Theorem \ref{thm:conn} implies the necessity part of Theorem \ref{thm:biconn} as well.

We are left to take care of sufficiency. Since $\mathcal{B}_L$ is monotone, this amounts to showing that $G^\mathcal{B}$ is biconnected. Observe that $G^\mathcal{B}$ is already connected by the sufficiency part of Theorem \ref{thm:conn}. We must demonstrate that it remains connected when we delete an arbitrary vertex from it.

To that end, first we verify that $G^\mathcal{B}$ remains connected when we delete an arbitrary edge from it.

Let $e$ be an edge of $G^\mathcal{B}$. Let also $A^\star$ be a subboard of $A^\mathcal{B}$ of size $(p + q) \times (p + q)$ that contains $e$ and let $G^\star$ be the complete leaper graph of $L$ on $A^\star$. Since each vertex of $G^\star$ is of degree either zero or two, the connected component of $G^\star$ which contains $e$ is a cycle. Thus $e$ is part of a cycle in $G^\mathcal{B}$, and so indeed its deletion cannot disconnect $G^\mathcal{B}$.

Consider now an arbitrary cell $a$ of $A^\mathcal{B}$. It is straightforward to see that the degree of $a$ in $G^\mathcal{B}$ does not exceed three. On the other hand, since the deletion of an edge cannot disconnect $G^\mathcal{B}$, at least two edges of $G^\mathcal{B}$ must join $a$ to each connected component of $G^\mathcal{B} \setminus a$. Therefore, $G^\mathcal{B} \setminus a$ cannot consist of two or more connected components. \end{proof*}

\medskip

Let $k$ be a positive integer. A graph $G$ is \emph{$k$-connected} when it contains more than $k$ vertices and, however we delete fewer than $k$ vertices from it together with all edges incident with them, the result is a connected subgraph of $G$.

For complete leaper graphs on rectangular boards, the cases of $k = 1$ and $k = 2$ are characterised by Theorems \ref{thm:conn} and \ref{thm:biconn}, respectively. On the other hand, a leaper graph of this kind cannot be $k$-connected with $k \ge 3$ because in it each corner cell of the board is of degree at most two.

%% file: induction-on-descent-16-weave-ii.tex
\section{Weaves II} \label{weave-ii}

We go on to describe the projection graphs with parameters $p$ and $q$ which contain different kinds of weaves.

When $p = 1$, there exists a weave in $\Pi(1, q, s)$ if and only if $s \le q$. This weave is always simple and it is a non-singleton if and only if $2 \le s \le q$.

Throughout the rest of this section, we consider the case when $p \ge 2$.

\medskip

\begin{proposition} \label{sws} Suppose that $p \ge 2$. Then $\Pi(p, q, s)$ contains a simple weave if and only if $s < p + q - (q \mmod p)$. \end{proposition} 

\medskip

\begin{proof} Let $I$ be an integer interval of size $s$.

(Necessity) Let $w_0$, $w_1$, \ldots, $w_d$ be the vertices of a simple weave in $\Pi(p, q, I)$ with $w_i + p = w_{i + 1}$ for all $i$. Then $w_0 - p < \min I$ and $w_d + p > \max I$ imply $s + 1 \le (d + 2)p = (p + q) - q + (d + 1)p$ whereas $w_d - q < \min I$ and $w_0 + q > \max I$ similarly imply $s + 1 \le 2q - dp = (p + q) + q - (d + 1)p$. Therefore, $p + q - s - 1 \ge |q - (d + 1)p| \ge q \mmod p$.

(Sufficiency) Choose $d_\text{I}$ so that $q \mmod p = |q - (d_\text{I} + 1)p|$ and let $J$ be an integer interval of size $p + q - (q \mmod p) - 1$. It is straightforward to verify that $\Pi(p, q, J)$ contains a simple weave $\Omega$ with $d_\text{I}$ short edges such that $J$ and $\Omega$ share the same center of symmetry.

Suppose, without loss of generality, that $I$ contains at least one vertex of $\Omega$ and $I \subseteq J$. Then the vertices of $\Omega$ in $I$ form a simple weave in $\Pi(p, q, I)$. \end{proof}

\medskip

\begin{proposition} \label{ws} Suppose that $p \ge 2$. Then $\Pi(p, q, s)$ contains a weave if and only if $s < p + q - (p \mmod [q \mmod 2p])$. \end{proposition} 

\medskip

The proof proceeds mostly along the same lines as with Proposition \ref{sws}.

\medskip

\begin{proof} Let $I$ be an integer interval of size $s$.

(Necessity) When $\Pi(p, q, I)$ contains a simple weave, we are done by the necessity part of Proposition \ref{sws}.

Suppose, otherwise, that $\Pi(p, q, I)$ contains a compound weave $W$. We let $q = 2cp + \varepsilon r$, we define $d$ to be the number of long edges in $W$, and we use the same notation for the vertices of $W$ as in Section \ref{weave}.

We consider the case when $\varepsilon = -1$ in detail. The opposite case, when $\varepsilon = 1$, is analogous. In it, we work with $w_{0, 0}$ and $w_{d, 2c}$ instead of $w_{0, 1}$ and $w_{d, 2c - 1}$.

Thus suppose that $\varepsilon = -1$. Then $w_{0, 1} - p < \min I$ and $w_{d, 2c - 1} + p > \max I$ imply $s + 1 \le 2p - w_{0, 1} + w_{d, 2c - 1} = 2p + dr + (2c - 2)p = (p + q) - p + (d + 1)r$ whereas $w_{d, 2c - 1} - q < \min I$ and $w_{0, 1} + q > \max I$ similarly imply $s + 1 \le 2q + w_{0, 1} - w_{d, 2c - 1} = 2q - dr - (2c - 2)p = (p + q) + p - (d + 1)r$. Therefore, $p + q - s - 1 \ge |p - (d + 1)r| \ge p \mmod r$.

(Sufficiency) Choose $d_\text{II}$ so that $p \mmod r = |p - (d_\text{II} + 1)r|$, where for convenience $r$ denotes $q \mmod 2p$ as above.

When $d_\text{II} = 0$, it follows also that $p \mmod r = q \mmod p$. This case is covered by the sufficiency part of Proposition \ref{sws}.

Suppose, throughout the rest of the proof, that $d_\text{II} \ge 1$.

Let $J$ be an integer interval of size $p + q - (p \mmod r) - 1$. It is straightforward to verify that $\Pi(p, q, J)$ contains a compound weave $\Omega$ with $d_\text{II}$ long edges such that $J$ and $\Omega$ share the same center of symmetry.

Suppose, without loss of generality, that $I$ contains at least one vertex of $\Omega$ and $I \subseteq J$. Then the vertices of $\Omega$ in $I$ form some number of connected components of $\Pi(p, q, I)$ each one of which, being a subgraph of $\Omega$, is a weave. \end{proof}

\medskip

\begin{proposition} \label{nws} Suppose that $p \ge 2$. (a) When $4p < 3q$, let $s_\vee = p + q - (p \mmod [q \mmod 2p])$. (b) Otherwise, when $4p \ge 3q$, let $s_\vee = 3q - 2p$. Then $\Pi(p, q, s)$ contains a non-singleton weave if and only if $p < s < s_\vee$. \end{proposition} 

\medskip

\begin{proof} (Necessity) When $4p < 3q$, we are done by Proposition \ref{ws}.

Otherwise, when $4p \ge 3q$ and there exists a non-singleton weave $W$ in $\Pi(p, q, s)$, it follows by Proposition \ref{ws} that $s < 2p$. Consequently, the signature of $W$ must be either $\mathtt{s}$ or $\mathtt{s}\mathtt{l}\mathtt{s}$. The former case implies $s < 2q - p$ and the latter one implies $s < 3q - 2p$.

(Sufficiency) Let $J$ be an integer interval of size $s_\vee - 1$.

When $4p < 3q$, either $J$ is of the same size as in the proof of Proposition \ref{sws} and $d_\text{I}$ is positive or $J$ is of the same size as in the proof of Proposition \ref{ws} and $d_\text{II}$ is positive. Both cases yield a non-singleton weave $\Omega$ in $\Pi(p, q, J)$.

Otherwise, when $4p \ge 3q$, it is straightforward to verify that $J$ contains a weave $\Omega$ with signature $\mathtt{s}\mathtt{l}\mathtt{s}$ such that $J$ and $\Omega$ share the same center of symmetry.

Consider an integer interval $I$ of size $s$ such that a short edge $e$ of $\Omega$ is in $I$ and $I \subseteq J$. Then the connected component of $\Pi(p, q, I)$ which contains $e$, being a subgraph of $\Omega$, is also a non-singleton weave. \end{proof}

%% file: induction-on-descent-17-ang-ii.tex
\section{Angularity II} \label{ang-ii}

Here we collect a number of additional results regarding angularity.

We define a \emph{zigzag} of $L$ on $A$ to be a path of $L$ on $A$ all of whose angles are laterally acute and which cannot be extended so that this property is preserved.

Thus all edges of a zigzag are of the same incline. A zigzag of slight edges alternates between two columns of $A$ spaced $q$ units apart, and similarly a zigzag of steep edges alternates between two rows of $A$ spaced $q$ units apart. Every edge of $L$ on $A$ belongs to exactly one zigzag. Furthermore, all edges of a zigzag are in the same angular component of $L$ on $A$.

\medskip

\begin{proposition} \label{ai} Let $C$ be an angular component of $L$ on $A$. Then $C$ is an induced leaper graph. \end{proposition} 

\medskip

\begin{proof} Let $a'a''$ be any edge of $L$ on $A$ such that both of $a'$ and $a''$ are cells of $C$. We must show that $a'a''$ is an edge of $C$.

Let $a'b'$ and $a''b''$ be edges of $C$. When angle $b'a'a''$ is zero, acute, right, or diagonally obtuse, we are done by Lemma \ref{rdo}. Same goes for angle $a'a''b''$. We are left to consider the case when each one of angles $b'a'a''$ and $a'a''b''$ is either laterally obtuse or straight.

Up to symmetry, all three edges $a'a''$, $a'b'$, and $a''b''$ are slight. Let $Z$, $Z'$, and $Z''$ be their zigzags, respectively. Observe that $Z$ and $Z'$ are reflections of one another with respect to the vertical line through $a'$, and similarly $Z$ and $Z''$ are reflections of one another with respect to the vertical line through $a''$.

Suppose, for the sake of contradiction, that no cell of $Z$ is incident with a steep edge of $L$ on $A$. Then the same must be true of its reflections $Z'$ and $Z''$. Consequently, $Z'$ is the entire angular component of edge $a'b'$ and $Z''$ is the entire angular component of edge $a''b''$. We have arrived at a contradiction with both edges $a'b'$ and $a''b''$ being in $C$.

Thus some cell $c$ of $Z$ must be incident with a steep edge $cd$ of $L$ on $A$. Up to symmetry, $c$ is in the column of $a'$.

Then $c$ is in $Z'$ as well. By Lemma \ref{rdo}, edge $cd$ is in the same angular component as $Z$, and also in the same angular component as $Z'$. Therefore, all edges of $Z$ and $Z'$ are in the same angular component and $a'a''$ is an edge of $C$. \end{proof}

\medskip

Observe that every non-angular connected component of $L$ on $A$ decomposes uniquely into pairwise edge-disjoint angular components of $L$ on $A$. We go on to study the structure of these decompositions somewhat more closely.

\medskip

\begin{proposition} \label{nas} Let $C$ be a non-angular connected component of $L$ on $A$. Suppose, for concreteness, that $C$ is a horizontal weave. Suppose also that $\CompY(C)$ is a simple weave. Then the angular components of $L$ on $A$ that $C$ consists of are zigzags of the same incline and the same length. \end{proposition} 

\medskip

The proof is straightforward.

Clearly, in the setting of Proposition \ref{nas} the number of angular components in the decomposition of the corresponding non-angular connected component of $L$ on $A$ can become arbitrarily large. The situation changes significantly with compound weaves.

\medskip

\begin{proposition} \label{nac} Let $C$ be a non-angular connected component of $L$ on $A$. Suppose, for concreteness, that $C$ is a horizontal weave. Suppose also that $\CompY(C)$ is a compound weave. Then $C$ consists of at most $p$ angular components of $L$ on $A$. Furthermore, the upper bound of $p$ is attained for all $n$ with $n \ge 2(p + q)$. Thus, in particular, it is attained whenever $n$ is sufficiently large. \end{proposition} 

\medskip

\begin{proof} Let $W = \CompY(C)$ and let $e_1$, $e_2$, \ldots, $e_k$ be all long edges of $W$ so that $e_i = w'_iw''_{i + 1}$ with $w'_i < w''_i$ for all $i$ and a path of short edges in $W$ connects $w''_i$ and $w'_{i + 1}$ for all $i$ with $1 \le i < k$.

Since $C$ is a non-angular horizontal weave, it follows as in Section \ref{ang} that $n \ge 2q + 1$.

Let $s$ be a remainder modulo $2p$.

For all $i$, we define $Z_{i, s}$ to be the steep zigzag of $L$ on $A$ induced by all cells $(jp + s, w'_i)$ of $A$ with $i + j$ odd and all cells $(jp + s, w''_i)$ of $A$ with $i + j$ even. Observe that $Z_{i, s}$ is nonempty for all $i$ since $n \ge 2q + 1$.

We also define $\Theta_s$ to be the leaper graph of $L$ on $A$ formed as the union of all $Z_{i, s}$ with $1 \le i \le k$ and all slight zigzags of $L$ on $A$ that contain a cell of some $Z_{i, s}$ with $1 \le i \le k$.

By Lemma \ref{rdo}, each steep zigzag $Z_{i, s}$ is in the same angular component as all slight zigzags that contain one of its cells.

We claim that, furthermore, $Z_{i, s}$ and $Z_{i + 1, s}$ are in the same angular component for all $i$ with $1 \le i < k$. Indeed, let $a = (x, w''_i)$ be a cell of $Z_{i, s}$. Then $b = (x, w'_{i + 1})$ is a cell of $Z_{i + 1, s}$. Since $w''_i \equiv w'_{i + 1} \pmod {2p}$ and $n \ge 2q + 1$, there exists a slight zigzag $Z$ of $L$ on $A$ that contains both of $a$ and $b$. But then both of $Z_{i, s}$ and $Z_{i + 1, s}$ must be in the same angular component as $Z$.

Thus $\Theta_s$ is an angular leaper graph of $L$ on $A$.

Conversely, it is straightforward to verify that each edge of $L$ on $A$ which makes an acute angle with an edge of $\Theta_s$ is also an edge of $\Theta_s$.

Or, in summary, $\Theta_s$ is an angular component of $L$ on $A$.

The union of the vertex sets of the $\Theta_s$ as $s$ ranges over all remainders modulo $2p$ coincides with the union of all rows of $A$ corresponding to a vertex of $W$. Consequently, all angular components in the decomposition of $C$ are of the form $\Theta_s$ for some $s$.

On the other hand, by Lemmas \ref{pp} and \ref{wb}, all $s$ such that $\Theta_s$ is a subgraph of $C$ are of the same parity. Therefore, $C$ does indeed consist of at most $p$ angular components of $L$ on $A$.

We are left to demonstrate that the upper bound of $p$ is attained for all $n$ with $n \ge 2(p + q)$.

We define an \emph{$s$-joint} to be a cell of $A$ of the form $(x, w'_1 + p)$ with $x_\text{Min} + q \le x \le x_\text{Max} - q$ and $x \equiv s \pmod {2p}$. Thus an $s$-joint belongs to two slight zigzags of $L$ on $A$ contained one each in $\Theta_{s - q}$ and $\Theta_{s + q}$. (For convenience, we let $\Theta_t$ denote the same leaper graph of $L$ on $A$ as $\Theta_{t \bmod 2p}$ for all integers $t$.)

When $n \ge 2(p + q)$, an $s$-joint exists on $A$ for all $s$. Consequently, $\Theta_{s - q}$ and $\Theta_{s + q}$ are in the same connected component of $L$ on $A$ for all $s$ as well. Since $p$ and $q$ are relatively prime, it follows that in fact all $\Theta_s$ with $s$ even are in the same connected component of $L$ on $A$ and so are all $\Theta_s$ with $s$ odd. \end{proof}

\medskip

We proceed to describe the rectangular boards which contain non-angular connected components of $L$.

The case when $p = 1$ is straightforward. Suppose, for concreteness, that $m \le n$. Then there exists a non-angular connected component of $L$ on $A$ if and only if $2 \le m \le q$ and $n \ge 2q + 1$.

Throughout the rest of this section, we consider the case when $p \ge 2$.

\medskip

\begin{proposition} \label{nab} Suppose that $p \ge 2$. Suppose also, for concreteness, that $m \le n$. (a) When $4p < 3q$, let $m_\measuredangle = p + q - (p \mmod [q \mmod 2p])$. (b) Otherwise, when $4p \ge 3q$, let $m_\measuredangle = 3q - 2p$. Then there exists a non-angular connected component of $L$ on $A$ if and only if $p < m < m_\measuredangle$ and $n \ge 2q + 1$. \end{proposition} 

\medskip

Conversely, on boards which do not satisfy the conditions of Proposition \ref{nab}, the angular components of $L$ coincide with the non-singleton connected components of $L$.

\medskip

\begin{proof} It suffices to show that there exists a non-angular connected component of $L$ on $A$ if and only if $n \ge 2q + 1$ and $\Pi_Y$ contains a non-singleton weave. The result would then follow by Proposition \ref{nws}.

(Necessity) Let $C$ be a non-angular connected component of $L$ on $A$. Since $m \le n$, our observations concluding Section \ref{ang} imply that $n \ge 2q + 1$ and $\CompY(C)$ is a non-singleton weave in $\Pi_Y$.

(Sufficiency) Let $n \ge 2q + 1$ and let $W$ be a non-singleton weave in $\Pi_Y$.

When $W$ is simple, since $n \ge 2q + 1$ there exist two slight zigzags $Z'$ and $Z''$ of $L$ on $A$ with a common cell such that both of their $y$-projections coincide with $W$. Then each one of $Z'$ and $Z''$ is an angular component of $L$ on $A$, and so the connected component of $L$ on $A$ which contains both of them is non-angular.

Otherwise, when $W$ is a compound weave, since $n \ge 2q + 1$ we obtain in the setting of Proposition \ref{nac} and its proof that there exists an $s$-joint on $A$ for some remainder $s$ modulo $2p$. Therefore, the two angular components $\Theta_{s - q}$ and $\Theta_{s + q}$ of $L$ on $A$ are in the same non-angular connected component of $L$ on $A$. \end{proof}